\documentclass[twoside,11pt,reqno]{amsart}
\usepackage{amsmath,amssymb,amscd,mathrsfs,epic,wasysym,latexsym,tikz,mathrsfs,cite,hyperref,enumerate}
\usepackage{pb-diagram}
\usepackage{tikz-cd}
\usepackage[all]{xy}
\usepackage{mathdots}
\usepackage[new]{old-arrows}

\makeatletter

\numberwithin{equation}{section}

\allowdisplaybreaks

\newtheorem{Proposition}[equation]{Proposition}
\newtheorem{Lemma}[equation]{Lemma}
\newtheorem{Theorem}[equation]{Theorem}
\newtheorem{Corollary}[equation]{Corollary}
\newtheorem{MainTheorem}{Theorem}

\theoremstyle{definition}  
\newtheorem{Definition}[equation]{Definition}
\newtheorem{Remark}[equation]{Remark}

\newtheorem{Example}[equation]{Example}

\let\<\langle
\let\>\rangle

\newcommand\Comment[2][\relax]{\space\par\medskip\noindent%
   \fbox{\begin{minipage}{\textwidth}\textbf{Comment\ifx\relax#1\else---#1\fi}\newline%
        #2\end{minipage}}\medskip
}

\def\bs{\text{\boldmath$s$}}

\def\br{\text{\boldmath$r$}}

\def\bA{\text{\boldmath$A$}}
\def\bB{\text{\boldmath$B$}}

\def\bz{\mathrm{\mathbf z}}

\def\bb{\mathrm{\mathbf b}}

\def\bla{\text{\boldmath$\lambda$}}
\def\bmu{\text{\boldmath$\mu$}}

\def\pmod#1{\text{ }(\text{\rm mod } #1)\,}

\newcommand{\Hom}{\operatorname{Hom}}
\newcommand{\Ext}{\operatorname{Ext}}

\newcommand{\ext}{\operatorname{ext}}
\newcommand{\End}{\operatorname{End}}

\newcommand{\id}{\operatorname{id}}
\def\sgn{\mathtt{sgn}}

\newcommand{\res}{\operatorname{res}}
\newcommand{\soc}{\operatorname{soc}}
\newcommand{\head}{\operatorname{head}}

\newcommand{\Z}{\mathbb{Z}}

\newcommand{\0}{{\bar 0}}
\renewcommand{\1}{{\bar 1}}
\def\eps{{\varepsilon}}
\def\phi{{\varphi}}

\newcommand{\zc}{{\textsf{c}}}

\newcommand{\zz}{{\textsf{z}}}
\newcommand{\zb}{{\textsf{b}}}

\newcommand{\ze}{{{\mathsf{e}}}}

\newcommand{\za}{{\mathsf{a}}}

\newcommand{\Ga}{\Gamma}
\newcommand{\la}{\lambda}
\newcommand{\La}{\Lambda}
\newcommand{\al}{\alpha}
\newcommand{\be}{\beta}

\def\Si{\mathfrak{S}}
\newcommand{\si}{\sigma}

\newcommand{\Om}{\Omega}

\newcommand{\de}{\delta}
\newcommand{\De}{\Delta}
\newcommand{\ka}{\kappa}

\newcommand{\Mull}{{\tt M}}

\newcommand{\quot}{\bold{quot}}

\def\id{\mathop{\mathrm {id}}\nolimits}
\renewcommand{\Im}{{\operatorname{Im}}}

\newcommand{\rad}{\operatorname{rad}}

\newcommand{\D}{{\mathscr D}}

\newcommand{\da}{{\downarrow}}
\newcommand{\ua}{{\uparrow}}

\renewcommand{\mod}{\bmod \,}

\newcommand{\Zig}{{\sf Z}}
\newcommand{\zH}{{\sf H}}

\def\Par{{\mathscr P}}

\def\k{\Bbbk}

\def\Parreg{{\mathscr P}^{\text{\rm $p$-reg}}}
\def\Parres{{\mathscr P}^{\text{\rm $p$-res}}}

\def\spa{\operatorname{span}}

\def\wt{{\operatorname{wt}}}

\def\onto{{\twoheadrightarrow}}
\def\into{{\hookrightarrow}}

\def\mod#1{#1\!\operatorname{-mod}}

\def\iso{\stackrel{\sim}{\longrightarrow}}
\def\B{{\mathsf B}}

\def\col{{\tt col}}
\def\row{{\tt row}}

\def\lan{\langle}
\def\ran{\rangle}

\def\Seq{\operatorname{Seq}}

\def\bx{\text{\boldmath$x$}}
\def\by{\text{\boldmath$y$}}

{\catcode`\|=\active
  \gdef\set#1{\mathinner{\lbrace\,{\mathcode`\|"8000%
  \let|\midvert #1}\,\rbrace}}
}
\def\midvert{\egroup\mid\bgroup}

\colorlet{darkgreen}{green!50!black}
\tikzset{dots/.style={very thick,loosely dotted},
         greendot/.style={fill,circle,color=darkgreen,inner sep=1.5pt,outer sep=0},
         blackdot/.style={fill,circle,color=black,inner sep=1.5pt,outer sep=0},
         graydot/.style={fill,circle,color=gray,inner sep=1.1pt,outer sep=0}
}
\def\greendot(#1,#2){\node[greendot] at(#1,#2){}}
\def\blackdot(#1,#2){\node[blackdot] at(#1,#2){}}
\def\graydot(#1,#2){\node[graydot] at(#1,#2){}}

\newenvironment{braid}{
  \begin{tikzpicture}[baseline=6mm,black,line width=1pt, scale=0.32,
                      draw/.append style={rounded corners},
                      every node/.append style={font=\fontsize{5}{5}\selectfont}]%
  }{\end{tikzpicture}
}

\def\Grid(#1,#2){
  \draw[very thin,gray,step=2mm] (0,0)grid(#1,#2);
  \draw[very thin,darkgreen,step=10mm] (0,0)grid(#1,#2);
}

\newcommand\Tableau[2][\relax]{
  \begin{tikzpicture}[scale=0.5,draw/.append style={thick,black}]
    \ifx\relax#1\relax%
    \else 
      \foreach\box in {#1} { \filldraw[blue!30]\box+(-.5,-.5)rectangle++(.5,.5); }
    \fi
    \newcount\row\newcount\col
    \row=0
    \foreach \Row in {#2} {
       \col=1
       \foreach\k in \Row {
          \draw(\the\col,\the\row)+(-.5,-.5)rectangle++(.5,.5);
          \draw(\the\col,\the\row)node{\k};
          \global\advance\col by 1
       }
       \global\advance\row by -1
    }
  \end{tikzpicture}
}

\newcommand\YoungDiagram[2][\relax]{
  \begin{tikzpicture}[scale=0.5,draw/.append style={thick,black}]
    \ifx\relax#1\relax%
    \else 
    \foreach\box in {#1} {
      \filldraw[blue!30]\box rectangle ++(1,1);
    }
    \fi
    \newcount\row
    \row=0
    \foreach \col in {#2} {
       \draw(1,\the\row)grid ++(\col,1);
       \global\advance\row by -1
    }
  \end{tikzpicture}
}

\newdimen\hoogte    \hoogte=12pt    
\newdimen\breedte   \breedte=14pt  
\newdimen\dikte     \dikte=0.5pt 

\newenvironment{Young}{\begingroup
       \def\vr{\vrule height0.89\hoogte width\dikte depth 0.2\hoogte}
       \def\fbox##1{\vbox{\offinterlineskip
                    \hrule height\dikte
                    \hbox to \breedte{\vr\hfill##1\hfill\vr}
                    \hrule height\dikte}}
       \vbox\bgroup \offinterlineskip \tabskip=-\dikte \lineskip=-\dikte
            \halign\bgroup &\fbox{##\unskip}\unskip  \crcr }
       {\egroup\egroup\endgroup}
\def\Youngdiagram#1{\relax\ifmmode\vcenter{\,\begin{Young}#1\end{Young}\,}\else%
              $\vcenter{\,\begin{Young}#1\end{Young}\,}$\fi}

\newcommand{\cont}{{\operatorname{cont}}}
\newcommand{\te}{{\tilde e}}
\newcommand{\tf}{{\tilde f}}
\newcommand{\he}{{\hat e}}
\newcommand{\hf}{{\hat f}}

\newcommand{\laM}{{\la^{\Mull}}}
\newcommand{\muM}{{\mu^{\Mull}}}
\newcommand{\mods}{\textrm{-mod}}

\newcommand{\beps}{\eps'}
\newcommand{\bphi}{\phi'}

\newcommand{\reg}{{\tt R}}
\newcommand{\Rem}{{\operatorname{Rem}}}
\newcommand{\Add}{{\operatorname{Add}}}

\begin{document}
\voffset = -30 pt

\title[Self-extensions over symmetric groups]{{\bf On self-extensions of irreducible modules over symmetric groups}}

\author{\sc Haralampos Geranios}
\address
{Department of Mathematics\\ University of York\\ York YO10 5DD, U.K.} 
\email{haralampos.geranios@york.ac.uk}

\author{\sc Alexander Kleshchev}
\address{Department of Mathematics\\ University of Oregon\\Eugene\\ OR 97403, USA}
\email{klesh@uoregon.edu}

\author{\sc Lucia Morotti}
\address
{Institut f\"{u}r Algebra, Zahlentheorie und Diskrete Mathematik\\ Leibniz Universit\"{a}t Hannover\\ 30167 Hannover\\ Germany} 
\email{morotti@math.uni-hannover.de}



\thanks{The first author gratefully acknowledges the support of The Royal Society through a University Research Fellowship. 
The second author was supported by the NSF grant DMS-1700905.  
The third author was supported by the DFG grants MO 3377/1-1 and MO 3377/1-2.}


\begin{abstract}
A conjecture going back to the eighties claims that there are no non-trivial self-extensions of irreducible modules over symmetric groups if the characteristic of the ground field is not equal to $2$. We obtain some partial positive results on this conjecture. 
\end{abstract}


\maketitle

\section{Introduction}\label{SIntro}
Let $\k$ be a field of characteristic $p>0$ and $\Si_n$ be the symmetric group on $n$ letters. In this paper we are concerned with the following conjecture:

\vspace{2mm}
\noindent
{\bf Self-extensions Conjecture for Symmetric Groups.}
{\em Let $p>2$ and $D$ be an irreducible $\k\Si_n$-module. Then $\Ext^1_{\Si_n}(D,D)=0$. 
}

\vspace{2mm}
This folklore conjecture, sometimes referred to as Kleshchev-Martin's conjecture, goes back to the late eighties. As even the case of the trivial module $D=\k_{\Si_n}$ shows, the assumption $p>2$ is necessary. The conjecture seems to be wide open. 

As in \cite[\S11]{JamesBook}, the irreducible $\k\Si_n$-modules are $\{D^\la\mid\la\in\Parreg_n\}$ where $\Parreg_n$ denotes the set of $p$-regular partitions of $n$. We also have the Specht modules $\{S^\la\mid \la\in\Par_n\}$ where $\Par_n$ is the set of all partitions of $n$, see \cite[\S4]{JamesBook}. We denote by $h(\la)$ the number of non-zero parts of a partition $\la$, and by $\unrhd$ the usual dominance order on $\Par_n$, see \cite[3.2]{JamesBook}. 

In \cite[Theorem 2.10]{KS} it is proved that under the assumptions 
$p>2$, $\la\,{\hspace{.1mm}\not}{\rhd}\,\mu$ and $h(\la),h(\mu)\leq p-1$, we have $\Ext^1_{\Si_n}(D^\la,D^\mu)\cong \Hom_{\Si_n}(\rad S^\la, D^\mu)$. In view of \cite[12.2]{JamesBook}, this immediately implies:

\begin{Proposition} \label{PKS} 
Let $p>2$ and $\la\in\Parreg_n$. If $h(\la)\leq p-1$ then $\Ext^1_{\Si_n}(D^\la,D^\la)=0$. 
\end{Proposition}

If $D^\la$ is isomorphic to a Specht module we immediately get from \cite[Theorem 3.3(c)]{KN}:

\begin{Proposition} \label{PKN} 
Let $p>3$ and $\la\in\Parreg_n$. If $D^\la\cong S^\nu$ for some $\nu\in\Par_n$, then $\Ext^1_{\Si_n}(D^\la,D^\la)=0$. 
\end{Proposition}

These are the only general results about self-extensions of irreducible modules over symmetric groups that we are aware of. In this paper we obtain several new positive results. 

\begin{MainTheorem}\label{TA}
Let $p>2$ and $D^\la$ be in a RoCK block. Then $\Ext^1_{\Si_n}(D^\la,D^\la)=0$. 
\end{MainTheorem}

Theorem~\ref{TA} is proved in Section~\ref{SRoCK} using a Morita equivalence, established in \cite{EK2}, between weight~$d$ RoCK blocks of symmetric groups and zigzag Schur algebras $T^\Zig(m,d)$ with $m\geq d$ which were defined by Turner \cite{T}, see also \cite{EK1}. We establish in Corollary~\ref{CSelExtT} that $\Ext^1_{T^\Zig(m,d)}(L,L)=0$ for any irreducible $T^\Zig(m,d)$-module $L$. This implies Theorem~\ref{TA}. Our results on extensions in RoCK blocks are actually stronger, and we refer the reader to  Corollary~\ref{CZeroShift}, Theorem~\ref{TExtT} and Remark~\ref{R100420} for more details on that.

By \cite{CR}, every block of a symmetric group is derived equivalent to a RoCK block, and so one might hope to extend Theorem~\ref{TA} to an arbitrary block using Chuang-Rouquier's perverse equivalences. We were unable to do that, so we had to resort to a less powerful approach employing translation functors. 
We  review the fundamental properties of the translation functors $e_i^{(r)}, f_i^{(r)}$ in \S\ref{subsec:tran.fun}, and in Section~\ref{sec:bra.ru} we establish some of their new properties. This allows us to prove in Section~\ref{SSmallWeight} at least the following:

\begin{MainTheorem}\label{TSmallWeight}
Let $p>2$ and $D^\la$ be in a block of weight $\leq 7$. Then $\Ext^1_{\Si_n}(D^\la,D^\la)=0$. 
\end{MainTheorem}

Using translation functors and some information about Specht modules, in Section~\ref{SSmallHeight} we also improve on  Proposition~\ref{PKS} above:

\begin{MainTheorem}\label{TSmallHeight}
Let $p>2$ and $\la\in\Parreg_n$. If $h(\la)\leq p+2$ then $\Ext^1_{\Si_n}(D^\la,D^\la)=0$.  
\end{MainTheorem}

The following result, proved in section~\ref{sec:res.irr.Sp}, 
verifies the Self-extensions Conjecture for some additional cases:

\begin{MainTheorem}\label{TB}
Let $p>2$, $\la\in\Parreg_n$ and $i\in \Z/p\Z$. If  $e_i^{(\eps_i(\la))}D^\la$ is isomorphic to an (irreducible) Specht module then $\Ext^1_{\Si_n}(D^\la,D^\la)=0$. 
 \end{MainTheorem}

If $D^\la$ itself is isomorphic to a Specht module, it is easy to see that it satisfies the assumptions of Theorem~\ref{TB}, so in particular the theorem generalizes Proposition~\ref{PKN}.   
We also note that the assumption that $e_i^{(\eps_i(\la))}D^\la$ is an (irreducible) Specht module in Theorem~\ref{TB} is equivalent to the assumption that $f_i^{(\phi_i(\la))}D^\la$ is an (irreducible) Specht module, see \S\ref{subsec:Sp.result}. We refer to Example~\ref{exm:res.irr.Sp} for a concrete example of an application of Theorem~\ref{TB}.

The new results on translation functors obtained in Section~\ref{sec:bra.ru} 
might be of independent interest, so we cite some of them here. The main point is that the divided power $i$-restriction functor $e_i^{(r)}$, when applied to an irreducible module $D^\la$, has a simple socle $D^{\te_i^r\la}$, and it is crucially important to know that the quotient $(e_i^{(r)}D^\la)/D^{\te_i^r\la}$ has no $D^{\te_i^r\la}$ in the socle (and similarly for the $i$-induction functor $f_i^{(r)}$). Unfortunately, we cannot prove this in general. But at least we establish the following: 
 
\begin{MainTheorem}\label{TC}
Let $\la\in\Parreg_n$, $i\in\Z/p\Z$, $B_1,\ldots,B_{\phi_i(\la)}$ be the $i$-conormal nodes of $\la$ counted from top to bottom, and $A_1,\ldots,A_{\eps_i(\la)}$ be the $i$-normal nodes of $\la$ counted from bottom to top.
\begin{enumerate}
\item[{\rm (i)}] If $0\leq r\leq \phi_i(\la)$ and $D^{\tf_i^r\la}\subseteq (f_i^{(r)}D^\la)/D^{\tf_i^r\la}$ then $0< r<\phi_i(\la)$ and the partition $\la^{B_{r+1}}$, obtained by adding $B_{r+1}$ to $\la$, is not $p$-regular.
\item[{\rm (ii)}] If $0\leq r\leq \eps_i(\la)$ and $D^{\te_i^r\la}\subseteq (e_i^{(r)}D^\la)/D^{\te_i^r\la}$ then $0< r<\eps_i(\la)$ and the partition $\la_{A_{r+1}}$, obtained by removing $A_{r+1}$ from $\la$, is not $p$-regular.
\end{enumerate}
\end{MainTheorem}

Some consequences of Theorem~\ref{TC} for self-extensions are obtained in \S\ref{SSConsSelfExt}.

\section{Preliminaries}\label{sec:prel}
We review some notions related to representation theory of the symmetric group $\Si_n$ referring the reader to \cite{JamesBook,JK} for details. 
We stick with the notation already introduced in Section~\ref{SIntro}. In particular, we work over the ground field  $\k$ of characteristic $p>0$. Throughout the paper we assume that $\k$ is algebraically closed (for symmetric groups this does not reduce any generality). 

\subsection{Generalities on representations}
\label{subsec:Irred}
Let $G$ be a finite group. We denote by $\mod{\k G}$ the category of finite dimensional $\k G$-modules. Let  $V,V_1,\dots,V_s\in\mod{\k G}$. We write $V\sim V_1\mid\ldots\mid V_s$ to indicate that $V$ has a filtration as a $\k G$-module with factors $V_1,\ldots,V_s$, listed from bottom to top. 
For an irreducible $\k G$-module $D$ we write $[V:D]$ for the multiplicity of $D$ as a composition factor of $V$. We denote by $\k_G$ the trivial $\k G$-module. 

\begin{Lemma} \label{LUniquenessSubmodule} 
Let $G$ be a finite group, $W\in\mod{\k G}$, and $V,V'\subseteq W$ be submodules with $V\cong V'$. If $\soc W\cong D$ is irreducible and $\dim \End_G(V)=[V:D]$ then $V= V'$. 
\end{Lemma}
\begin{proof} Let $J$ be the injective hull of  $D$. 
If $V\not=V'$ then 
\begin{align*}
\dim\Hom_{G}(V,W)>\dim\End_{G}(V)=[V:D]=\dim\Hom_{G}(V,J)\geq \dim\Hom_{G}(V,W),
\end{align*}
giving a contradiction.
\end{proof}

\begin{Lemma}\label{L100221}
Let $G$ be a finite group and $V,U_1,\dots,U_m\in\mod{\k G}$. If\, $\soc V$ is simple and there is an injective homomorphism $f:V\to  U_1\oplus\dots\oplus U_m$, then there is an injective homomorphism $V\to U_j$ for some $1\leq j\leq m$.
\end{Lemma}

\begin{proof}
There exist homomorphisms  $\pi_j:U_1\oplus\dots \oplus U_m\to U_j,\ \iota_j:U_j\to U_1\oplus\dots \oplus U_m$ with $\sum_j\iota_j\circ \pi_j=\id$.  Now $f=\sum_j\iota_j\circ \pi_j\circ f$ is injective and $\soc V$ is simple, so $\iota_j\circ \pi_j\circ f |_{\soc V}\neq 0$ for some $j$. Then
$\pi_j\circ f:V\to U_j$ is injective.
\end{proof}

\begin{Lemma}\label{L110221}
Let $G$ be a finite group and $M,N,V,W\in\mod{\k G}$. Suppose:
\begin{enumerate}
\item[{\rm (a)}] $M\cong V^{\oplus k}$ and $N\cong W^{\oplus l}$ for some $k\leq l$;
\item[{\rm (b)}] $V$ is a submodule of $W$ and $M$ is a submodule of $N$;
\item[{\rm (c)}] $D:=\soc W$ is simple and $\dim\End_G(W)=[W:D]$.
\end{enumerate}
Then $N/M\cong (W/V)^{\oplus k}\oplus W^{\oplus l-k}$.
\end{Lemma}

\begin{proof}
Let $J$ be the injective hull of $D$. Since $\soc W=D$ by (c), we have $W\subseteq J$. So, using (a), we get for any $X\in\mod{\k G}$:
\begin{equation}\label{E250521}
\dim\Hom_G(X,N)=l\dim\Hom_G(X,W)\leq  l\dim\Hom_G(X,J)= l[X:D].
\end{equation}

The short exact sequence $0\to V\to W\to W/V\to 0$ yields an exact sequence 
\begin{equation}\label{E250521_2}
0\longrightarrow \Hom_G(W/V,N)\longrightarrow\Hom_G(W,N)\stackrel{\phi}{\longrightarrow}\Hom_G(V,N).
\end{equation}
Using (c), (\ref{E250521_2}), and (\ref{E250521}) with $X=V$ and $X=W/N$, we now get 
\begin{align*}
l[W:D]&=l\dim\End_G(W)=\dim\Hom_G(W,N)
\\
&\leq\dim\Hom_G(V,N)+\dim\Hom_G(W/V,N)
\\
&\leq l[V:D]+l[W/V:D]=l[W:D].
\end{align*}
Hence 
\[\dim\Hom_G(W,N)=\dim\Hom_G(V,N)+\dim\Hom_G(W/V,N).\]
It follows that the map $\phi$ in (\ref{E250521_2}) is surjective. In other words, 
every homomorphism $f\in \Hom_G(V,N)$ extends to a homomorphism $\hat f\in \Hom_G(W,N)$.

By (a), we can write $M= M_1\oplus\ldots\oplus M_k$ and $N= N_1\oplus\ldots\oplus N_l$ for submodules 
$M_1,\dots,M_k$ of $M$ isomorphic to $V$ and submodules $N_1,\dots,N_l$ of $N$ isomorphic to $W$. It suffices to show  that we can also write $N= N_1'\oplus\ldots\oplus N_l'$ for submodules $N_1',\dots,N_l'$ of $N$ isomorphic to $W$ and such that $M_i\subseteq N_i'$ for $i=1,\dots,k$. We may assume that $V\neq 0$. We have the embeddings 
$$\iota_i:V\iso M_i\,\longhookrightarrow\, M\,\longhookrightarrow\, N\qquad(1\leq i\leq k).$$ 
By the previous paragraph, these can be extended to homomorphisms $\hat\iota_i:W\to N$ which are necessarily injective since $\soc V=\soc W$ is simple. Let $N_i'$ be the image of $\hat \iota_i$ for $i=1,\dots,k$. Using the simple socles and the fact that the sum $\sum_{i=1}^k M_i$ is direct we deduce that the sum $\sum_{i=1}^k N_i'$ is also direct. 
Now, up to relabeling the modules $N_1,\dots,N_l$, we may assume that
$$\soc(N_1'\oplus \ldots\oplus N_k')\cap\soc(N_{k+1}\oplus\ldots\oplus N_l)=0.$$
Then 
\[(N_1'\oplus \ldots\oplus N_k')\cap(N_{k+1}\oplus\ldots\oplus N_l)=0.\]
By dimensions, we conclude 
\[N=N_1'\oplus \ldots\oplus N_k'\oplus N_{k+1}\oplus\ldots\oplus N_l,\]
and so we can take $N_j':=N_j$ for $j=k+1,\dots,l$. 
\end{proof}

\begin{Lemma} \label{LSubmoduleStructure} 
Let $G$ be a finite group, $D$ be an irreducible $\k G$-module,  and suppose that $V\in\mod{\k G}$ has a filtration
$0= V_0\subseteq V_1\subseteq\dots\subseteq V_b=V$ such that $\soc V_a\cong \head V_a\cong D$ and $[V_a:D]=a$ for  for all $a=1,\dots b$. 

\begin{enumerate}
\item[{\rm (i)}] If $Z\subseteq V$ is a submodule with $\head Z\cong D^{\oplus m}$ for some $m\in\Z_{\geq 0}$, then $Z=V_a$ for some $a$. 
\item[{\rm (ii)}] If $X\subseteq Y\subseteq V$ are submodules such that $[X:D]=a-1$ and $[Y:D]=a$ for some $a\in\Z_{>0}$, then  $V_a\not\subseteq X$, $V_a\subseteq Y$, and  
$V_a/(V_a\cap X)\into Y/ X$. In particular, if $D:=\soc (Y/X)$ is  simple, then $[V_a:D]\neq 0$. 
\end{enumerate}
\end{Lemma}
\begin{proof} 
(i) Let $a$ be minimal with $Z\subseteq V_a$. Then either $Z=V_a$ or $Z\subseteq \rad V_a$, whence $Z\subseteq V_{a-1}$ since $[(\rad V_a)/V_{a-1}:D]=0$, giving a contradiction. 

(ii) We only need to prove that $V_a\subseteq Y$. Note that $Y$ has a submodule $Z$ such that $[Z:D]=a$ and $\head Z\cong D^{\oplus m}$ for some $m\in\Z_{\geq 0}$. By (i), we have $Z=V_a$.
\end{proof}

\subsection{Partitions and abaci}\label{subsec:comb}
Let $\la=(\la_1,\la_2,\dots)\in\Par_n$. We denote $|\la|:=n$. The transpose partition is denoted $\la'$. Recall that $\la$ is called 
{\em $p$-restricted} if $\la_k-\la_{k+1}<p$ for all $k$. Then $\la$ is $p$-regular if and only if $\la'$ is $p$-restricted. We denote the sets of $p$-regular (resp. $p$-restricted) partitions of $n$ by $\Parreg_n$ (resp. $\Parres_n$). A partition $\la$ that is not $p$-regular is  called \emph{$p$-singular}. For every $\la\in \Par_n$, James \cite{J2} defines its {\em regularization} $\la^\reg\in\Parreg_n$. Note that $\la^\reg\unrhd\la$, and $\la^\reg=\la$ if and only if $\la\in\Parreg_n$. We refer the reader to \cite{Mull,FK, BO, K2} for the Mullineux involution 
$$\Parreg_n\to \Parreg_n,\ \la\mapsto \la^\Mull.$$ 
We denote by $\varnothing$ the trivial partition of $0$, thus $\Par_0=\{\varnothing\}$. 

We identify $\la$ with its  \emph{Young diagram} $\{(k,l)\in \Z_{>0}\times \Z_{>0}\mid l \leq \la_k\}$. The elements of $\Z_{>0}\times \Z_{>0}$ are called {\em nodes}. We set $I:=\Z/p\Z$, identified with $\{0,1,\dots,p-1\}$.
Let $A=(k,l)$ be a node. The \emph{residue} of $A$ is $\res A:=l-k\pmod{p}\in I$. 
For $k',l'\in\Z$ such that $k+k',l+l'>0$ we have the node $A+(k',l'):=(k+k',l+l')$. 

Consider the free $\Z$-module $Q:= \bigoplus_{i\in I}\Z\cdot \al_i$ with basis $\{\al_i\mid i\in I\}$. We have the subsets $Q^+:=\{\sum_{i\in I} c_i\al_i\in Q\mid c_i\geq 0\ \text{for all $i\in I$}\}$ and 
$Q^+_n:=\{\sum_{i\in I} c_i\al_i\in Q^+\mid \sum_{i\in I}c_i=n\}$. 
 The {\em residue content} of a partition $\la\in\Par_n$ is $\cont(\la):=\sum_{i\in I}a_i\al_i\in Q^+_n$, where $a_i$ is the number of nodes of $\la$ of residue $i$.

We assume familiarity with the \emph{abacus} notation for  partitions, see \cite[\S 2.7]{JK}. Recall that positions on the abacus are labeled with non-negative integers, so that for $i\in I$, the positions $\{i+pa\mid a\in\Z_{\geq 0}\}$ form the runner $i$ of the abacus. We denote an abacus display for $\la$ by $\Ga(\la)$. Recall that $\Ga(\la)$ is not unique and depends on the number of beads chosen, so we will need to make sure that the number of beads is agreed upon. In particular, we will always make sure that position $0$ in $\Ga(\la)$ is occupied, and if some operation with abaci creates an abacus with position $0$ unoccupied we will simply pass to the equivalent abacus with $p$ more beads.  
A position $k>0$ in $\Ga(\la)$ is {\em removable} (resp. {\em addable}) if it is occupied (resp. unoccupied) and position $k-1$ is unoccupied (resp. occupied).

Note that $\la$ is $p$-regular (resp. $p$-restricted) if and only if there is no unoccupied (resp. occupied) position $r$ in $\Ga(\la)$ such that the positions $r+1,\ldots,r+p$ (resp. $r-1,\ldots,r-p$) are occupied (resp. unoccupied). By replacing each bead in $\Ga(\la)$ with an empty space and vice versa, and then rotating the abacus through $180^\circ$, we obtain the abacus display $\Ga(\la)'$ for $\la'$, so 
\begin{equation}\label{EAbTr}
\Ga(\la)' =\Ga(\la').
\end{equation} 

We refer the reader to \cite[\S2.7]{JK} for the notions of  the \emph{core} and \emph{weight} of a partition. The weight of $\la$ is denoted $\wt(\la)$. 
Let $\rho\in\Par_r$ be a core, $d\in\Z_{\geq 0}$, and set $n:=r+pd$. Denote by $\Par_{\rho,d}$ the set of all partitions of $n$ with core $\rho$ (and weight $d$), and denote $\Parres_{\rho,d}:=\Par_{\rho,d}\cap\Parres_n$. By \cite[Theorem 2.7.41]{JK}, two partitions $\la,\mu\in \Par_n$ have the same core if and only if $\cont(\la)=\cont(\mu)$.

We will also use the notion of the \emph{quotient}  of a partition $\la\in\Par_{\rho,d}$ denoted $\quot(\la):=(\la^{(0)},\dots,\la^{(p-1)})$. This is a multipartition of $d$, with each $\la^{(i)}$ being the partition corresponding to the moves made on the runner $i$ to go from $\Ga(\rho)$ to $\Ga(\la)$, see \cite[2.7.29]{JK}. Note that, unlike in \cite{JK}, we do not insists on using an abacus $\Ga(\la)$ with a multiple of $p$ beads, and so $\quot(\la)=(\la^{(0)},\dots,\la^{(p-1)})$ depends on $\Ga(\la)$ and is only defined in general up to a cyclic permutation of $\la^{(0)},\dots,\la^{(p-1)}$. This is in agreement with \cite{F1}. For $i\in I$, we denote $\wt_i(\Ga):=|\la^{(i)}|$. We have $\sum_{i\in I}\wt_i(\la)=\wt(\la)$. We also denote by $\Ga_j$ the bead configuration on the runner $j$ of $\Ga$ and by $r_j(\Ga)$ the number of beads on $\Ga_j$.

\subsection{Removable and addable nodes}\label{SSRA}
Let $\la\in \Par_n$. A node $A\in \la$ (resp $B\not\in \la$) is called \emph{removable} (resp. {\em addable}) for $\la$ if  
$\la_A:=\la\setminus\{A\}$ (resp. $\la^B:=\la\cup\{B\}$) is a diagram of a partition. Fix $i\in I$. 
A removable (resp. addable) node is called {\em $i$-removable} (resp. {\em $i$-addable}) if it has residue $i$. 
If $A_1,A_2,\ldots, A_k$ are removable (resp. addable) nodes of $\lambda$ then we denote $\la_{A_1,\dots,A_k}:=\la\setminus\{A_1,\dots,A_k\}$ (resp. $\la^{A_1,\dots,A_k}:=\la\cup\{A_1,\dots,A_k\}$). 

For every removable node $A$ of $\la$, there exists a unique removable position $k$ in $\Ga(\la)$ such that $\Ga(\la_A)$ is obtained by moving a bead from position $k$ to position $k-1$. For every addable node $B$ for $\la$, there exists a unique addable position $k$ in $\Ga(\la)$ such that $\Ga(\la^B)$ is obtained by moving a bead from position $k-1$ to position $k$. 
Removable (resp. addable) nodes $A$ and $A'$ have the same residue if and only if the corresponding removable (resp. addable) positions $k$ and $k'$ are on the same runner. We denote by $i_j(\Ga)$ the residue of the removable/addable nodes of $\la$ corresponding to removable/addable positions on runner $j$. 

Labelling the $i$-addable
nodes of $\la$ by $+$ and the $i$-removable nodes of $\la$ by $-$, the {\em $i$-signature} of 
$\la$ is the sequence of pluses and minuses obtained by going along the 
rim of the Young diagram from bottom left to top right and reading off
all the signs.
The {\em reduced $i$-signature} of $\la$ is obtained 
from the $i$-signature
by successively erasing all neighboring 
pairs of the form $-+$. 
The nodes corresponding to  $-$'s (resp. $+$'s) in the reduced $i$-signature are
called {\em $i$-normal} (resp. {\em $i$-conormal}) for $\la$.
The leftmost $i$-normal (resp. rightmost $i$-conormal) node is called {\em $i$-good} (resp. {\em $i$-cogood}) for $\la$. 
We write
\begin{align*}
  \label{count.nodes'}
\eps_i(\la):=\sharp\{i\text{{-normal nodes of}}\ \la\}\quad\text{and}\quad
 \phi_i(\la):=\sharp \{i\text{{-conormal nodes of}}\ \la\}.
 \end{align*}
 

\begin{Lemma}\label{L110320}
Let $\la$ be a $p$-regular partition, $\Gamma=\Ga(\la)$ and $1\leq j\leq p-1$. If $r_j(\Gamma)-r_{j-1}(\Gamma)\geq \wt_{j-1}(\Gamma)+\wt_j(\Gamma)$, 
then 
$\eps_{i_j(\Ga)}(\la)=r_j(\Gamma)-r_{j-1}(\Gamma)$ 
and $\phi_{i_j(\Ga)}(\la)=0$.
\end{Lemma}

\begin{proof}
If a position $(j-1)+pa$ on runner $j-1$ is occupied then $a<r_{j-1}(\Gamma)+\wt_{j-1}(\Ga)$. So, by assumption, $a<r_j(\Gamma)-\wt_j(\Gamma)$ and the position $j+pa$ on runner $j$ is also occupied. 
The result follows. 
\end{proof}

Let $\la\in \Parreg_n$ and $i\in I$. Let $A_1,A_2,\ldots,A_{\eps_i(\la)}$ (resp. $B_1,B_2,\ldots,B_{\phi_i(\la)}$) be the $i$-normal (resp. $i$-conormal) nodes for $\la$, labelled from bottom to top (resp. from  top to bottom). We set 
\begin{equation*}
\te^r_i\la:=\la_{A_1,\ldots,A_r},\quad \tf^r_i\la:=\la^{B_1,\ldots,B_r},
\end{equation*}
where $\te^r_i\la$ (resp. $\tf^r_i\la$) is interpreted as $0$ if $r>\eps_i(\la)$ (resp. $r>\phi_i(\la)$). It is well known that the partitions $\te^r_i\la$ and $\tf^r_i\la$ are $p$-regular. Moreover, for $r\leq \phi_i(\la)$ (resp. $r\leq \eps_i(\la)$) we have $\te_i^r\tf_i^r\la=\la$ (resp. $\tf_i^r\te_i^r\la=\la$).

\begin{Lemma}\label{Lwt}
Let $\la\in\Parreg_n$, $i\in I$, and $\mu=\tilde f_i^r\la$ for some $0\leq r\leq\phi_i(\la)$. Then $\wt(\mu)=\wt(\la)+r(\phi_i(\la)-\eps_i(\la)-r)$. In particular, if $\eps_i(\la)=0$ then $\wt(\mu)=\wt(\la)+r(\phi_i(\la)-r)$, and $\wt(\mu)\geq\wt(\la)$, with equality holding if and only if $r=0$ or $r=\phi_i(\la)$.
\end{Lemma}

\begin{proof}
Let $\cont(\la):=\sum_{j\in I}a_j\al_j$. Note that $\cont(\mu)=\sum_{j\in I}(a_j+r\de_{i,j})\al_j$. In view of \cite[Lemmas 11.1.4, 11.1.5]{K3} we have that
\begin{align*}
\wt(\mu)-\wt(\la)&=r\de_{i,0}-\sum_j((a_j+r\de_{i,j})^2-a_j^2)+\sum_j((a_j+r\de_{i,j})(a_{j+1}+r\de_{i,j+1})-a_ja_{j+1})\\
&=r(\de_{i,0}-2a_i+a_{i-1}+a_{i+1}-r).
\end{align*}
The result then follows from $\phi_i(\la)-\eps_i(\la)=
\de_{i,0}-2a_i+a_{i-1}+a_{i+1}$ (this comes from \cite[Lemma 8.5.8]{K3} and can also be seen by induction, starting with the empty partition and considering when the nodes at the right or below $A$ are addable in $\nu^A$).
\end{proof}

There is also a purely combinatorial proof of Lemma~\ref{Lwt} by comparing 
$\Ga(\la)$ and $\Ga(\mu)$. 

\begin{Lemma}\label{L110320_2}
Let $\la\in\Parreg_{n}$, and $i\in I$ be such that $\eps_i(\la)>0$, $\phi_i(\la)>0$ and $\la^B_A\not\in\Parreg_n$ for the $i$-good node $A$ and the $i$-cogood node $B$ for $\la$. Suppose that $a$ is the removable position on $\Ga(\la)$ corresponding to $A$ and $b$ is the addable position on $\Ga(\la)$ corresponding to $B$. Then $a=b+p$ and the positions $c$ satisfying $b<c<a-1$ are all occupied. 
\end{Lemma}

\begin{proof}
Since $\la_A$ and $\la^B$ are $p$-regular, it follows that $A=B+(1-p,1)$. This is equivalent to the required property of  $\Ga(\la)$. 
\end{proof}

The following two lemmas can be easily checked and are left as an exercise.

\begin{Lemma}
\label{lem:re.nod}
Let $\la\in\Parreg_n$, $i\in I$ and $A_1,\ldots,A_{\eps_i(\la)}$ be the $i$-normal nodes of $\la$ labeled from bottom to top. For  $1\leq r\leq\eps_i(\la)$, the following are equivalent:
\begin{itemize}
\item[(i)] $\la_{A_r}\not\in\Parreg_{n-1}$;
\item[(ii)] $r\geq 2$ and for some $j\leq r-2$ we have $\la_{A_1,\ldots,A_j,A_r}\not\in\Parreg_{n-j-1}$;
\item[(iii)] $r\geq 2$ and $A_r=A_{r-1}+(1-p,1)$.
\end{itemize}
\end{Lemma}

\begin{Lemma}
\label{lem:ad.nod}
Let $\la\in\Parreg_n$, $i\in I$ and $B_1,\ldots,B_{\phi_i(\la)}$ be the $i$-conormal nodes of $\la$ labeled from top to bottom. Fix $1\leq r\leq\phi_i(\la)$. The following are equivalent:
\begin{itemize}
\item[(i)] $\la^{B_r}\not\in\Parreg_{n+1}$;
\item[(ii)] $r\geq 2$ and for some $j\leq r-2$ we have $\la^{B_1,\ldots,B_j,B_r}\not\in\Parreg_{n+j+1}$;
\item[(iii)] $r\geq 2$ and $B_r=B_{r-1}+(p-1,-1)$.
\end{itemize}
\end{Lemma}

Let $\la\in \Par_n$ and $i\in I$. Define
\begin{align*}
 \label{count.nodes}
 \beps_i(\la):=\sharp\{i\text{{-removable nodes of}}\ \la\}\quad\text{and}\quad
 \bphi_i(\la):=\sharp \{i\text{{-addable nodes of}}\ \la\}. 
 \end{align*}
Let $A_1,A_2,\ldots,A_{\beps_i(\la)}$ (resp. $B_1,B_2,\ldots,B_{\bphi_i(\la)}$) be the $i$-removable (resp. $i$-addable) nodes for $\la$, labelled from bottom to top (resp. from  top to bottom). We set 
\begin{equation*}
\he^r_i\la:=\la_{A_1,\ldots,A_r},\quad \hf^r_i\la:=\la^{B_1,\ldots,B_r},
\end{equation*}
where $\he^r_i\la$ (resp. $\hf^r_i\la$) is interpreted as $0$ if $r>\beps_i(\la)$ (resp. $r>\bphi_i(\la)$). 

The following is easy to see:

\begin{Lemma}
\label{lem:core2}
Let $\la$ be a core, $i\in I$,  and $B_1,\ldots,B_{\phi_i'(\la)}$ be the $i$-addable nodes of $\la$ labeled from top to bottom. If some $\la^{B_r}$ is not $p$-restricted then 
$\la^{B_s}$ is not $p$-restricted for all $s=1,\dots, r$ and $B_1=(1,\la_1+1)$ is the top addable node of $\la$.
\end{Lemma}

\subsection{Representations of symmetric groups}\label{SSSG}
In addition to the notation introduced in Section~\ref{SIntro}, 
we denote by $\sgn$ the sign representation of $\Si_n$. Let $\la\in\Parreg_n$. 
By \cite{FK,BO,K2}, we have $D^\la\otimes \sgn\cong D^{\laM}.$  Moreover $\eps_i(\la)=\eps_{-i}(\laM)$ and $\phi_i(\la)=\phi_{-i}(\laM)$ for all $i\in I$. Recall that $(D^\la)^*\cong D^\la$. Passing to duals and tensoring with $\sgn$, we deduce for all $\la,\mu\in\Parreg_n$ and $k\geq0$:
\begin{equation}\label{ESignExt}
\Ext^k_{\Si_n}(D^\mu,D^\la)\cong \Ext^k_{\Si_n}(D^\la,D^\mu)\cong\Ext^k_{\Si_n}(D^\laM,D^\muM).
\end{equation}

Let $\rho\in\Par_r$ be a core, $d\in\Z_{\geq 0}$, and $n=r+dp$. Denote by $B_{\rho,d}$ the block of the symmetric group algebra $\k\Si_{n}$ corresponding to $\rho$, cf. \cite[6.1.21]{JK}. The corresponding central idempotent will be denoted $b_{\rho,d}$, so $B_{\rho,d}=\k\Si_{n}b_{\rho,d}$. 
The irreducible  $B_{\rho,d}$-modules are $\{D^\la\mid \la\in \Parreg_{\rho,d}\}$, cf. \cite[7.1.13, 7.2.13]{JK}. We also have $b_{\rho,d}S^\la=S^\la$ for all $\la\in\Par_{\rho,d}$. 
For $\la\in\Par_{\rho,d}$, we have $[S^\la:D^{\la^\reg}]=1$ and $[S^\la:D^\mu]\neq 0$ implies $\mu\unrhd\la^\reg$ and $\mu\in\Parreg_{\rho,d}$, see \cite{J2}. If $\la$ is $p$-regular, we have $\head S^\la\cong D^\la$, see \cite[$\S$11]{JamesBook}.

Let $\theta=\cont(\rho)+\sum_{i\in I}d\al_i\in Q^+_n$. We can recover $\rho$ and $d$ from $\theta$, so it is unambiguous to write 
$B_\theta$ for $B_{\rho,d}$.  Note that 
$\cont(\la)=\theta$ for all $\la\in\Par_{\rho,d}$. Now, for a general $\theta\in Q^+_n$ we set $B_\theta:=B_{\rho,d}$, $b_\theta:=b_{\rho,d}$ if $\theta=\cont(\rho)+\sum_{i\in I}d\al_i$ for some core $\rho$ and $d\in\Z_{\geq 0}$, and set $B_\theta:=0$, $b_\theta:=0$ otherwise. Then we have $
\k\Si_n=\bigoplus_{\theta\in Q^+_n}B_\theta,
$
and the corresponding decomposition $1=\sum_{\theta\in Q^+_n}b_\theta$. 


\begin{Lemma} \label{LRowRem} {\rm \cite{J3}} 
Let $\la=(l,\la_2,\la_3\dots)\in\Par_n$, $\mu=(l,\mu_2,\mu_3,\dots)\in\Parreg_n$, and set $\bar\la:=(\la_2,\la_3,\dots)\in\Par_{n-l}$, $\bar\mu:=(\mu_2,\mu_3,\dots)\in\Parreg_{n-l}$. Then $[S^\la:D^\mu]=[S^{\bar\la}:D^{\bar\mu}]$. 
\end{Lemma}

Recalling the notation of \S\ref{subsec:comb}, especially $\quot(\la)=(\la^{(0)},\dots,\la^{(p-1)})$, we have: 

 \begin{Lemma} 
 \label{lem:irr.Sp.aba} 
 {\rm \cite[Proposition 2.1]{F1}} 
 Let $\la\in \Par_n$. The Specht module $S^\la$ is irreducible if and only if $\la$ has an abacus display such that for  some $j,k\in I$ we have:
 \begin{itemize}
 \item[(i)] $\la^{(l)}=\varnothing$ for $j\neq l\neq k$;
 \item[(ii)] If position $j+pa$  on runner $j$ is unoccupied, then any position $b>j+pa$ not on runner $j$ is unoccupied;
 \item[(iii)] If position $k+pc$ on runner $k$ is occupied, then any  position $d<k+pc$ not on runner $k$ is occupied;
 \item[(iv)] the partition $\la^{(j)}$ is $p$-regular and the Specht module $S^{\la^{(j)}}$ is irreducible;
 \item[(v)] the partition $\la^{(k)}$ is $p$-restricted and the Specht module $S^{\la^{(k)}}$ is irreducible.
\end{itemize}
\end{Lemma}

Suppose $S^\la$ is irreducible and choose an abacus display for $\la$. In view of Lemma~\ref{lem:irr.Sp.aba}, if $\la$ is not a core, i.e. $\la^{(l)}\neq \varnothing$ for some $l\in I$, then $\la$ is non-$p$-regular or non-$p$-restricted.
Moreover, if $\la$ is non-$p$-restricted then there is a unique runner $j$ as in Lemma~\ref{lem:irr.Sp.aba}---this runner will be called {\em non-restricted}. 
Similarly, if $\la$ is non-$p$-regular then there is a unique runner $k$ as in Lemma~\ref{lem:irr.Sp.aba}---this runner will be called {\em non-regular}. 

For an arbitrary $p$-singular partition $\la$, Fayers provides in \cite{F2} an algorithm for going from an abacus $\Ga(\la)$ to an abacus $\Ga(\la)^\reg$ which is an abacus of the regularization $\la^\reg$ of $\la$, so we can write $\Ga(\la^\reg)=\Ga(\la)^\reg$.  
Using \cite[\S2]{F2} one can easily verify the following: 

\begin{Lemma}
\label{lem:reg1}
Let $\la$ be a $p$-singular partition such that the Specht module  $S^\la$ is irreducible. Let $k$ be the non-regular runner of an abacus display $\Ga(\la)$, and $\quot(\la^\reg)=((\la^\reg)^{(0)},\dots, (\la^\reg)^{(p-1)})$ is defined using $\Ga(\la^\reg)=\Ga(\la)^\reg$. Then:
\begin{itemize}
\item[(i)] $(\la^\reg)^{(k)}=\varnothing$;
\item[(ii)] If position $k+pa$ on runner $k$ of $\Ga(\la)^\reg$ is occupied, then every position $b<k+pa$ of $\Ga(\la)^\reg$ is occupied.
\end{itemize}
\end{Lemma}

\subsection{Translation functors}
\label{subsec:tran.fun}
We review the {\em $i$-induction} and {\em $i$-restriction}  (translation) functors, referring the reader to \cite{K3} for more details. 
Let $i\in I$, $\theta\in Q^+_n$, $r\in\Z_{\geq 0}$, and $V$ be a module over the block $B_\theta=\k\Si_n b_\theta$. Extending $V$ to a $\k\Si_n$-module, we define 
\begin{align*}
e_i^{(r)}V&:= b_{\theta-r\al_i}(V \da^{\Si_n}_{\Si_{n-r}\times\Si_r})^{\Si_r}\in\mod{B_{\theta-r\al_i}},\\  
f_i^{(r)}V&:= b_{\theta+r\al_i}(V\boxtimes\k_{\Si_r})\ua_{\Si_n\times\Si_r}^ {\Si_{n+r}}\in\mod{B_{\theta+r\al_i}}.
\end{align*}
Note by \cite[(8.13)]{K3}, the functorial isomorphism
\begin{equation}\label{ECoinv}
e_i^{(r)}V\cong b_{\theta-r\al_i}(V \da^{\Si_n}_{\Si_{n-r}\times\Si_r})_{\Si_r},
\end{equation}
where $(-)_{\Si_r}$ stands for $\Si_r$-coinvariants. 

We then extend the definition of $e_i^{(r)}V$ and $f_i^{(r)}V$ to any $\k\Si_n$-module $V$ additively and obtain the functors
$
e_i^{(r)}: \k\Si_n\mods \to \k\Si_{n-r}\mods$ and $f_i^{(r)}: \k\Si_n\mods\to \k\Si_{n+r}\mods.
$
We write $e_i:=e_i^{(1)}$ and $f_i:=f_i^{(1)}$. Then 
$
V\da_{\Si_{n-1}}\cong \bigoplus_{i\in I}e_iV$ and $ V\ua^{\Si_{n+1}}\cong \bigoplus_{i\in I}f_iV$.

\begin{Lemma}
\label{lem:div.func}
{\rm \cite[Lemma 8.2.2(ii), Theorem 8.3.2]{K3}} 
The functors $e_i^{(r)}$ and $f_i^{(r)}$ are exact, biadjoint and commute with duality. Moreover, $e_i^{r}\cong (e_i^{(r)})^{\oplus r!}$ and $f_i^{r}\cong (f_i^{(r)})^{\oplus r!}$.
\end{Lemma}

\begin{Lemma}
\label{lem:shap}
Let $k\in\Z_{\geq 0}$. For $V\in\mod{\k\Si_n}$ and $W\in\mod{\k\Si_{n-r}}$, we have 
\begin{align*}
\Ext^k_{\Si_{n-r}}(e_i^{(r)}V, W)\cong\Ext^k_{\Si_{n}}(V,f_i^{(r)}W) 
\quad {\textrm {and}}\quad 
\Ext^k_{\Si_{n}}(f_i^{(r)}W, V)\cong\Ext^k_{\Si_{n-r}}(W,e_i^{(r)}V).
\end{align*}
\end{Lemma}

\begin{proof}
This follows immediately by Lemma~\ref{lem:div.func} and Shapiro's lemma.
\end{proof}

We now record some results on the application of  $e_i^{(r)}$ and $f_i^{(r)}$ to irreducible   modules. 


\begin{Lemma}
\label{lem:res.irre}
{\rm \cite[Theorems 11.2.10, 11.2.11]{K3}} 
Let $\la\in \Parreg_n$, $i\in I$ and $r\in\Z_{\geq 1}$. Then: 
\begin{itemize}
\item[(i)] $e_i^{(r)}D^\la\neq 0$ (resp. $f_i^{(r)}D^\la\neq 0$) if and only if $r\leq \eps_i(\la)$ (resp. $r\leq \phi_i(\la)$), in which case $e_i^{(r)}D^\la$ (resp. $f_i^{(r)}D^\la$) is a self-dual indecomposable module with simple socle and head both isomorphic to $D^{ \te_i^r\la}$ (resp. $D^{ \tf_i^r\la})$;
\item[(ii)] $[e_i^{(r)}D^\la:D^{ \te_i^r\la}]={\eps_i(\la)\choose r} 
=\dim\End_{\Si_{n-r}}(e_i^{(r)}D^\la)$;
and $[f_i^{(r)}D^\la:D^{ \tf_i^r\la}]={\phi_i(\la)\choose r}=\dim\End_{\Si_{n+r}}(f_i^{(r)}D^\la)$;
\item[(iii)] If $D^\mu$ is a composition factor of $e_i^{(r)}D^\la$ (resp. $f_i^{(r)}D^\la$), then $\eps_i(\mu)\leq\eps_i(\la)-r$ (resp. $\phi_i(\mu)\leq\phi_i(\la)-r$), with equality holding if and only if $\mu= \te_i^r\la$ (resp. $\mu= \tf_i^r\la$). In particular, $e_i^{(\eps_i(\la))}D^\la\cong D^{\te_i^{\eps_i(\la)}}$ and 
$f_i^{(\phi_i(\la))}D^\la\cong D^{\tf_i^{\phi_i(\la)}}$. 
\item[(iv)] $e_i^{(r)}D^\la$ (resp $f_i^{(r)}D^\la$) is irreducible if and only if $r=\eps_i(\la)$ (resp. $r=\phi_i(\la)$).
\end{itemize}
\end{Lemma}

\vspace{2mm}
Let $\la\in\Parreg_n$, $\mu\in\Parreg_m$ and $i\in I$. 
We say that $\mu$ is an {\em $i$-reflection} of $\la$ if  
 $\eps_i(\la)=0$ and $\mu=\tf_i^{\phi_i(\la)}\la$, or $\phi_i(\la)=0$ and $\mu=\te_i^{\eps_i(\la)}\la$. If $\mu$ is an $i$-reflection for some $i$ we say that $\mu$ is a {\em reflection} of $\la$.

\begin{Lemma} \label{L200509}
Let $\mu$ be a reflection of $\la$, then $\Ext^1_{\Si_n}(D^\la,D^\la)\cong \Ext^1_{\Si_{m}}(D^\mu,D^\mu)$. 
\end{Lemma}
\begin{proof}
Follows from Lemmas~\ref{lem:shap} and \ref{lem:res.irre}. 
\end{proof}

\begin{Lemma}
 \label{lem:irr.Sp.mul}
Let $\la\in\Parreg_n$ and $i\in I$. Let $\mu:=\te_i^{\eps_i(\la)}\la$ and assume that $D^\mu\cong S^\nu$ for some  $\nu \in \Par_{n-r}$. Then
$e_{-i}^{(\eps_i(\la))}D^{\laM}\cong D^{\muM}\cong S^{\nu'}.$ 
\end{Lemma}
 \begin{proof}
By \cite[Theorem 8.15]{JamesBook} and using the self-duality of irreducible modules over symmetric groups, we get $S^\nu\otimes \sgn\cong S^{\nu'}$. By Lemma~\ref{lem:res.irre}(iii) and the assumption, we have $e_i^{(\eps_i(\la))}D^\la\cong D^\mu\cong S^\nu$. Tensoring with $\sgn$ and using the functorial isomorphism $(e_i^{(r)}-)\otimes \sgn\cong e_{-i}^{(r)}(-\otimes \sgn)$, we get the result. 
 \end{proof}


Let $\la\in \Par_n$, $i\in I$, and $\Rem(\la,i)$ (resp. $\Add(\la,i)$) be the set of all $i$-removable (resp. $i$-addable) nodes for $\la$. Let $r\in \Z_{\geq 0}$, and for a set $X$, denote by $\Om^r(X)$ the set of all $r$-element subsets of $X$. If $\bA=\{A_1,\ldots,A_r\}\in\Om^r(\Rem(\la,i))$ (resp. $\bB=\{B_1,\ldots,B_r\}\in\Om^r(\Add(\la,i))$), define $\la_\bA:=\la_{A_1,\ldots,A_r}$ (resp. $\la^\bB:=\la^{B_1,\ldots,B_r}$). 

We say that $V\in\mod{\k\Si_n}$  has a \emph{Specht filtration} if  $V\sim S^{\la^1}\mid\ldots\mid S^{\la^s}$ for some Specht modules $S^{\la^j}$ with $\la^j\in \Par_n$.

\begin{Lemma}
\label{lem:res.Sp}
Let $\la\in \Par_n$, $i\in I$ and $r\in\Z_{\geq 0}$. 
Then  $e_i^{(r)}S^\la$ (resp. $f_i^{(r)}S^\la$) has a Specht filtration with factors $\{S^{\la_\bA}\mid\bA\in \Om^r(\Rem(\la,i))\}$ (resp. $\{S^{\la^\bB}\mid\bB\in \Om^r(\Add(\la,i))\}$), each appearing once, such that the factor $S^{\la_{\bA}}$ (resp. $S^{\la^{\bB}}$) occurs above the factor $S^{\la_{\bA'}}$ (resp. $S^{\la^{\bB'}}$) whenever $\la_{\bA}\rhd\la_{\bA'}$ (resp. $\la^{\bB}\rhd\la^{\bB'}$). 
In particular, $e_i^{(r)}S^\la\neq 0$ (resp. $f_i^{(r)}S^\la\neq 0$) if and only if $r\leq \beps_i(\la)$ (resp. $r\leq \bphi_i(\la)$), in which case $S^{\he_i^r\la}$ (resp. $S^{\hf_i^r\la}$) is the top Specht factor.
\end{Lemma}

\begin{proof}
By \cite[Theorem 9.3]{JamesBook}, $e_i^rS^\la\neq 0$ if and only if $r\leq \beps_i(\la)$. Now by  Lemma~\ref{lem:div.func} it follows that $e_i^{(r)}S^\la\neq 0$ if and only if $r\leq \beps_i(\la)$ and so, when working with $e_i^{(r)}S^\la$, we may assume that $r\leq \beps_i(\la)$; in particular, $r\leq \la_1$, and the skew shape $\la/(r)$ makes sense. 

Moreover, by \cite[Theorem 9.3]{JamesBook} again, $e_i^rS^\la$ has a Specht filtration with factors $\{S^{\la_{\bA}}\mid \bA\in \Om^r(\Rem(\la,i))\}$, each appearing $r!$ times. So by Lemma~\ref{lem:div.func},  
\begin{equation}
\label{proof:dim1}
\dim e_i^{(r)}S^\la=\sum_{\bA\in \Om^r(\Rem(\la,i))} \dim S^{\la_{\bA}}.
\end{equation}
 
By \cite[Theorem 3.1]{JP} (cf. \cite[Lemma 1.3.9]{DG}), $S^\la \da_{\Si_{n-r}\times\Si_r}$ has a filtration with factors  $S^{\la/\tau}\boxtimes S^\tau$, where $\tau\in \Par_r$ and $S^{\la/\tau}$ is the Specht module corresponding to the skew shape $\la/\tau$; moreover in this filtration the factor $S^{\la/\tau}\boxtimes S^\tau$ appears above the factor $S^{\la/\sigma}\boxtimes S^\sigma$ whenever $\tau\rhd\si$. In particular,  $S^{\la/(r)}\boxtimes S^{(r)}$ is a quotient of $S^\la \da_{\Si_{n-r}\times\Si_r}^{\Si_n}$. So 
$S^{\la/(r)}$ is a quotient of coinvariants $S^\la_{\Si_r},$ and, using (\ref{ECoinv}), we deduce that 
$
b_{\theta-r\al_i}S^{\la/(r)}$ is a quotient of $e_i^{(r)}S^\la$. 

By \cite[Theorem 5.5]{JP}, the module $S^{\la/(r)}$  has a Specht filtration and its factors are given by the Littelwood-Richardson rule \cite[(9.2)]{Mac}. It follows that 
$b_{\theta-r\al_i}S^{\la/(r)}$ has a Specht filtration with factors $\{S^{\la_{\bA}}\mid \bA\in \Om^r(\Rem(\la,i))\}$ each appearing once. Using (\ref{proof:dim1}), we deduce
$
\dim b_{\theta-r\al_i}S^{\la/(r)}=\dim  e_i^{(r)}S^\la,
$
so 
$e_i^{(r)}S^\la\cong b_{\theta-r\al_i}S^{\la/(r)}$, 
which implies the result using \cite[Theorem 5.5]{JP}.

The argument for $f_i^{(r)}$ is similar but uses 
\cite[Corollary 17.14]{JamesBook} instead of \cite{JP}.
\end{proof}

\begin{Corollary} \label{CDomFactor} 
Let $\la\in \Parreg_n$, $i\in I$ and $r\in\Z_{\geq 0}$. 
If $D^\mu$ is a composition factor of 
 $e_i^{(r)}D^\la$ (resp. $f_i^{(r)}D^\la$) 
then $\mu\unrhd \la_\bA$ for some $\bA\in \Om^r(\Rem(\la,i))$ (resp. $\mu\unrhd \la^\bB$ for some $\bB\in \Om^r(\Add(\la,i))$).
\end{Corollary}
\begin{proof}
Since $D^\la$ is composition factor of $S^\la$ and composition factors of $S^{\nu}$ are of the form $D^\mu$ for $\mu\unrhd \nu$, the result follows from Lemma~\ref{lem:res.Sp}.
\end{proof}

\section{Self-extensions for RoCK blocks}
\label{SRoCK}

Throughout this section we assume that $p>2$. We prove that  there are no self-extensions for irreducible modules lying in a RoCK block.

\subsection{Notation} 
In this section `graded' always means `$\Z$-graded'. 
For a graded vector space $V=\bigoplus_{r\in \Z}V^r$ and $s\in \Z$, we denote by $q^s V$ the same vector space with grading shifted by $s$, i.e. $(q^s V)^r=V^{r-s}$. Given a finite dimensional graded $\k$-algebra $A$, the irreducible $A$-modules are gradable uniquely up to grading shift. For graded $A$-modules, we use the notation $\hom_A(V,W)$ and $\ext^t_A(V,W)$ to denote homomorphism and extension spaces in the category of graded $A$-modules. For example, $\hom_A(V,W)$ means degree $0$ homomorphisms. It is well-known, see e.g. \cite[2.4.7]{NVO} that for finite dimensional $V$ and $W$, we have 
\begin{equation}\label{E280521}
\Ext^t_A(V,W)\cong \bigoplus_{s\in\Z}\ext^t_A(V,q^sW)
\end{equation} 
where $\Ext^t_A(V,W)$ is the usual $\Ext$ in the ungraded category.

Let $\Gamma$ be the 
quiver with vertex set 
$$J:=\{1,2,\dots,p-1\}\subset I=\{0,\dots,p-1\}$$ 
and arrows
 $\za_{k,j}$ from $j$ to $k$ for all $(k,j)\in J^2$ such that $|k-j|=1$:
\begin{align*}
\begin{braid}\tikzset{baseline=3mm}
\coordinate (1) at (0,0);
\coordinate (2) at (4,0);
\coordinate (3) at (8,0);
\coordinate (4) at (12,0);
\coordinate (6) at (16,0);
\coordinate (L1) at (20,0);
\coordinate (L) at (24,0);
\draw [thin, black,->,shorten <= 0.1cm, shorten >= 0.1cm]   (1) to[distance=1.5cm,out=100, in=100] (2);
\draw [thin,black,->,shorten <= 0.25cm, shorten >= 0.1cm]   (2) to[distance=1.5cm,out=-100, in=-80] (1);
\draw [thin,black,->,shorten <= 0.25cm, shorten >= 0.1cm]   (2) to[distance=1.5cm,out=80, in=100] (3);
\draw [thin,black,->,shorten <= 0.25cm, shorten >= 0.1cm]   (3) to[distance=1.5cm,out=-100, in=-80] (2);
\draw [thin,black,->,shorten <= 0.25cm, shorten >= 0.1cm]   (3) to[distance=1.5cm,out=80, in=100] (4);
\draw [thin,black,->,shorten <= 0.25cm, shorten >= 0.1cm]   (4) to[distance=1.5cm,out=-100, in=-80] (3);
\draw [thin,black,->,shorten <= 0.25cm, shorten >= 0.1cm]   (6) to[distance=1.5cm,out=80, in=100] (L1);
\draw [thin,black,->,shorten <= 0.25cm, shorten >= 0.1cm]   (L1) to[distance=1.5cm,out=-100, in=-80] (6);
\draw [thin,black,->,shorten <= 0.25cm, shorten >= 0.1cm]   (L1) to[distance=1.5cm,out=80, in=100] (L);
\draw [thin,black,->,shorten <= 0.1cm, shorten >= 0.1cm]   (L) to[distance=1.5cm,out=-100, in=-100] (L1);
\blackdot(0,0);
\blackdot(4,0);
\blackdot(8,0);
\blackdot(20,0);
\blackdot(24,0);
\draw(0,0) node[left]{$1$};
\draw(4,0) node[left]{$2$};
\draw(8,0) node[left]{$3$};
\draw(14,0) node {$\cdots$};
\draw(17.45,0) node[right]{$p-2$};
\draw(21.45,0) node[right]{$p-1$};
\draw(2,1.2) node[above]{$\za_{2,1}$};
\draw(6,1.2) node[above]{$\za_{3,2}$};
\draw(10,1.2) node[above]{$\za_{4,3}$};
\draw(18,1.2) node[above]{$\za_{p-3,p-2}$};
\draw(22,1.2) node[above]{$\za_{p-1,p-2}$};
\draw(2,-1.2) node[below]{$\za_{1,2}$};
\draw(6,-1.2) node[below]{$\za_{2,3}$};
\draw(10,-1.2) node[below]{$\za_{3,4}$};
\draw(18,-1.2) node[below]{$\za_{p-3,p-2}$};
\draw(22,-1.2) node[below]{$\za_{p-2,p-1}$};
\end{braid}
\end{align*}
The {\em zigzag algebra $\Zig$} is the integral path algebra $\Z\Gamma$ modulo the following relations:
\begin{enumerate}
\item All paths of length three or greater are zero.
\item All paths of length two that are not cycles are zero.
\item All cycles of length $2$ based at the same vertex are equal.
\end{enumerate}
Length zero paths yield the standard idempotents $\{ \ze_j\mid j\in J\}$ with $ \ze_i  \za_{i,j} \ze_j= \za_{i,j}$ for all admissible $i,j$. For every $j\in J$, define
$
 \zc_j:= \za_{j,j\pm1} \za_{j\pm1,j}.
$

The algebra $\Zig$ is graded by the path length: 
$\Zig=\Zig^0\oplus \Zig^1\oplus \Zig^2.
$ 
We consider $\Zig$ as a superalgebra with 
$\Zig_\0=\Zig^0\oplus \Zig^2\quad \text{and}\quad \Zig_\1=\Zig^1.
$ For $\eps\in\Z/2\Z$ and $a\in \Zig_\eps\setminus \{0\}$ we denote $\bar a:=\eps$. We have a basis 
$\B_{\0}:=\{ \ze_i,\zc_j\mid i\in J\}$ of $\Zig_\0$, a basis 
$\B_\1:=\{ \za_{ij} \mid |i-j|=1\}$ of $\Zig_\1$, and a basis $\B:=\B_\0\sqcup \B_\1$ of $\Zig$. 

Let $d\in\Z_{\geq 0}$ and $m\in\Z_{>0}$. We set $[m]:=\{1,2,\dots,m\}$. 
For a set $X$ and $d\in\Z_{\geq 0}$ we often write $x_1\cdots x_d:=(x_1,\dots,x_d)\in X^d$. For $x\in X$, we often denote  $x^d:=x\cdots x\in X^d$. 
The symmetric group $\Si_d$ acts on the right on $X^d$ by place permutations:
$
(x_1\cdots x_d)\si=x_{\si 1}\cdots x_{\si d}.
$
If $X_1,\dots,X_N$ are sets, then $\Si_d$ acts on $X_1^d\times\dots\times X_N^d$ diagonally. 
We write $(\bx^1,\dots,\bx^N)\sim (\by^1,\dots,\by^N)$ if $(\bx^1,\dots,\bx^N)\si=(\by^1,\dots,\by^N)$ for some $\si\in \Si_d$.

Let $\zH\subseteq \Zig$ be a set of non-zero homogeneous elements of $\Zig$; in particular, $\zH=\zH_\0\sqcup \zH_\1$  where $\zH_\eps:=\zH\cap \Zig_\eps$ for $\eps\in\Z/2\Z$. Define $\Seq^\zH (m,d)$ to be the set of all triples 
$$
(\bz,\br, \bs) = ( \zz_1\cdots \zz_d,\, r_1\cdots r_d,\, s_1\cdots s_d ) \in  \zH^d\times[m]^d\times[m]^d
$$
such that for all $1\leq k\neq l\leq d$ we have 
$(\zz_k,r_k,s_k)=(\zz_l,r_l, s_l)$ 
only if $\zz_k\in \zH_\0$. Then $\Seq^\zH (m,d)\subseteq \zH^d\times[m]^d\times[m]^d$ is $\Si_d$-invariant.  
We choose a set $\Seq^\zH (m,d)/\Si_d$ of $\Si_d$-orbit representatives (and identify it with the set of $\Si_d$-orbits on $\Seq^\zH (m,d)$). 

We fix a total order `$<$' on $\zH\times[m]\times[m]$. 
Then we also have a total order on $\Seq^\zH(m,d)$ defined as follows: $(\bz,\br,  \bs)< (\bz',\br',  \bs')$ if and only if there exists $l\in\{1,\dots, d\}$ such that $(\zz_k,r_k,s_k)=(\zz_k',r_k',s_k')$ for all $k<l$ and $(\zz_l,r_l,s_l)<(\zz_l',r_l',s_l')$. 
For $(\bz,\br, \bs) \in \Seq^\zH(m,d)$ and $\si\in\Si_d$, 
we define
\begin{align*}
\lan\bz, \br,  \bs\ran
&:=\sharp\{(k,l)\in[d]^2\mid k<l,\ \zz_k,\zz_l\in \zH_\1,\ (\zz_k,r_k,s_k)> (\zz_l,r_l, s_l)\},
\\
\lan \si;\bz\ran&:=\sharp\{(k,l)\in[d]^2\mid k<l,\  \si^{-1}k>\si^{-1}l,\ \zz_k,\zz_l\in \zH_\1\}.
\end{align*}

We denote by 
$\La(m,d)$ the set of compositions of $d$ with $m$ (non-negative) parts. 
Set
\begin{align*}
\La_+(m,d)&:=\{\la=(\la_1,\dots,\la_m)\in\La(m,d)\mid \la_1\geq\dots\geq \la_m\}.
\\
\La^J_+(m,d)&:=\bigsqcup_{d_1+\dots+d_{p-1}=d}\La_+(m,d_1)\times\dots\times \La_+(m,d_{p-1}).
\end{align*}
Let $S(m,d)$ be the classical Schur algebra, see \cite{Green}. 
The irreducible $S(m,d)$-modules are 
\begin{equation}\label{EIrrClSchur}
\{L(\la)\mid\la\in \La_+(m,d)\},
\end{equation}
where $L(\la)$ is the irreducible $S(m,d)$-module with highest weight $\la$, see \cite[3.5a]{Green}. 

\subsection{Zigzag Schur algebras}
Let $M_m(\Zig)$ be the superalgebra of $m\times m$ matrices with entries in $\Zig$. For $\zz\in \Zig$, we denote by  
$
\xi_{r,s}^\zz\in M_m(\Zig)
$ 
the matrix with $\zz$ in the position $(r,s)$ and zeros elsewhere. 
We write 
$$
\xi_{r,s}:=\xi_{r,s}^{1_\Zig}=\sum_{j\in J}\xi_{r,s}^{\ze_j}.
$$

The group $\Si_d$ acts on $M_m(\Zig)^{\otimes d}$ on the right by algebra automorphisms, such that for all $r_1,s_1,\dots,r_d,s_d\in [m]$, $\si\in \Si_d$ and homogeneous $\zz_1,\dots,\zz_d\in \Zig$, we have 
$$(\xi_{r_1,s_1}^{\zz_1}\otimes\dots\otimes \xi_{r_d,s_d}^{\zz_d})^\si=
(-1)^{\lan\si;\bz\ran} \xi_{r_{\si1},s_{\si1}}^{\zz_{\si1}}\otimes\dots\otimes \xi_{r_{\si d},s_{\si d}}^{\zz_{\si d}}.
$$

We consider the superalgebra of invariants 
$
S^\Zig(m,d)_\Z:=\big(M_m(\Zig)^{\otimes d}\big)^{\Si_d}.
$
Note that $S^{\Zig^0}(m,d)_\Z$ is a natural subalgebra of the even part  $S^\Zig(m,d)_{\Z,\0}$. Moreover, the algebra $S^\Zig(m,d)_\Z$ inherits a (non-negative) grading from $\Zig$,  with the degree zero component 
$S^{\Zig}(m,d)_\Z^0$ being exactly $S^{\Zig^0}(m,d)_\Z.$  
For $(\bz,\br,\bs)\in\Seq^\zH(m,d)/\Si_d$, we have elements 
\begin{equation*}\label{EXiDef}
\xi_{\br,\bs}^\bz:= \sum_{(\bz',\br',\bs')\sim(\bz,\br,\bs)} 
(-1)^{\lan\bz,\br,\bs\ran+\lan\bz',\br',\bs'\ran}
\xi_{r_1',s_1'}^{\zz_1'}\otimes\dots\otimes \xi_{r_d',s_d'}^{\zz_d'}
\in S^\Zig(m,d)_\Z.
\end{equation*}
Note that the similarly defined 
$
S(m,d)_\Z:=\big(M_m(\Z)^{\otimes d}\big)^{\Si_d}
$
is a $\Z$-form of the classical Schur algebra $S(m,d)$ with standard Schur's basis elements $\xi_{\br,\bs}$ as in \cite{Green}.

Let $\la=(\la_1,\dots,\la_l)$ be a composition of $d$. We have the standard parabolic subgroup $\Si_\la:=\Si_{\la_1}\times\dots\times\Si_{\la_l}\leq \Si_d,$ and we denote by\, ${}^{\la}\D$ the set of the shortest coset representatives for $\Si_{\la}\backslash\Si_d$. 
Given $\xi_1\in M_m(\Zig)^{\otimes \la_1},\dots,\xi_l\in M_m(\Zig)^{\otimes \la_l}$, we define 
\begin{equation}\label{EStarNotation}
\xi_1*\cdots* \xi_l:=\sum_{\si\in{}^{\la}\D}(\xi_1\otimes \dots\otimes \xi_l)^\si.
\end{equation}

\begin{Lemma} {\rm \cite[Lemma 3.3]{KMTwo})}  \label{LBasis'} 
The set $
\{\xi_{\br,\bs}^\bb\mid (\bb,\br,\bs)\in\Seq^\B(m,d)/\Si_d\}
$ is a $\Z$-basis of $S^\Zig(m,d)_\Z$. In particular, 
$$
\bigsqcup_{d_1+\dots+d_{p-1}=d}\{\xi_{\br^1,\bs^1}^{\ze_1^{d_1}}*\dots* \xi_{\br^{p-1},\bs^{p-1}}^{\ze_{p-1}^{d_{p-1}}}\mid (\br^j,\bs^j)\in([m]^{d_j}\times[m]^{d_j})/\Si_{d_j}\ \text{for all $j\in J$}\}
$$ 
is a $\Z$-basis of the degree zero part $S^{\Zig^0}(m,d)_\Z$, and there is an isomorphism of algebras 
\begin{align*}
S^{\Zig^0}(m,d)_\Z
&\stackrel{{{}_\sim}}{\rightarrow}  \bigoplus_{d_1+\dots+d_{p-1}=d}S(m,d_1)_\Z\otimes\dots\otimes S(m,d_{p-1})_\Z,
\\
\xi_{\br^1,\bs^1}^{\ze_1^{d_1}}*\dots* \xi_{\br^{p-1},\bs^{p-1}}^{\ze_{p-1}^{d_{p-1}}}&\mapsto \xi_{\br^1,\bs^1}\otimes\dots\otimes \xi_{\br^{p-1},\bs^{p-1}}.
\end{align*}
\end{Lemma}

Let $(\bb,\br,  \bs) \in \Seq^\B(m,d)/\Si_d$. We denote 
\begin{align*}
\eta^\bb_{\br,\bs}:=\left(\prod_{j\in J,\, r,s\in [m]}
|\{k\in[d]\mid  (\zb_k,r_k,s_k)=(\zc_j,r,s)\}|!\right) \xi^\bb_{\br,  \bs}. 
\end{align*}
Let 
$$
T^{\Zig}(m,d)_\Z:=\spa_\Z\big\{\,\eta^\bb_{\br,\bs}\mid (\bb,\br,\bs)\in\Seq^\B(m,d)/\Si_d\,\big\}\subseteq S^{\Zig}(m,d)_\Z.
$$ 
By Lemma~\ref{LBasis'}, $\big\{\,\eta^\bb_{\br,\bs}\mid (\bb,\br,\bs)\in\Seq^\B(m,d)/\Si_d\,\big\}$ is a $\Z$-basis of $T^{\Zig}(m,d)_\Z$. By \cite[Proposition 3.11]{KMTwo}, $T^{\Zig}(m,d)_\Z$ is a unital graded $\Z$-subalgebra of $S^{\Zig}(m,d)_\Z$.

\begin{Theorem} \label{PSubalgebra} {\rm \cite[Theorem 7.4]{EK1}} Let $m\geq d$. Then 
$T^{\Zig}(m,d)_\Z$ is the unital $\Z$-subalgebra of $S^\Zig(m,d)_\Z$ generated by $S^{\Zig^0}(m,d)_\Z$ and the set 
$$
\{\xi^{\zz}_{1,1}*\xi_{2,2}^{\otimes \la_2}*\dots*\xi_{m,m}^{\otimes \la_m}\mid \zz\in\Zig,\,\la_2,\dots,\la_m\in\Z_{\geq 0},\ \la_2+\dots+\la_m=d-1\}.
$$
\end{Theorem}

We now extend scalars to $\k$ and denote $T^\Zig(m,d):=\k\otimes_\Z T^\Zig(m,d)_\Z$, 
$T^\Zig(m,d)^0=\k\otimes_\Z T^\Zig(m,d)^0_\Z$, 
$\eta^\bb_{\br,\bs}:=1\otimes\eta_{\br,\bs}^\bb$, etc. 
The algebra $T^\Zig(m,d)$  inherits the (non-negative) grading from $T^\Zig(m,d)_\Z$, and by Lemma~\ref{LBasis'}, we have \begin{equation}\label{EDeg0}
T^\Zig(m,d)^0=S^{\Zig^0}(m,d)\cong \bigoplus_{d_1+\dots+d_{p-1}=d}S(m,d_1)\otimes\dots\otimes S(m,d_{p-1}).
\end{equation}

\begin{Proposition} \label{PGen}
Let $m\geq d$. Then 
$T^{\Zig}(m,d)$ is the unital subalgebra of $S^\Zig(m,d)$ generated by the degree $0$ elements of the form 
\begin{equation}\label{EDeg0Gen}
\xi_{\br^1,\bs^1}^{\ze_1^{d_1}}*\dots* \xi_{\br^{p-1},\bs^{p-1}}^{\ze_{p-1}^{d_{p-1}}}\qquad(\br^j,\bs^j\in[m]^{d_j}\ \text{for all $j\in J$})
\end{equation}
and the degree $1$ elements of the form 
\begin{equation}\label{EDeg1Gen}
\xi^{\za_{i,j}}_{1,1}*\xi_{\br^1,\br^1}^{\ze_1^{d_1}}*\dots*\xi_{\br^{p-1},\br^{p-1}}^{\ze_{p-1}^{d_{p-1}}}
\qquad(\br^j\in\{2,\dots,m\}^{d_j}\ \text{for all $j\in J$}).
\end{equation}
\end{Proposition}
\begin{proof}
We use Theorem~\ref{PSubalgebra}. Note that the degree of the generator $\xi^{\zz}_{1,1}*\xi_{2,2}^{\otimes \la_2}*\dots*\xi_{m,m}^{\otimes \la_m}$ is the degree of $\zz$. So such generators with $\deg(\zz)=0$ belong to $T^\Z(m,d)^0$, which in view of (\ref{EDeg0}) and Lemma~\ref{LBasis'} is generated by the elements of the form (\ref{EDeg0Gen}).  On the other hand, if $\deg(\zz)=1$, we may assume that $\zz$ is of the form $\za_{ij}$, and we can write   $\xi^{\za_{ij}}_{1,1}*\xi_{2,2}^{\otimes \la_2}*\dots*\xi_{m,m}^{\otimes \la_m}$ as a linear combination of  generators of the form (\ref{EDeg1Gen}).

Finally, suppose $\deg(\zz)=2$, in which case we may assume that $\zz=\zc_j$ for some $j\in J$. Note that for $i$ with $|i-j|=1$ we have $\zc_j=\za_{j,i}\za_{i,j}$, therefore
\begin{align*}\xi^{\zc_j}_{1,1}*\xi_{2,2}^{\otimes \la_2}*\dots*\xi_{m,m}^{\otimes \la_m}
=(\xi^{\za_{j,i}}_{1,1}*\xi_{2,2}^{\otimes \la_2}*\dots*\xi_{m,m}^{\otimes \la_m})(\xi^{\za_{i,j}}_{1,1}*\xi_{2,2}^{\otimes \la_2}*\dots*\xi_{m,m}^{\otimes \la_m}).
\end{align*}
Since we can write   $\xi^{\za_{j,i}}_{1,1}*\xi_{2,2}^{\otimes \la_2}*\dots*\xi_{m,m}^{\otimes \la_m}$ and $\xi^{\za_{i,j}}_{1,1}*\xi_{2,2}^{\otimes \la_2}*\dots*\xi_{m,m}^{\otimes \la_m}$ as  linear combinations of generators of the form (\ref{EDeg1Gen}), the result follows. 
\end{proof}

Recall from (\ref{EDeg0}) that the algebra $T^\Zig(m,d)$ is non-negatively graded with the degree zero component being a direct sum of tensor products of classical Schur algebras. Denoting $T^\Zig(m,d)^{>0}:=\bigoplus_{m>0}T^\Zig(m,d)^m$, we have 
\begin{equation}\label{ET0}
T^\Zig(m,d)/T^\Zig(m,d)^{>0}\cong T^\Zig(m,d)^0\cong  \bigoplus_{d_1+\dots+d_{p-1}=d}S(m,d_1)\otimes\dots\otimes S(m,d_{p-1}).
\end{equation}
So the modules over the algebra in the right hand side of (\ref{ET0}) can be considered as modules over $T^\Zig(m,d)$ by inflation.

Let $\bla=(\la^{(1)},\dots,\la^{(p-1)})\in \La^J_+(m,d)$, so for each $j\in J$, we have $\la^{(j)}\in\La_+(m,d_j)$ for some $d_j\in\Z_{\geq 0}$. 
Recalling (\ref{EIrrClSchur}), consider the irreducible $(S(m,d_1)\otimes\dots\otimes S(m,d_{p-1}))$-module 
$L(\la^{(1)})\boxtimes\dots\boxtimes L(\la^{(p-1)})$, extend it trivially to the module over the right hand side of (\ref{ET0}), and then inflate to $T^\Zig(m,d)$ to get the irreducible $T^\Zig(m,d)$-module denoted $L(\bla)$. Note that $L(\bla)$ is concentrated in degree $0$. As $T^\Zig(m,d)$ is non-negatively graded, we get:

\begin{Lemma} \label{LIrrT} 
Up to isomorphism, 
$
\{q^s L(\bla)\mid s\in\Z,\,
\bla\in \La^J_+(m,d)\}
$
is a complete irredundant set of irreducible graded $T^\Zig(m,d)$-modules. 
\end{Lemma}

\subsection{Extensions of irreducible modules over zigzag Schur algebras}
We now study the extensions of irreducible modules over $T^\Zig(m,d)$. The trivial shift case is easily reduced to the extensions over classical Schur algebras in view of the following general lemma:

\begin{Lemma} \label{LGenExt}
Let $t\in\Z_{\geq 0}$, $A=\bigoplus_{r\geq 0}A^r$ be a non-negatively graded finite dimensional $\k$-algebra and $V,W$ be $A^0$-modules considered as graded $A$-modules concentrated in degree $0$ via inflation along $A^{>0}$. Then: 
\begin{enumerate}
\item[{\rm (i)}] $\ext^t_A(V,q^sW)\neq 0$ implies $s\geq 0$;
\item[{\rm (ii)}] $\ext^t_A(V,W)\cong\Ext^t_{A^0}(V,W)$. 
\end{enumerate}
\end{Lemma}
\begin{proof}
By assumption, there exists a projective resolution $\dots\to P_1\to P_0\to V$ with each $P_t$ concentrated in non-negative degrees, which already implies (i). Considering the degree $0$ component $P_t^0$ of each $P_t$, we get an exact sequence of $A^0$-modules $\dots\to P_1^0\to P_0^0\to V$, with $\hom_A(P_t,W)\cong \Hom_{A^0}(P_t^0,W)$. To complete the proof of (ii), it  remains to notice that each $P_t^0$ is a projective $A^0$-module. 
\end{proof}

\begin{Corollary} \label{CZeroShift} 
Let $\bla=(\la^{(1)},\dots,\la^{(p-1)})$ and $\bmu=(\mu^{(1)},\dots,\mu^{(p-1)})$ be elements of $\La^J(m,d)$ with $\la^{(j)}\in\La_+(m,d_j)$ and $\mu^{(j)}\in\La_+(m,c_j)$ for all $j\in J$. 
Then for any $t\geq 0$, we have $\ext^t_{T^\Zig(m,d)}(L(\bla),L(\bmu))=0$, unless $c_j=d_j$ for all $j\in J$, in which case 
$$\ext^t_{T^\Zig(m,d)}(L(\bla),L(\bmu))\cong
\bigoplus_{t_1+\dots+t_{p-1}=t}\bigotimes_{j\in J}\Ext^{t_j}_{S(m,d_j)}(L(\la^{(j)}),L(\mu^{(j)})).
$$
\end{Corollary}
\begin{proof}
This follows from Lemma~\ref{LGenExt}(ii), (\ref{ET0}) and the K\"unneth theorem.
\end{proof}

\begin{Theorem} \label{TExtT}
Suppose $m\geq d$. Let $\bla=(\la^{(1)},\dots,\la^{(p-1)})$ and $\bmu=(\mu^{(1)},\dots,\mu^{(p-1)})$ be elements of $\La^J(m,d)$ with $\la^{(j)}\in\La_+(m,d_j)$ and $\mu^{(j)}\in\La_+(m,c_j)$ for all $j\in J$. Then $\ext^1_{T^\Zig(m,d)}(L(\bla),q^sL(\bmu))\neq 0$ only if one of the following two conditions holds:
\begin{enumerate}
\item[{\rm (i)}] $s=0$, $c_j=d_j$ for all $j\in J$, and there exists $j\in J$ such that the following two conditions hold: (a) $\Ext^1_{S(m,d_j)}(L(\la^{(j)}),L(\mu^{(j)}))\neq 0$, (b) $\la^{(i)}=\mu^{(i)}$ for all $i\neq j$. In this case we have 
$$
\ext^1_{T^\Zig(m,d)}(L(\bla),L(\bmu))\cong \Ext^1_{S(m,d_j)}(L(\la^{(j)}),L(\mu^{(j)})).
$$
\item[{\rm (ii)}] $s=1$ and there exist $i,j\in J$ such that the following four conditions hold: (a) $|i-j|=1$, (b) $c_i=d_i+1$, (c) $c_j=d_j-1$, (d) $c_k=d_k$ for all $k\neq i,j$. 
\end{enumerate}
\end{Theorem}
\begin{proof}
If $s=0$, the result comes from Corollary~\ref{CZeroShift}. Let $s\neq 0$. 
Consider an extension
$$
0\to q^sL(\bmu)\stackrel{\iota}{\to} E\to L(\bla)\to 0
$$
in the category of graded $T^\Zig(m,d)$-modules. We prove that the extension splits unless the condition (ii) holds. Recall that $E^n$ denotes the degree $n$ component of the graded vector space $E$. We have $E=E^s\oplus E^0$ as vector spaces, with $E^s=\iota(q^sL(\bmu))$ being a $T^\Zig(m,d)$-submodule and it suffices to prove that $E^0$ is a $T^\Zig(m,d)$-submodule. By Proposition~\ref{PGen}, $T^\Zig(m,d)$ is generated by degree $0$ elements together with degree $1$ elements of the form (\ref{EDeg1Gen}). Degree zero elements leave $E^0$ invariant, and degree $1$ elements send $E^0$ to $E^1$, so we may assume that $s=1$, in which case all elements of the form (\ref{EDeg1Gen}) still annihilate $E^0$ unless there exist $i,j\in J$ such that the conditions (a)--(d) in (ii) hold. 
\end{proof}

It is a classical fact going back to \cite{GreenComb} that the module category over $S(m,d)$ is a highest weight category, cf.  \cite[(2.5.3)]{Parshall}. In particular, by \cite[Lemma 3.2(b)]{CPS}, $S(m,d)$ has no non-trivial self-extensions. So from the theorem and (\ref{E280521}) we get:

\begin{Corollary} \label{CSelExtT} 
Suppose $m\geq d$. Then $\Ext^1_{T^\Zig(m,d)}(L,L)= 0$ for any irreducible  $T^\Zig(m,d)$-module $L$. 
\end{Corollary}

\subsection{RoCK blocks} We refer the reader to \cite{CK,T,EK2} for the information and notation concerning RoCK blocks of the symmetric groups. Our conventions are as in \cite[\S5]{EK2}. Let $d\in\Z_{\geq 0}$ and $\rho\in\Par_r$ be a {\em $d$-Rouquier core}. This means that $\rho$ is a core and there is an abacus display for $\rho$ which has at least $d-1$ more beads on runner $i+1$ than on runner $i$ for all $i=0,\dots,p-2$. 
Let $n=r+dp$. Recalling the notation of \S\ref{SSSG}, 
the block $B_{\rho,d}$ is then called a {\em RoCK block}. The algebra $B_{\rho,d}$ has a {\em KLR grading}, see \cite{BK,Ro}.

\begin{Theorem} \label{TEquiv} {\rm \cite{EK2}}
Let $d\in\Z_{\geq 0}$, $\rho\in\Par_r$ be a $d$-Rouquier core and $m\geq d$. Then  $B_{\rho,d}$ and $T^\Zig(m,d)$ are Morita equivalent as graded algebras.
\end{Theorem}

We can now prove Theorem~\ref{TA}:

\begin{Corollary} \label{CRoCK} 
Let $\la\in\Parreg_{\rho,d}$ for a $d$-Rouquier core $\rho$. Then 
$\Ext^1_{\Si_n}(D^\la,D^\la)=0.$
\end{Corollary}
\begin{proof}
The result follows immediately from Theorem~\ref{TEquiv} and Corollary~\ref{CSelExtT}.
\end{proof}

\begin{Remark} \label{R100420} 
{\rm 
The Morita equivalence of Theorem~\ref{TEquiv} can be used to translate the rest of Theorem~\ref{TExtT} into the language of symmetric groups, using the observation that under the Morita equivalence the irreducible $T^\Zig(m,d)$-module $L(\bla)$ with $\bla=(\la^{(1)},\dots,\la^{(p-1)})$ corresponds to $D^\la$ where $\la\in\Par_{\rho,d}$ is the partition with $\quot(\la)=(\varnothing,\la^{(1)},\dots,\la^{(p-1)})$. We sketch the proof of the latter fact. One needs to observe, using the formal characters of Specht modules of \cite{BKW} and \cite[Corollary 6.23]{EK2}, that under the Morita equivalence the Specht module $S^\la$ corresponds to a $T^\Zig(m,d)$-module $\De^\la$ such that the weight $\bla:=(\la^{(1)},\dots,\la^{(p-1)})$ appears in the formal character of $\De^\la$, and $\bmu\in \La^J(m,d)$ appears in the formal character of $\De^\la$ only if $\bmu\unrhd\bla$. 
Here the dominance order $\unrhd$ on $p$-multipartitions is defined by moving boxes up within a component or to the bigger component. 
Then the result follows by induction on $\unrhd$ starting with the largest multipartition $(\varnothing,\dots,\varnothing,(d))$. 
}
\end{Remark}

\section{Translation functors}
\label{sec:bra.ru}
In this section, we do not assume $p>2$.

\subsection{On the structure of $e_i^{(r)}D^\la$ and $f_i^{(r)}D^\la$}
\label{sec:div.power}
Throughout the subsection, we fix $\la\in \Parreg_n$ and $i\in I$. Recall from Lemma~\ref{lem:res.irre}(i) that for $r\leq \eps_i(\la)$ (resp. $r\leq \phi_i(\la)$)  the module $e_i^{(r)}D^\la$ (resp. $f_i^{(r)}D^\la$)  has simple socle and head both isomorphic to $D^{ \te_i^r\la}$ (resp. $D^{ \tf_i^r\la}$). It will be crucial for us to analyze the quotient $(e_i^{(r)}D^\la)/D^{\te_i^r\la}$ (resp. $(f_i^{(r)}D^\la)/D^{\tf_i^r\la}$). Results regarding the structure of these quotients have been provided in \cite[\S3.2]{KMT} for the special case of $r=1$. In this subsection we generalize these results to $r>1$.

\begin{Lemma} \label{LSpecialCaseofLuciaLemma} 
Let $\nu\in\Parreg_m$, $0< t<\phi_i(\nu)$, and $B_1,\dots,B_{\phi_i(\nu)}$ be the $i$-conormal nodes of $\nu$ labeled from top to bottom. If $\nu^{B_{t+1}}$ is $p$-regular, then
$$
\dim\Hom_{\Si_{m+1}}(e_i^{(t-1)} D^{\tilde f_i^t\nu},f_i D^\nu/D^{\tilde f_i \nu})\leq t-1.
$$ 
\end{Lemma}
\begin{proof}
By \cite[Lemma 3.11]{KMT} and duality, there exist submodules $$0\subseteq V_1\subseteq \dots\subseteq V_{\phi_i(\nu)}= f_i D^\nu$$  
such that 
$[V_a:D^{\tilde{f}_i \nu}]=a$ and $\soc V_a\cong \head V_a\cong D^{\tilde{f}_i \nu}$ for all $a$.
On the other hand, by \cite[Remark on p.83]{BK2}, there exist submodules 
$$0=W_0\subseteq W_1\subseteq \dots\subseteq W_{\phi_i(\nu)}= f_i D^\nu$$  
such that for all $a=1,\dots,\phi_i(\nu)$, we have that $W_{a}/W_{a-1}$ is a non-zero submodule of the dual Specht module $(S^{\nu^{B_{a}}})^*$ and 
$[W_{a}/W_{a-1}:D^{\tilde{f}_i \nu}]=1.$

By the assumption that $\nu^{B_{t+1}}$ is $p$-regular, we have $D^{\nu^{B_{t+1}}}\cong\soc (W_{t+1}/W_{t})$. 
Taking $X=W_t$ and $Y=W_{t+1}$ in Lemma~\ref{LSubmoduleStructure}(ii), we deduce that 
 $V_{t+1}$ (and then any $V_a$ with $a\geq t+1$) has a composition factor $D^{\nu^{B_{t+1}}}$. On the other hand, if $B_{t+1}$ is in row $r$, then $\sum_{l\geq r}(\nu^{B_{t+1}})_l=1+\sum_{l\geq r}(\tilde f_i^t\nu)_l$. So, by Corollary~\ref{CDomFactor}, $D^{\nu^{B_{t+1}}}$ is not a composition factor of $e_i^{(t-1)} D^{\tilde f_i^t\nu}$. 

Let $\psi:e_i^{(t-1)} D^{\tilde f_i^t\nu}\to f_i D^\nu/D^{\tilde f_i \nu}$ be a non-zero homomorphism. 
We have $D^{\tf_i\nu}\cong \head e_i^{(t-1)} D^{\tilde f_i^t\nu}$, 
so $D^{\tf_i\nu}\cong \head \Im\psi$. By Lemma~\ref{LSubmoduleStructure}(i), we have $\Im\psi=V_a/D^{\tilde f_i \nu}$ for some $a$. Since $D^{\nu^{B_{t+1}}}$ is a composition factor of $V_{t+1}, V_{t+2},\dots$ but it is not a composition factor of $e_i^{(t-1)} D^{\tilde f_i^t\nu}$, it follows that $a\leq t$. Thus 
\[\dim\Hom_{\Si_{m+1}}(e_i^{(t-1)} D^{\tilde f_i^t\nu},f_i D^\nu/D^{\tilde f_i \nu})=\dim\Hom_{\Si_{m+1}}(e_i^{(t-1)} D^{\tilde f_i^t\nu}, V_{t}/D^{\tf_i\nu}).\]
Moreover since $[V_{t}/D^{\tf_i\nu}:D^{\tf_i\nu}]=t-1$ and $D^{\tf_i\nu}\cong \head e_i^{(t-1)} D^{\tilde f_i^t\nu}$, we deduce  that
\[\dim\Hom_{\Si_{m+1}}(e_i^{(t-1)} D^{\tilde f_i^t\nu}, V_{t}/D^{\tf_i\nu})\leq t-1,\]
completing the proof.
\end{proof}

\begin{Lemma} \label{lem:seq} 
For $0\leq s<  r\leq \phi_i(\la)$, there exists a unique submodule of $f_i^{(s+1)} D^{\tf_i^{r-s-1}\la}$ which is isomorphic to $f_i^{(s)} D^{\tf_i^{r-s}\la}$.
\end{Lemma}
\begin{proof}
We will repeatedly use Lemma~\ref{lem:res.irre}. Apply the exact functor $f_i^{(s)}$ to the embedding $D^{\tf_i^{r-s}\la}\into f_iD^{\tf_i^{r-s-1}\la}$ to get
\[f_i^{(s)} D^{\tf_i^{r-s}\la}\into f_i^{(s)}f_iD^{\tf_i^{r-s-1}\la}\cong (f_i^{(s+1)} D^{\tf_i^{r-s-1}\la})^{\oplus s+1},\]
and note that $\soc (f_i^{(s)} D^{\tf_i^{r-s}\la})\cong D^{\tf_i^{r}\la}\cong \soc (f_i^{(s+1)} D^{\tf_i^{r-s-1}\la})$ is simple. Thus $f_i^{(s)} D^{\tf_i^{r-s}\la}$ is isomorphic to a submodule of $f_i^{(s+1)} D^{\tf_i^{r-s-1}\la}$ by Lemma \ref{L100221}. For uniqueness apply Lemma~\ref{LUniquenessSubmodule}.
\end{proof}

Let $0\leq s<  r\leq \phi_i(\la)$. In view of Lemma~\ref{lem:seq}, we can write unambiguously $f_i^{(s)} D^{\tf_i^{r-s}\la}\subseteq f_i^{(s+1)} D^{\tf_i^{r-s-1}\la}$ and define
$$
M_{\la,i,r,s}:=(f_i^{(s+1)} D^{\tf_i^{r-s-1}\la})/(f_i^{(s)} D^{\tf_i^{r-s}\la}).
$$
Applying the exact functor $f_i^{s}$ to the embedding $D^{\tf_i^{r-s}\la}\into f_iD^{\tf_i^{r-s-1}\la}$, we define the quotient 
\[N_{\la,i,r,s}:=(f_i^{s+1} D^{\tf_i^{r-s-1}\la})/(f_i^s D^{\tf_i^{r-s}\la})
\cong f_i^s((f_iD^{\tf_i^{r-s-1}\la})/D^{\tf_i^{r-s}\la})=f_i^s M_{\la,i,r-s,0}
.\]
Moreover, 
$$
f_i^s D^{\tf_i^{r-s}\la}\cong 
(f_i^{(s)} D^{\tf_i^{r-s}\la})^{\oplus s!}
\quad\text{and}\quad 
f_i^{s+1} D^{\tf_i^{r-s-1}\la}\cong (f_i^{(s+1)} D^{\tf_i^{r-s-1}\la})^{\oplus (s+1)!}.
$$
Taking $M=f_i^s D^{\tf_i^{r-s}\la}$, $N=f_i^{s+1} D^{\tf_i^{r-s-1}\la}$, $V=f_i^{(s)} D^{\tf_i^{r-s}\la}$, and $W=f_i^{(s+1)} D^{\tf_i^{r-s-1}\la}$ in Lemma~\ref{L110221} we get:
\begin{equation}\label{EMN}
N_{\la,i,r,s}\cong M_{\la,i,r,s}^{\oplus s!}\oplus (f_i^{(s+1)} D^{\tf_i^{r-s-1}\la})^{\oplus (s+1)!-s!}.
\end{equation}

\begin{Proposition}
\label{the:quot}
Let $B_1,\ldots,B_{\phi_i(\la)}$ be the $i$-conormal nodes of $\la$ labeled from top to bottom. If $1\leq r<\phi_i(\la)$ and $\la^{B_{r+1}}\in\Parreg_{n+1}$ then $D^{\tf_i^r\la}\not\subseteq M_{\la,i,r,s}$ for all $0\leq s<r$.
\end{Proposition}

\begin{proof}
By remarks preceding the theorem, using Lemmas~\ref{lem:div.func} and \ref{lem:res.irre}, we get:
\begin{align*}
&\dim\Hom_{\Si_{n+r}}(D^{\tf_i^r\la}, M_{\la,i,r,s})\\
=\,&\frac{1}{s!}\dim\Hom_{\Si_{n+r}}(D^{\tf_i^r\la}, N_{\la,i,r,s})-s\dim\Hom_{\Si_{n+r}}(D^{\tf_i^r\la},f_i^{(s+1)} D^{\tf_i^{r-s-1}\la})\\
=\,&\frac{1}{s!}\dim\Hom_{\Si_{n+r}}(D^{\tf_i^r\la}, f_i^sM_{\la,i,r-s,0})-s\\
=\,&\frac{1}{s!}\dim\Hom_{\Si_{n+r-s}}(e_i^sD^{\tf_i^r\la}, M_{\la,i,r-s,0})-s\\
=\,&\dim\Hom_{\Si_{n+r-s}}(e_i^{(s)}D^{\tf_i^r\la}, M_{\la,i,r-s,0})-s.
\end{align*}
So it is enough to prove that $\dim\Hom_{\Si_{n+r-s}}(e_i^{(s)}D^{\tf_i^r\la}, M_{\la,i,r-s,0})\leq s$. By Lemma~\ref{lem:ad.nod} we have that $(\tf_i^{r-s-1}\la)^{B_{r+1}}\in\Parreg_{n+r-s}$, and the required inequality comes by taking $\nu=\tilde f_i^{r-s-1}\la$ and $t=s+1$ in Lemma~\ref{LSpecialCaseofLuciaLemma}. 
\end{proof}

\begin{Theorem}
\label{prop:sKMT.ind}
Let $\la\in\Parreg_n$ and $B_1,\ldots,B_{\phi_i(\la)}$ be the $i$-conormal nodes of $\la$ labeled from top to bottom. If $r\leq \phi_i(\la)$ and $D^{\tf_i^r\la}\subseteq (f_i^{(r)}D^\la)/D^{\tf_i^r\la}$ then $r<\phi_i(\la)$ and $\la^{B_{r+1}}\not\in\Parreg_{n+1}$.
\end{Theorem}

\begin{proof}
If $r=\phi_i(\la)$ then by Lemma~\ref{lem:res.irre}, $f_i^{(r)}D^\la\cong D^{\tf_i^r\la}$, so we may assume $1\leq r<\phi_i(\la)$. 
By Lemma~\ref{lem:seq}, we have a filtration 
$$
D^{\tf_i^r\la}\subseteq f_i D^{\tf_i^{r-1}\la}\subseteq \dots\subseteq f_i^{(r)} D^\la.
$$
So if $D^{\tf_i^r\la}\subseteq (f_i^{(r)}D^\la)/D^{\tf_i^r\la}$ then 
$D^{\tf_i^r\la}\subseteq (f_i^{(s+1)} D^{\tf_i^{r-s-1}\la})/(f_i^{(s)} D^{\tf_i^{r-s}\la})=M_{\la,i,r,s}$ for some $0\leq s<r$. Now 
$\la^{B_{r+1}}\not\in\Parreg_{n+1}$ by Proposition~\ref{the:quot}.
\end{proof}

\begin{Theorem}
\label{prop:sKMT.res}
Let $\la\in\Parreg_n$ and $A_1,\ldots,A_{\eps_i(\la)}$ be the $i$-normal nodes of $\la$, labeled from bottom to top. If $r\leq \eps_i(\la)$ and $D^{\te_i^r\la}\subseteq (e_i^{(r)}D^\la)/D^{\te_i^r\la}$ then  $r<\eps_i(\la)$ and $\la_{A_{r+1}}\not\in\Parreg_{n-1}$.
\end{Theorem}

\begin{proof}
If $r=\eps_i(\la)$ then by Lemma~\ref{lem:res.irre}, $e_i^{(r)}D^\la\cong D^{\te_i^r\la}$, so we may assume $1\leq r<\eps_i(\la)$. By Lemma~\ref{lem:re.nod} and \ref{lem:ad.nod}, we have  $\la_{A_{r+1}}\not\in\Parreg_{n-1}$ if and only if $(\te_i^{\eps_i(\la)}\la)^{A_r}\not\in\Parreg_{n-\eps_i(\la)+1}$. Since $A_r$ is the $(\eps_i(\la)-r+1)$-th $i$-conormal nodes of $\te_i^{\eps_i(\la)}\la$ counting from the top, we have $1\leq \eps_i(\la)-r<\phi_i(\te_i^{\eps_i(\la)}\la)$, and then $e_i^{(r)}D^\la\subseteq f_i^{(\eps_i(\la)-r)}D^{\te_i^{\eps_i(\la)}\la}$, thanks to \cite[Lemma 3.3]{M}. 
The result now follows from Theorem~\ref{prop:sKMT.ind}.
\end{proof}

\subsection{Some consequences for self-extensions}
\label{SSConsSelfExt}
The following result explains our interest in the socle of $(f_i^{(r)}D^\mu)/D^{\tilde f_i^r\mu}$. 

\begin{Lemma}
\label{LES}
Let $\la\in \Parreg_n$, $i\in I$, $r:=\eps_i(\la)$, $s:=\phi_i(\la)$, $\mu:=\te_i^{r}\la$, $\nu:=\tf_i^s\la$. We have exact sequences:
\begin{align*}
&0\to\Hom_{\Si_n}(D^\la, (f_i^{(r)}D^\mu)/D^\la)\to \Ext^1_{\Si_n}(D^\la,D^\la)\to\Ext^1_{\Si_{n-r}}(D^\mu,D^\mu),
\\
&0\to\Hom_{\Si_n}(D^\la, (e_i^{(s)}D^\nu)/D^\la)\to \Ext^1_{\Si_n}(D^\la,D^\la)\to\Ext^1_{\Si_{n+s}}(D^\nu,D^\nu).
\end{align*} 
\end{Lemma} 
\begin{proof} 
We obtain the first sequence, the argument for the second one being dual. We may assume that $r>0$. By Lemma~\ref{lem:res.irre} we have a short exact sequence $0\to D^\la\to f_i^{(r)}D^\mu\to (f_i^{(r)}D^\mu)/D^\la\to 0$. 
Apply $\Hom_{\Si_{n}}(D^\la,-)$ we get an exact sequence 
\begin{align*} 
0&\to \Hom_{\Si_n}(D^\la,D^\la)\to \Hom_{\Si_n}(D^\la,f_i^{(r)} D^\mu)\to\Hom_{\Si_n}(D^\la, (f_i^{(r)}D^\mu)/D^\la)\cr
&\to\Ext^1_{\Si_n}(D^\la,D^\la)\to \Ext^1_{\Si_n}(D^\la,f_i^{(r)} D^\mu).
\end{align*} 
By Lemmas~\ref{lem:shap} and \ref{lem:res.irre} we have 
$$\Ext^k_{\Si_n}(D^\la,f_i^{(r)} D^\mu)\cong \Ext^k_{\Si_{n-r}}(e_i^{(r)}D^\la,D^\mu)\cong \Ext^k_{\Si_{n-r}}(D^\mu,D^\mu)$$
for $k\geq 0$, and the result follows.
\end{proof}

The following result developing Lemma~\ref{L200509} is a very useful tool:

\begin{Lemma}
\label{lem:str}
Let $\la\in \Parreg_n$, $i\in I$, $r:=\eps_i(\la)$, $s:=\phi_i(\la)$, $\mu:=\te_i^{r}\la$, $\nu:=\tf_i^s\la$, $B_1,\ldots,B_{r+s}$ be the $i$-conormal nodes of $\mu$ labeled from top to bottom, and $A_1,\dots,A_{r+s}$ be the $i$-normal nodes of $\nu$ labeled from bottom to top. 
\begin{enumerate}
\item[{\rm (i)}] We have an embedding 
$\Ext^1_{\Si_n}(D^\la,D^\la)\into\Ext^1_{\Si_{n-r}}(D^\mu,D^\mu)$, unless $s>0$ and $\mu^{B_{r+1}}$ is not $p$-regular. 

\item[{\rm (ii)}] We have an embedding 
$\Ext^1_{\Si_n}(D^\la,D^\la)\into\Ext^1_{\Si_{n+s}}(D^\nu,D^\nu)$, unless $r>0$ and $\nu_{A_{s+1}}$ is not $p$-regular. 
\end{enumerate} 
\end{Lemma} 
\begin{proof} 
We prove (i), the proof of (ii) being dual. In view of Lemma~\ref{L200509} we may assume that $s>0$. 
If  $\mu^{B_{r+1}}$ is $p$-regular then taking $\la=\mu$ in  Theorem~\ref{prop:sKMT.ind}, we get $\Hom_{\Si_n}(D^\la, (f_i^{(r)}D^\mu)/D^\la)=0$. The result now follows from Lemma~\ref{LES}.
\end{proof}

\begin{Definition}\label{DDiff}
Let $\la$ be a $p$-regular partition and $i\in I$. We say that $\la$ is {\em $i$-difficult} if $\eps_i(\la),\phi_i(\la)>0$ and $\la_A^B$ is not $p$-regular for $A$ being the $i$-good and $B$ being the $i$-cogood nodes for $\la$. 
\end{Definition}

\begin{Corollary} \label{C200509} 
Let $\la\in \Parreg_n$, $i\in I$, $r:=\eps_i(\la)$, $s:=\phi_i(\la)$, $\mu:=\te_i^{r}\la$ and $\nu:=\tf_i^s\la$. If $\la$ is not $i$-difficult then we have embeddings
$$
\Ext^1_{\Si_n}(D^\la,D^\la)\into\Ext^1_{\Si_{n-r}}(D^\mu,D^\mu)\quad\text{and}\quad \Ext^1_{\Si_n}(D^\la,D^\la)\into\Ext^1_{\Si_{n+s}}(D^\nu,D^\nu).
$$ 
\end{Corollary}
\begin{proof}
We prove the first embedding, the proof of the second one is dual. Let $B_1,\ldots,B_{r+s}$ be the $i$-conormal nodes of $\mu$ labeled from top to bottom. Note that $\la=\mu^{B_1,\dots,B_r}$, $B_r$ is the $i$-good node for $\la$, and $B_{r+1}$ is the $i$-cogood node for $\la$. Moreover, in view of Lemma~\ref{lem:ad.nod}, $\mu^{B_{r+1}}$ is not $p$-regular only if $\la_{B_r}^{B_{r+1}}$ is not $p$-regular. So the result follows from Lemma~\ref{lem:str}(i).
\end{proof}

Given a partition $\la\in \Parreg_n$ and $i\in I$ with $r:=\eps_i(\la)>0$  one might hope to use Corollary~\ref{C200509}  to obtain an embedding
$\Ext^1_{\Si_n}(D^\la,D^\la)\into\Ext^1_{\Si_{n-r}}(D^{\te^r_i\la}, D^{\te^r_i\la})$ 
and proceed by induction on the degree $n$. By Lemma~\ref{lem:ad.nod}, for the $i$-good node $A$ and the $i$-cogood node $B$ for $\la$ we have that $\la_A^B$ is not $p$-regular if and only if $B=A+(p-1,-1)$. 
Therefore the critical cases are partitions of the form 
\begin{equation}\label{E200511}
\la=(\la_1,\ldots,\la_m,a+1,a^{p-2},a-1,\la_{m+p+1},\ldots),
\end{equation}
where $A:=(m+1,a+1)$ is good and $B:=(m+p,a)$ is cogood for $\la$. For some of these we still have a degree reduction:

\begin{Lemma}
\label{lem:fixed.top}
Let $p>2$ and $\la\in\Parreg_n$ be of the form $\la=((a+1)^c,a^{p-2},a-1,\ldots)$ with $a,c\geq 1$. If 
$A:=(c,a+1)$ is good and $B:=(c+p-1,a)$ is cogood for $\la$, then we have an embedding
$\Ext^1_{\Si_n}(D^\la,D^\la)\into\Ext^1_{\Si_{n-1}}(D^{\la_A},D^{\la_A}).$
\end{Lemma}

\begin{proof}
Note that $\res A=\res B$. We denote this  residue by $i$. Note that $\eps_i(\la)=1$, so $\eps_{-i}(\laM)=1$, $\phi_{-i}(\laM)=\phi_i(\la)>0$. Let $A'$ (resp. $B'$) be the $(-i)$-good ($(-i)$-cogood) node of $\laM$. If we can show that $(\laM)_{A'}^{B'}$ is $p$-regular then by Corollary~\ref{C200509}, 
there exists an embedding $\Ext^1_{\Si_n}(D^{\la^\Mull},D^{\la^\Mull})\into\Ext^1_{\Si_{n-1}}(D^{(\la_A)^\Mull},D^{(\la_A)^\Mull}),$ and the lemma follows by 
(\ref{ESignExt}). Recall from \cite{Mull,FK,BO} the Mullineux algorithm for computing $\laM$. Let $\la=\la^0,\la^1,\ldots,\la^t=\varnothing$ be the partitions obtained from $\la$ by recursively removing the $p$-rim. If $t=1$, then by assumption we have that $\la=(2,1^{p-2})=(n)^\Mull$ and then $\Ext^1_{\Si_n}(D^{\la},D^{\la})=0$, so we may assume that $t>1$.

Note that $\la^u=((a(p-1)+c-u p)^\Mull,\ldots)$ for all $0\leq u<t$ (this holds by induction by definition of the $p$-rim, the case $u=0$ holding by assumption). So by \cite[Lemma 2.2]{BKZ} we have that the first row of $(\la^u)^\Mull$ is $(\la^u)^\Mull_1=a(p-1)+c-u p$ for all $0\leq u<t$. 
In particular, $(\la^0)^\Mull_1-(\la^1)^\Mull_1=p$. Hence $\la^\Mull_1-\la^\Mull_2\geq p-1$ by definition of the Mullineux bijection (the first segment of the $p$-rim of $\la^\Mull$ has length $p$ and is contained in the first row, in particular the last $p$ nodes on the first row of $\la^\Mull$ are contained in the rim of $\la^\Mull$). So the top removable node of $\la^\Mull$ is $(1,a(p-1)+c)=(1,\la^\Mull_1)$. Since this node has residue $-i$ and is normal in $\la^\Mull$ we have from $\eps_{-i}(\laM)=1$ that $(1,a(p-1)+c)=A'$. Since $\la^\Mull_1-\la^\Mull_2\geq p-1\geq 2$ it   follows that $(\laM)_{A'}^{B'}$ is $p$-regular. 
\end{proof}

\subsection{Additional results on translation functors and some consequences}

\begin{Lemma}
\label{lem:rem.ad1}
If $D^\mu\cong S^\nu$ for some $\mu\in\Parreg_n$ and $\nu\in\Par_n$, then $\eps_i(\mu)=\beps_i(\nu)$ and  $\phi_i(\mu)=\bphi_i(\nu)$ for all $i\in I$.
\end{Lemma}
\begin{proof}
By Lemma~\ref{lem:res.irre}, $\eps_i(\mu)=\max\{r|e_i^{(r)}D^\mu\neq 0\}$. By Lemma~\ref{lem:res.Sp}, $\beps_i(\nu)=\max\{r|e_i^{(r)}S^\nu\neq 0\}$. The result for $\eps_i$'s follows. The argument for $\phi_i$'s is similar. 
\end{proof}

 \begin{Lemma}
\label{lem:irr.Sp.dfilt}
Suppose $D^\mu\cong S^\nu$ for some $\mu\in\Parreg_n$ and $\nu\in\Par_n$.  Let $i\in I$ and $0\leq r\leq \beps_i(\nu)$ (resp. 
$0\leq r\leq \bphi_i(\nu)$). Then for the Specht filtration $e_i^{(r)}S^\nu\sim S^{\nu^1}\mid\ldots \mid S^{\nu^s}$ (resp. $f_i^{(r)}S^\nu\sim S^{\nu^1}\mid\ldots \mid S^{\nu^s}$) from Lemma~\ref{lem:res.Sp} we have:
\begin{itemize}
\item[(i)] $\nu^s=\he_i^r\nu$ (resp. $\nu^s=\hf_i^r\nu$);
\item[(ii)]  $s={\beps_i(\nu)\choose r}={\eps_i(\mu)\choose r}$ (resp. $s={\bphi_i(\nu)\choose r}={\phi_i(\mu)\choose r}$);
\item[(iii)] $[S^{\nu^k}:D^{\te^r_i\mu}]=1$ for $1\leq k\leq s$ (resp. $[S^{\nu^k}:D^{\tf^r_i\mu}]=1$) for all $1\leq k\leq s$;
\item[(iv)] $\head S^{\he_i^r\nu}\cong D^{\te^r_i\mu}$ (resp. $\head S^{\hf_i^r\nu}\cong D^{\tf^r_i\mu}$). 
\end{itemize}
\end{Lemma}

\begin{proof}
We give the argument for $e$'s, the argument for $f$'s being similar. Part (i) is clear from Lemma~\ref{lem:res.Sp}. Part (ii) follows from Lemmas~\ref{lem:res.Sp} and~\ref{lem:rem.ad1}. Part (iv) follows immediately from the fact that $\head e_i^{(r)}D^\mu\cong D^{\te^r_i\mu}$, see Lemma~\ref{lem:res.irre}. To prove part (iii), let $1\leq k \leq s$. By Lemma~\ref{lem:res.Sp}, we have $e_i^{(\beps_i(\nu)-r)}S^{\nu^k}\cong S^\si$ where $\si:=\he_i^{\beps_i(\nu)}\nu$. By Lemma~\ref{lem:res.irre}, $e_i^{(\eps_i(\mu)-r)}D^{\te^r_i\mu}\cong D^\rho$ where $\rho:=\te_i^{\eps_i(\mu)}\mu$. But $\beps_i(\nu)=\eps_i(\mu)$, so $S^\si\cong D^\rho$. Hence $[S^{\nu^k}:D^{\te^r_i\mu}]\leq[S^\si: D^\rho]=1.$ However, by Lemma~\ref{lem:res.irre}, we have 
$$[e_i^{(r)}S^\nu: D^{\te^r_i\mu}]=[e_i^{(r)}D^\mu: D^{\te^r_i\mu}]={\eps_i(\mu)\choose r}=s,$$
Hence $[S^{\nu^k}:D^{\te^r_i\mu}]=1$ for every $1\leq k\leq s$. 
\end{proof}

\begin{Lemma}
\label{lem:head.nodes}
Let $\mu\in \Parreg_n$, $\nu\in\Par_n$, and $D^\mu$ be a composition factor of $S^\nu$. Then
\begin{itemize}
\item[(i)] $\eps_i(\mu)\leq\beps_i(\nu)$ and $\phi_i(\mu)\leq \bphi_i(\nu)$ for all $i\in I$;

\item[(ii)] If $D^\mu$ appear in the head of $S^\nu$ and $\eps_i(\nu)=\beps_i(\nu)$ (equiv. $\phi_i(\nu)=\bphi_i(\nu)$) for some $i\in I$, then $\eps_i(\mu)=\beps_i(\nu)$ (equiv. $\phi_i(\mu)=\bphi_i(\nu)$).
\end{itemize}
\end{Lemma} 

\begin{proof}
By Lemma~\ref{lem:res.irre}, $e_i^{(\eps_i(\mu))}D^\mu\not=0$. In particular $e_i^{(\eps_i(\mu))}S^\nu\not=0$ and then $\eps_i(\mu)\leq\beps_i(\nu)$ by Lemma~\ref{lem:res.Sp}. Similarly $\phi_i(\mu)\leq \bphi_i(\nu)$.

Clearly $\phi_i(\nu)-\eps_i(\nu)=\bphi_i(\nu)-\beps_i(\nu)$. Further $\phi_i(\mu)-\eps_i(\mu)=\phi_i(\nu)-\eps_i(\nu)$ by \cite[Lemma 8.5.8]{K3}. So $\eps_i(\nu)=\beps_i(\nu)$ and $\phi_i(\nu)=\bphi_i(\nu)$ are equivalent, as are $\eps_i(\mu)=\beps_i(\nu)$ and $\phi_i(\mu)=\bphi_i(\nu)$.

By part (i), it is enough to show that $\beps_i(\nu)\leq\eps_i(\mu)$. The assumption $\eps_i(\nu)=\beps_i(\nu)$ implies that the $i$-removable nodes of $\nu$ occur above any of its $i$-addable nodes. Let $r:=\eps_i'(\nu)$ and consider the partition $\he_i^r\nu$. By Lemma~\ref{lem:res.Sp} we have a surjection $f_i^{(r)}S^{\he_i^r\nu}\onto S^\nu,$ 
which yields an embedding
$\Hom_{\Si_n}(S^\nu,D^\mu)\into\Hom_{\Si_n}(f_i^{(r)}S^{\he_i^r\nu},D^\mu).$ Since $D^\mu$ is in the head of $S^\nu$ we deduce 
$$\Hom_{\Si_n}(S^{\he_i^r\nu},e_i^{(r)}D^\mu)\cong\Hom_{\Si_n}(f_i^{(r)}S^{\he_i^r\nu},D^\mu)\neq 0.$$ 
In particular, $e_i^{(r)}D^\mu\neq0$. Hence, by Lemma~\ref{lem:res.irre} we get that $\beps_i(\nu)=r\leq \eps_i(\mu)$.
\end{proof}

\begin{Lemma}
\label{lem:rem.node.trick}
Let $\mu\in \Parreg_n$ with $\eps_i(\mu)=0$, and $\nu\in\Par_n$ with $\beps_i(\nu)>0$. Let $A$ be the top $i$-removable node of $\nu$. If $A$ is normal for $\nu$ then
$\Hom_{\Si_n}(S^\nu,D^\mu)=0.$
\end{Lemma}

\begin{proof}
By assumption, $\nu$ does not have any $i$-addable nodes above $A$. So by Lemma~\ref{lem:res.Sp} we have a surjection 
$f_iS^{\nu_A}\onto S^\nu$, which yields an embedding 
$$\Hom_{\Si_n}(S^\nu,D^\mu)\into\Hom_{\Si_{n-1}}(f_iS^{\nu_A},D^\mu)\cong \Hom_{\Si_{n-1}}(S^{\nu_A},e_iD^\mu).$$
Since $\eps_i(\mu)=0$, we have $e_iD^\mu=0$, and the result follows.
\end{proof}

The following result is similar to Lemma~\ref{lem:rem.node.trick} and has a similar proof which we skip:

\begin{Lemma}
\label{lem:ad.node.trick}
Let $\mu\in \Parreg_n$ with $\phi_i(\mu)=0$, and $\nu\in\Par_n$ with $\bphi_i(\nu)>0$. Let $B$ be the lowest $i$-addable node of $\nu$. If $B$ is conormal for $\nu$ then
$\Hom_{\Si_n}(S^\nu,D^\mu)=0.$
\end{Lemma}

\section{Connecting to irreducible Specht modules}
\label{sec:res.irr.Sp}

In this section we prove Theorem~\ref{TB}. We assume that $p>2$, fix $\la\in\Parreg_n$, $i\in I$, and set 
$r:=\eps_i(\la),\ \mu:=\tilde e_i^r\la,$ 
so that $D^\mu\cong e_i^{(r)}D^\la$ by Lemma~\ref{lem:res.irre}. We assume that $D^\mu$ is isomorphic to a Specht module, i.e. $D^\mu\cong S^\nu$ for some $\nu\in\Par_{n-r}$. By Lemma~\ref{lem:rem.ad1}, we have 
$\phi_l(\mu)=\bphi_l(\nu)$ for all $l\in I$.  
Theorem~\ref{TB} follows from Propositions~\ref{prop:Sp.Core}, \ref{prop:Sp.reg}, \ref{prop:Sp.res} and \ref{prop:Sp.gen} proved below.

\subsection{Self-extensions for Irreducible Specht Modules}
\label{sec:irr.Sp}
Suppose for the moment that $D^\la$ itself is isomorphic to a Specht module. By Lemmas~\ref{lem:rem.ad1} and \ref{lem:res.Sp}, $D^\la$ satisfies the assumptions of Theorem~\ref{TB}, and so the equality $\Ext^1_{\Si_n}(D^\la,D^\la)=0$ is a special case of that theorem. This special case is known for  $p> 3$, see \cite[Theorem 3.3(c)]{KN}.  As a step towards the proof of Theorem~\ref{TB}, we give an independent proof that covers the case $p=3$.

\begin{Lemma}
\label{lem:seq.part}
Let $\ka\in\Par_n$ and $S^\ka$ be an irreducible Specht module. There exists a sequence of partitions $\ka=\ka^1,\ka^2,\ldots, \ka^s$ such that:
\begin{itemize}
\item[(i)]  $S^{\ka^u}$ is irreducible for $u=1,\dots,s$;
\item[(ii)] $S^{\ka^s}$  lies in a RoCK block;
\item[(iii)] For $u=1,\dots,s-1$, we have $S^{\ka^{u+1}}\cong f_{i_u}^{(\bphi_{i_u}(\ka^u))}S^{\ka^u}$ and $S^{\ka^{u}}\cong e^{(\bphi_{i_u}(\ka^u))}_{i_u}S^{\ka^{u+1}}$ for some  $i_u\in I$.
\end{itemize}
\end{Lemma}

\begin{proof}
This follows from \cite[Lemmas~3.1 and 3.2]{F1} and Lemma~\ref{lem:res.Sp}. For the irreducibility in (i) one may also argue in the following way. Since $S^{\ka^1}\cong D^\mu$ for some $\mu\in \Parreg_n$ we have that $\bphi_{i_1}(\ka^1)=\phi_{i_1}(\mu):=r$  by Lemma~\ref{lem:rem.ad1}. Now  Lemma~\ref{lem:res.irre}(iv) gives that  
$$S^{\ka^2}\cong f_{i_1}^{(r)}S^{\ka^1}\cong f_{i_1}^{(r)}D^\mu\cong D^{\tf^r_{i_1}\mu},$$
and so $S^{\ka^2}$ is irreducible. Repeating this argument one obtains immediately the result.
\end{proof}

\begin{Proposition}
\label{the:irr.Sp}
If $D^\la\cong S^\ka$ for some $\ka\in\Par_n$, then
$\Ext^1_{\Si_n}(D^\la,D^\la)=0.$  
\end{Proposition}

\begin{proof}
Let $\ka=\ka^1,\ka^2\in\Par_{n_2},\ldots, \ka^s\in\Par_{n_s}$ be as in Lemma~\ref{lem:seq.part}. 
By Lemma~\ref{lem:shap}, we have 
 $$\Ext^1_{\Si_n}(S^{\ka^1}, S^{\ka^1})\cong\Ext^1_{\Si_{n_2}}(S^{\ka^2},S^{\ka^2})\cong\ldots\cong\Ext^1_{\Si_{n_s}}(S^{\ka^s}, S^{\ka^s}).$$ Now $\Ext^1_{\Si_{n_s}}(S^{\ka^s}, S^{\ka^s})=0$ by Corollary~\ref{CRoCK} and we are done.
 \end{proof}

\begin{Corollary}
\label{CRepeated}
Let $\la\in\Parreg_n$, $r:=\eps_i(\la)$, $\mu:=\tilde e_i^r\la$, and let $B_1,\dots,B_{\phi_i(\la)}$ be the $i$-conormal nodes of $\mu$ labeled from top to bottom. Suppose that $D^\mu$ is isomorphic to a Specht module. Then 
$\Ext^1_{\Si_n}(D^\la,D^\la)=0$ unless $0<r<\phi_i(\mu)$ and $\mu^{B_{r+1}}$ is not $p$-regular. 
\end{Corollary}
\begin{proof}
By Proposition~\ref{the:irr.Sp} and Lemma~\ref{lem:shap}, we may assume that $0<r<\phi_i(\mu)$. Since $\Ext^1_{\Si_{n-r}}(D^\mu,D^\mu)=0$ by  Proposition~\ref{the:irr.Sp}, the result follows from Lemma~\ref{lem:str}. 
\end{proof}



\begin{Example}
\label{lem:hp,1}
Let $\la=(b,2,1^{p-2})$ with $b\equiv 1 \pmod{p}$. 
The Specht module $S^{(b,1^p)}$ is irreducible for example by  \cite[Theorem 23.7]{JamesBook}, and $(b,1^p)^\reg=(b,2,1^{p-2})$. Therefore, $D^\la\cong S^{(b,1^p)}$ and 
we have $\Ext^1_{\Si_{p+b}}(D^\la,D^\la)=0$ by Proposition~\ref{the:irr.Sp}.
\end{Example}

\subsection{The case $\mu=\nu$ is $p$-restricted.}
\label{subsec:Sp.core}
By Lemma~\ref{lem:irr.Sp.aba}, in this case $\mu$ is  a core.

\begin{Lemma}
\label{lem:core1}
Let $\mu$ be a core, $i\in I$, and assume that $\mu^B$ is $p$-restricted for every $i$-addable node $B$ for $\mu$. Then  
$\Ext^1_{\Si_n}(D^\la,D^\la)=0.$
\end{Lemma} 

\begin{proof}
By Lemma~\ref{lem:irr.Sp.mul}, we have  
$e_{-i}^{(\eps_{-i}(\laM))}D^{\laM}\cong S^{\mu'}.$ 
As $\mu^B$ is $p$-restricted for every $i$-addable node $B$ for $\mu$, we have that $(\mu')^C$ is $p$-regular for every $(-i)$-addable node $C$ for $\mu'$. Now the result follows from  Corollary~\ref{CRepeated} applied to $\la^\Mull$ instead of $\la$ and (\ref{ESignExt}).
\end{proof}

\begin{Proposition}
\label{prop:Sp.Core}
If $\mu$ is a core then $\Ext^1_{\Si_n}(D^\la,D^\la)=0.$
\end{Proposition}

\begin{proof}
Let $\phi:=\phi_i(\mu)=\phi_i'(\mu)$ and $B_1,\ldots,B_\phi$ be the $i$-addable nodes of $\mu$ labelled from top to bottom. By Lemmas~\ref{lem:core1} and \ref{lem:core2} we may assume that $B_1=(1,\mu_1+1)$ and $\mu^{B_1}$ is not $p$-restricted. Moreover, by Corollary~\ref{CRepeated} we may assume that 
$r<\phi$ and 
$\mu^{B_{r+1}}$ is $p$-singular.  

We choose an abacus display of $\mu$ so that the $i$-addable nodes correspond to addable positions on runner $0$. Recall from \S\ref{subsec:comb} that positions on the abacus are labeled with non-negative integers of the form $l+pa$ with $l\in I$ and $a\in\Z_{\geq 0}$. Assume that the last bead on runner $p-1$  occurs at position $p-1+pa$. Since $\mu^{B_1}$ is not $p$-restricted, we deduce that for $l\neq p-1$ the positions $l+pa$ are unoccupied. Moreover since $\mu^{B_{r+1}}$ is $p$-singular it follows that for $l\neq 0$ the positions $l+p(a-r+1)$ are occupied. Since $\mu$ is a core we 
we deduce that for $l\neq0$, the positions $l+pc$ with $c\leq a-r+1$ are occupied.  We write $s$ for the number of beads occurring on runner $p-2$ below the position $p-2+p(a-r)$. The above reasoning implies that $1\leq s<r$. An example of such a configuration is the following (it is enough to depict only the runners $p-2,p-1$ and $0$ for our considerations):
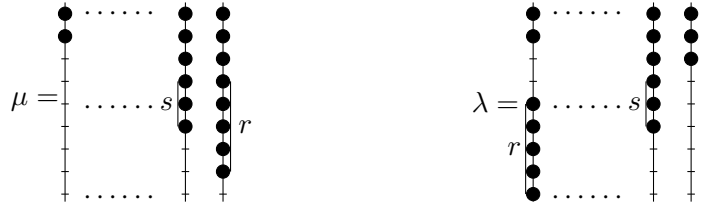
\vspace{2mm}

\hspace*{\fill}
{\begin{tikzpicture}
\draw(-.4,-1.325)node{$\mu=$};

 \draw (0,0)--(0,-2.65);
\draw (1.6,0)--(1.6,-2.65);
\draw(2.1,0)--(2.1,-2.65);

\draw(-.05,-.15)--(.05,-.15);
\draw(-.05,-.45)--(.05,-.45);
\draw(-.05,-.75)--(.05,-.75);
\draw(-.05,-1.05)--(.05,-1.05);
\draw(-.05,-1.35)--(.05,-1.35);
\draw(-.05,-1.65)--(.05,-1.65);
\draw(-.05,-1.95)--(.05,-1.95);
\draw(-.05,-2.25)--(.05,-2.25);
\draw(-.05,-2.55)--(.05,-2.55);

\draw(1.55,-.15)--(1.65,-.15);
\draw(1.55,-.45)--(1.65,-.45);
\draw(1.55,-.75)--(1.65,-.75);
\draw(1.55,-1.05)--(1.65,-1.05);
\draw(1.55,-1.35)--(1.65,-1.35);
\draw(1.55,-1.65)--(1.65,-1.65);
\draw(1.55,-1.95)--(1.65,-1.95);
\draw(1.55,-2.25)--(1.65,-2.25);
\draw(1.55,-2.55)--(1.65,-2.55);

\draw(2.05,-.15)--(2.15,-.15);
\draw(2.05,-.45)--(2.15,-.45);
\draw(2.05,-.75)--(2.15,-.75);
\draw(2.05,-1.05)--(2.15,-1.05);
\draw(2.05,-1.35)--(2.15,-1.35);
\draw(2.05,-1.65)--(2.15,-1.65);
\draw(2.05,-1.95)--(2.15,-1.95);
\draw(2.05,-2.25)--(2.15,-2.25);
\draw(2.05,-2.55)--(2.15,-2.55);

\draw(0.5,-.15)node{$\cdots$};
\draw(1,-.15)node{$\cdots$};
\draw(0.5,-1.4)node{$\cdots$};
\draw(1,-1.4)node{$\cdots$};
\draw(0.5,-2.6)node{$\cdots$};
\draw(1,-2.6)node{$\cdots$};

\filldraw [black] (0,-.15) circle (2.5pt);
\filldraw [black] (0,-.45) circle (2.5pt);

\filldraw [black] (1.6,-.15) circle (2.5pt);
\filldraw [black] (1.6,-.45) circle (2.5pt);
\filldraw [black] (1.6,-.75) circle (2.5pt);
\filldraw [black] (1.6,-1.05) circle (2.5pt);
\filldraw [black] (1.6,-1.35)circle (2.5pt);
\filldraw [black] (1.6,-1.65) circle (2.5pt);

\filldraw [black] (2.1,-.15) circle (2.5pt);
\filldraw [black] (2.1,-.45) circle (2.5pt);
\filldraw [black] (2.1,-.75) circle (2.5pt);
\filldraw [black] (2.1,-1.05) circle (2.5pt);
\filldraw [black] (2.1,-1.35)circle (2.5pt);
\filldraw [black] (2.1,-1.65) circle (2.5pt);
\filldraw [black] (2.1,-1.95) circle (2.5pt);
\filldraw [black] (2.1,-2.25) circle (2.5pt);

\draw(1.5,-1.05)--(1.5,-1.65);
\draw(1.35,-1.35)node{$s$};
\draw(2.2,-1.05)--(2.2,-2.25);
\draw(2.4,-1.65)node{$r$};
\end{tikzpicture}}
\hfill 
{\begin{tikzpicture}
\draw(-.5,-1.325)node{$\lambda=$};

 \draw (0,0)--(0,-2.65);
\draw (1.6,0)--(1.6,-2.65);
\draw(2.1,0)--(2.1,-2.65);

\draw(-.05,-.15)--(.05,-.15);
\draw(-.05,-.45)--(.05,-.45);
\draw(-.05,-.75)--(.05,-.75);
\draw(-.05,-1.05)--(.05,-1.05);
\draw(-.05,-1.35)--(.05,-1.35);
\draw(-.05,-1.65)--(.05,-1.65);
\draw(-.05,-1.95)--(.05,-1.95);
\draw(-.05,-2.25)--(.05,-2.25);
\draw(-.05,-2.55)--(.05,-2.55);

\draw(1.55,-.15)--(1.65,-.15);
\draw(1.55,-.45)--(1.65,-.45);
\draw(1.55,-.75)--(1.65,-.75);
\draw(1.55,-1.05)--(1.65,-1.05);
\draw(1.55,-1.35)--(1.65,-1.35);
\draw(1.55,-1.65)--(1.65,-1.65);
\draw(1.55,-1.95)--(1.65,-1.95);
\draw(1.55,-2.25)--(1.65,-2.25);
\draw(1.55,-2.55)--(1.65,-2.55);

\draw(2.05,-.15)--(2.15,-.15);
\draw(2.05,-.45)--(2.15,-.45);
\draw(2.05,-.75)--(2.15,-.75);
\draw(2.05,-1.05)--(2.15,-1.05);
\draw(2.05,-1.35)--(2.15,-1.35);
\draw(2.05,-1.65)--(2.15,-1.65);
\draw(2.05,-1.95)--(2.15,-1.95);
\draw(2.05,-2.25)--(2.15,-2.25);
\draw(2.05,-2.55)--(2.15,-2.55);

\draw(0.5,-.15)node{$\cdots$};
\draw(1,-.15)node{$\cdots$};
\draw(0.5,-1.4)node{$\cdots$};
\draw(1,-1.4)node{$\cdots$};
\draw(0.5,-2.6)node{$\cdots$};
\draw(1,-2.6)node{$\cdots$};

\filldraw [black] (0,-.15) circle (2.5pt);
\filldraw [black] (0,-.45) circle (2.5pt);
\filldraw [black] (0,-1.35) circle (2.5pt);
\filldraw [black] (0,-1.65)circle (2.5pt);
\filldraw [black] (0,-1.95) circle (2.5pt);
\filldraw [black] (0,-2.25) circle (2.5pt);
\filldraw [black] (0,-2.55) circle (2.5pt);

\filldraw [black] (1.6,-.15) circle (2.5pt);
\filldraw [black] (1.6,-.45) circle (2.5pt);
\filldraw [black] (1.6,-.75) circle (2.5pt);
\filldraw [black] (1.6,-1.05) circle (2.5pt);
\filldraw [black] (1.6,-1.35)circle (2.5pt);
\filldraw [black] (1.6,-1.65) circle (2.5pt);

\filldraw [black] (2.1,-.15) circle (2.5pt);
\filldraw [black] (2.1,-.45) circle (2.5pt);
\filldraw [black] (2.1,-.75) circle (2.5pt);

\draw(1.5,-1.05)--(1.5,-1.65);
\draw(1.35,-1.35)node{$s$};
\draw(-.1,-1.35)--(-.1,-2.55);
\draw(-.25,-1.95)node{$r$};
\end{tikzpicture}}
\hspace*{\fill}

Note that $\phi_{i-1}(\la)=s$ and $\beps_{i-1}(\la)=0$. Set $\xi:=\tf^s_{i-1}\la$. 
By Lemma~\ref{L200509} we have 
\begin{equation}\label{E050520}
\Ext^1_{\Si_n}(D^\la,D^\la)\cong\Ext^1_{\Si_{n+s}}(D^\xi,D^\xi).
\end{equation} 
 Note that $\phi_i(\xi)=\phi-r$ and set $\rho:=\tf_i^{\phi-r}\xi$. 
Then $\eps_{i-1}(\rho)=s$ and let $\si:=\te_{i-1}^{s}\rho$. Then $\si$ is a core and in fact $\si=\tf_i^\phi\mu$. Moreover, for any  $(i-1)$-addable node $C$ of $\si$ the partition $\si^C$ is $p$-restricted. So by Lemma~\ref{lem:core1}, we have $\Ext^1_{\Si_{n+s+\phi-r}}(D^\rho, D^\rho)=0$. Now, 
by Lemma~\ref{lem:shap},  
$$\Ext^1_{\Si_{n+s}}(D^\xi,e_i^{(\phi-r)}D^\rho)\cong 
\Ext^1_{\Si_{n+s+\phi-r}}(f_i^{(\phi-r)}D^\xi,D^\rho)\cong
\Ext^1_{\Si_{n+s+\phi-r}}(D^\rho, D^\rho)=0.$$  
So Lemma~\ref{LES} 
implies an isomorphism 
 \begin{equation}\label{E050520_2}
 \Hom_{\Si_{n+s}}(D^\xi, (e_i^{(\phi-r)}D^\rho)/D^\xi)\cong \Ext^1_{\Si_{n+s}}(D^\xi, D^\xi).
 \end{equation}
Now $\eps_i(\rho)=\phi-s$ and if $A_1,\dots,A_{\phi-s}$ are the corresponding $i$-normal nodes of $\rho$ labeled from bottom to top, one easily sees that $\rho_{A_{\phi-r+1}}$ is $p$-regular. So Theorem~\ref{prop:sKMT.res} implies that 
$\Hom_{\Si_{n+s}}(D^\xi, (e_i^{(\phi-r)}D^\rho)/D^\xi)=0$.
The result now follows from (\ref{E050520_2}) and (\ref{E050520}). 
\end{proof}

\subsection{The case $\mu=\nu$ is not $p$-restricted}
\label{subsec:Sp.reg}
Since $\mu$ is not $p$-restricted, there exists a non-restricted runner $j$ for $\Ga(\mu)$, see Lemma~\ref{lem:irr.Sp.aba}(ii). We choose $\Ga(\mu)$ so that $j=p-1$. In view of  Lemma~\ref{lem:irr.Sp.aba}, we then also have $\mu^{(l)}=\varnothing$ for all $l\neq p-1$. 
Let $\phi:=\phi_i(\mu)$ and $B_1,\ldots,B_{\phi}$ the $i$-conormal nodes of $\mu$ labelled from top to bottom. These correspond to addable positions on some runner $m$ of $\Ga(\mu)$. With this notation we have:

\begin{Lemma}
\label{lem:reg2}
If $m\neq 0$ then $\Ext^1_{\Si_n}(D^\la,D^\la)=0.$
\end{Lemma}

\begin{proof}
By (\ref{ESignExt}), it suffices to prove that 
$\Ext_{\Si_n}(D^{\la^\Mull},D^{\la^\Mull})=0$. 
Since $e_{-i}^{(r)}D^\laM\cong D^\muM\cong S^{\mu'}$ by Lemma~\ref{lem:irr.Sp.mul}, the desired equality will follow from 
Corollary~\ref{CRepeated}, once we check that $(\muM)^B$ is $p$-regular for any $(-i)$-addable node $B$ of $\muM$. 

The isomorphism $D^\muM\cong S^{\mu'}$ implies $\muM=(\mu')^\reg$. By (\ref{EAbTr}) and Lemma~\ref{lem:irr.Sp.aba}, we have $\Ga(\mu')=\Ga(\mu)'$, $0$ is the non-regular runner of $\Ga(\mu')$, and $(\mu')^{(l)}=\varnothing$ for all $l \neq 0$. Consider the abacus display $\Ga((\mu')^\reg)=\Ga(\mu')^\reg$ as in Lemma~\ref{lem:reg1}. Let $B$ be a $(-i)$-addable node for $(\mu')^\reg$, and let $t$ be the corresponding addable position in 
$\Ga((\mu')^\reg)$. Then $t$ is on runner $p-m\neq 0$. By Lemma~\ref{lem:reg1}, we have that $((\mu')^\reg)^{(0)}=\varnothing$, and if the last bead on runner $0$ in $\Ga((\mu')^\reg)$ occurs at position $pa$ then every position $b<pa$ is occupied in $\Ga((\mu')^\reg)$. Therefore $t>pa$. It follows that $((\mu')^\reg)^B$ is $p$-regular. 
\end{proof}

\begin{Proposition}
\label{prop:Sp.reg}
If $\mu=\nu$ is not $p$-restricted then 
$\Ext^1_{\Si_n}(D^\la,D^\la)=0.$
\end{Proposition}

\begin{proof}
We choose $\Ga(\mu)$ as in the beginning of this subsection. By Lemma~\ref{lem:reg2}, we may assume that the $i$-conormal nodes of $\mu$ correspond to addable positions on runner $0$. In view of Lemma~\ref{lem:irr.Sp.aba}, we must have  $B_1=(1,\mu_1+1)$. By Corollary~\ref{CRepeated} we may assume that $r<\phi$ and $\mu^{B_{r+1}}$ is $p$-singular. 

We may assume that the addable position in $\Ga(\mu)$ corresponding to $B_{r+1}$ is of the form $pa$ for some $a\in\Z_{>0}$. Then for all $l\neq0$ the positions $l+pa$ are occupied, since $\mu^{B_{r+1}}$  is $p$-singular. Since $\mu^{(l)}=\varnothing$ for $l\neq p-1$ we deduce that every position $b\leq p-1+pa$ not on runner $0$ is occupied.  Let the first unoccupied position on runner $p-1$ be $p-1+pc$. By Lemma~\ref{lem:irr.Sp.aba}(ii) every position $b>p-1+pc$ not on runner $p-1$ is unoccupied. 
Let $s$ be the number of beads on runner $p-2$ below position $p-2+p(a-1)$. Then $1\leq s\leq c-a+1$. 

{\sf{Case 1:}} $s\leq c-a$. An example of such a configuration is:

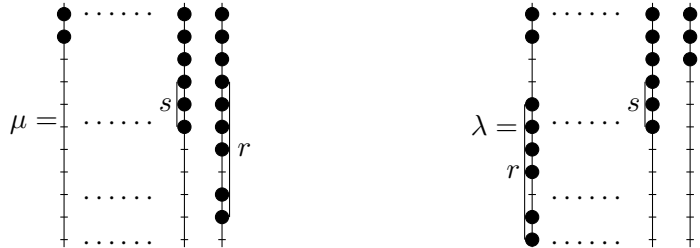
\vspace{2mm}

\hspace*{\fill}
{\begin{tikzpicture}
\draw(-.4,-1.625)node{$\mu=$};

 \draw (0,0)--(0,-3.25);
\draw (1.6,0)--(1.6,-3.25);
\draw(2.1,0)--(2.1,-3.25);

\draw(-.05,-.15)--(.05,-.15);
\draw(-.05,-.45)--(.05,-.45);
\draw(-.05,-.75)--(.05,-.75);
\draw(-.05,-1.05)--(.05,-1.05);
\draw(-.05,-1.35)--(.05,-1.35);
\draw(-.05,-1.65)--(.05,-1.65);
\draw(-.05,-1.95)--(.05,-1.95);
\draw(-.05,-2.25)--(.05,-2.25);
\draw(-.05,-2.55)--(.05,-2.55);
\draw(-.05,-2.85)--(.05,-2.85);
\draw(-.05,-3.15)--(.05,-3.15);

\draw(1.55,-.15)--(1.65,-.15);
\draw(1.55,-.45)--(1.65,-.45);
\draw(1.55,-.75)--(1.65,-.75);
\draw(1.55,-1.05)--(1.65,-1.05);
\draw(1.55,-1.35)--(1.65,-1.35);
\draw(1.55,-1.65)--(1.65,-1.65);
\draw(1.55,-1.95)--(1.65,-1.95);
\draw(1.55,-2.25)--(1.65,-2.25);
\draw(1.55,-2.55)--(1.65,-2.55);
\draw(1.55,-2.85)--(1.65,-2.85);
\draw(1.55,-3.15)--(1.65,-3.15);

\draw(2.05,-.15)--(2.15,-.15);
\draw(2.05,-.45)--(2.15,-.45);
\draw(2.05,-.75)--(2.15,-.75);
\draw(2.05,-1.05)--(2.15,-1.05);
\draw(2.05,-1.35)--(2.15,-1.35);
\draw(2.05,-1.65)--(2.15,-1.65);
\draw(2.05,-1.95)--(2.15,-1.95);
\draw(2.05,-2.25)--(2.15,-2.25);
\draw(2.05,-2.55)--(2.15,-2.55);
\draw(2.05,-2.85)--(2.15,-2.85);
\draw(2.05,-3.15)--(2.15,-3.15);

\draw(0.5,-.15)node{$\cdots$};
\draw(1,-.15)node{$\cdots$};
\draw(0.5,-1.6)node{$\cdots$};
\draw(1,-1.6)node{$\cdots$};
\draw(0.5,-2.6)node{$\cdots$};
\draw(1,-2.6)node{$\cdots$};
\draw(.5,-3.2)node{$\cdots$};
\draw(1,-3.2)node{$\cdots$};

\filldraw [black] (0,-.15) circle (2.5pt);
\filldraw [black] (0,-.45) circle (2.5pt);

\filldraw [black] (1.6,-.15) circle (2.5pt);
\filldraw [black] (1.6,-.45) circle (2.5pt);
\filldraw [black] (1.6,-.75) circle (2.5pt);
\filldraw [black] (1.6,-1.05) circle (2.5pt);
\filldraw [black] (1.6,-1.35)circle (2.5pt);
\filldraw [black] (1.6,-1.65) circle (2.5pt);

\filldraw [black] (2.1,-.15) circle (2.5pt);
\filldraw [black] (2.1,-.45) circle (2.5pt);
\filldraw [black] (2.1,-.75) circle (2.5pt);
\filldraw [black] (2.1,-1.05) circle (2.5pt);
\filldraw [black] (2.1,-1.35)circle (2.5pt);
\filldraw [black] (2.1,-1.65) circle (2.5pt);
\filldraw [black] (2.1,-1.95) circle (2.5pt);
\filldraw [black] (2.1,-2.55) circle (2.5pt);
\filldraw [black] (2.1,-2.85) circle (2.5pt);

\draw(1.5,-1.05)--(1.5,-1.65);
\draw(1.35,-1.35)node{$s$};
\draw(2.2,-1.05)--(2.2,-2.85);
\draw(2.4,-1.95)node{$r$};
\end{tikzpicture}}
\hfill 
{\begin{tikzpicture}
\draw(-.5,-1.625)node{$\lambda=$};

 \draw (0,0)--(0,-3.25);
\draw (1.6,0)--(1.6,-3.25);
\draw(2.1,0)--(2.1,-3.25);

\draw(-.05,-.15)--(.05,-.15);
\draw(-.05,-.45)--(.05,-.45);
\draw(-.05,-.75)--(.05,-.75);
\draw(-.05,-1.05)--(.05,-1.05);
\draw(-.05,-1.35)--(.05,-1.35);
\draw(-.05,-1.65)--(.05,-1.65);
\draw(-.05,-1.95)--(.05,-1.95);
\draw(-.05,-2.25)--(.05,-2.25);
\draw(-.05,-2.55)--(.05,-2.55);
\draw(-.05,-2.85)--(.05,-2.85);
\draw(-.05,-3.15)--(.05,-3.15);

\draw(1.55,-.15)--(1.65,-.15);
\draw(1.55,-.45)--(1.65,-.45);
\draw(1.55,-.75)--(1.65,-.75);
\draw(1.55,-1.05)--(1.65,-1.05);
\draw(1.55,-1.35)--(1.65,-1.35);
\draw(1.55,-1.65)--(1.65,-1.65);
\draw(1.55,-1.95)--(1.65,-1.95);
\draw(1.55,-2.25)--(1.65,-2.25);
\draw(1.55,-2.55)--(1.65,-2.55);
\draw(1.55,-2.85)--(1.65,-2.85);
\draw(1.55,-3.15)--(1.65,-3.15);

\draw(2.05,-.15)--(2.15,-.15);
\draw(2.05,-.45)--(2.15,-.45);
\draw(2.05,-.75)--(2.15,-.75);
\draw(2.05,-1.05)--(2.15,-1.05);
\draw(2.05,-1.35)--(2.15,-1.35);
\draw(2.05,-1.65)--(2.15,-1.65);
\draw(2.05,-1.95)--(2.15,-1.95);
\draw(2.05,-2.25)--(2.15,-2.25);
\draw(2.05,-2.55)--(2.15,-2.55);
\draw(2.05,-2.85)--(2.15,-2.85);
\draw(2.05,-3.15)--(2.15,-3.15);

\draw(0.5,-.15)node{$\cdots$};
\draw(1,-.15)node{$\cdots$};
\draw(0.5,-1.6)node{$\cdots$};
\draw(1,-1.6)node{$\cdots$};
\draw(0.5,-2.6)node{$\cdots$};
\draw(1,-2.6)node{$\cdots$};
\draw(.5,-3.2)node{$\cdots$};
\draw(1,-3.2)node{$\cdots$};

\filldraw [black] (0,-.15) circle (2.5pt);
\filldraw [black] (0,-.45) circle (2.5pt);
\filldraw [black] (0,-1.35)circle (2.5pt);
\filldraw [black] (0,-1.65) circle (2.5pt);
\filldraw [black] (0,-1.95) circle (2.5pt);
\filldraw [black] (0,-2.25) circle (2.5pt);
\filldraw [black] (0,-2.85) circle (2.5pt);
\filldraw [black] (0,-3.15) circle (2.5pt);

\filldraw [black] (1.6,-.15) circle (2.5pt);
\filldraw [black] (1.6,-.45) circle (2.5pt);
\filldraw [black] (1.6,-.75) circle (2.5pt);
\filldraw [black] (1.6,-1.05) circle (2.5pt);
\filldraw [black] (1.6,-1.35)circle (2.5pt);
\filldraw [black] (1.6,-1.65) circle (2.5pt);

\filldraw [black] (2.1,-.15) circle (2.5pt);
\filldraw [black] (2.1,-.45) circle (2.5pt);
\filldraw [black] (2.1,-.75) circle (2.5pt);

\draw(1.5,-1.05)--(1.5,-1.65);
\draw(1.35,-1.35)node{$s$};
\draw(-.1,-1.35)--(-.1,-3.15);
\draw(-.25,-2.25)node{$r$};
\end{tikzpicture}}
\hspace*{\fill}

Note that $\phi_{i-1}(\la)=s$ and $\beps_{i-1}(\la)=0$. Setting $\xi:=\tf^s_{i-1}\la$ we have by Lemma~\ref{L200509}: 
\begin{equation}\label{E060520}
\Ext^1_{\Si_n}(D^\la,D^\la)\cong\Ext^1_{\Si_{n+s}}(D^\xi,D^\xi).
\end{equation}  
Note that $\phi_i(\xi)=\phi-r$, and let $\rho:=\tf_i^{\phi-r}\xi$. 
Note that $\eps_i(\rho)=\phi-s$ and if $A_1,\dots,A_{\phi-s}$ are the $i$-normal nodes of $\rho$ labeled from bottom to top, one easily sees that $\rho_{A_{\phi-r+1}}$ is $p$-regular. So by  Theorem~\ref{prop:sKMT.res} we have $\Hom_{\Si_{n+s}}(D^\xi, (e_i^{(\phi-r)}D^\rho)/D^\xi)=0$. So Lemma~\ref{LES} yields an embedding 
\begin{equation}\label{E060520_2}
\Ext^1_{\Si_{n+s}}(D^\xi,D^\xi)\into\Ext^1_{\Si_{n+s+\phi-r}}(D^\rho,D^\rho).
\end{equation}

Note that $\eps_{i-1}(\rho)=s$. For $\si:=\te_{i-1}^{s}\rho$, we have that $\si=\tf_i^{\phi}\mu$. By Lemma~\ref{lem:rem.ad1} $\phi=\bphi_i(\mu)$ and so by Lemma~\ref{lem:res.Sp} we get that 
 $D^\si\cong f_i^{(\phi)}D^\mu\cong f_i^{(\phi)}S^\mu\cong S^\si,$
However, the $(i-1)$-addable nodes of $\si$ do not correspond to addable positions on the non-restricted runner $0$ of $\si$. Therefore, by Lemma~\ref{lem:reg2} applied to $\si$ instead of $\la$ we get $\Ext^1_{\Si_{n+s+\phi-r}}(D^\rho, D^\rho)=0$, and so 
 $\Ext^1_{\Si_n}(D^\la,D^\la)=0$ by (\ref{E060520}) and (\ref{E060520_2}).

\vspace{1mm}
{\sf{Case 2:}} $s=c-a+1$. An example of such a configuration is:

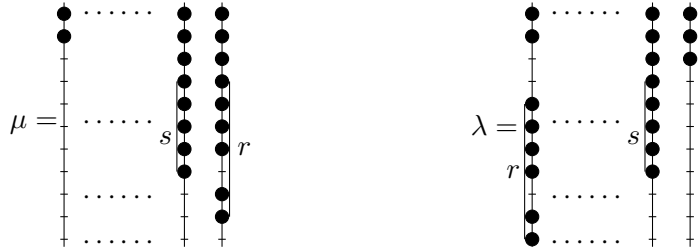
\vspace{2mm}

\hspace*{\fill}
{\begin{tikzpicture}
\draw(-.4,-1.625)node{$\mu=$};

 \draw (0,0)--(0,-3.25);
\draw (1.6,0)--(1.6,-3.25);
\draw(2.1,0)--(2.1,-3.25);

\draw(-.05,-.15)--(.05,-.15);
\draw(-.05,-.45)--(.05,-.45);
\draw(-.05,-.75)--(.05,-.75);
\draw(-.05,-1.05)--(.05,-1.05);
\draw(-.05,-1.35)--(.05,-1.35);
\draw(-.05,-1.65)--(.05,-1.65);
\draw(-.05,-1.95)--(.05,-1.95);
\draw(-.05,-2.25)--(.05,-2.25);
\draw(-.05,-2.55)--(.05,-2.55);
\draw(-.05,-2.85)--(.05,-2.85);
\draw(-.05,-3.15)--(.05,-3.15);

\draw(1.55,-.15)--(1.65,-.15);
\draw(1.55,-.45)--(1.65,-.45);
\draw(1.55,-.75)--(1.65,-.75);
\draw(1.55,-1.05)--(1.65,-1.05);
\draw(1.55,-1.35)--(1.65,-1.35);
\draw(1.55,-1.65)--(1.65,-1.65);
\draw(1.55,-1.95)--(1.65,-1.95);
\draw(1.55,-2.25)--(1.65,-2.25);
\draw(1.55,-2.55)--(1.65,-2.55);
\draw(1.55,-2.85)--(1.65,-2.85);
\draw(1.55,-3.15)--(1.65,-3.15);

\draw(2.05,-.15)--(2.15,-.15);
\draw(2.05,-.45)--(2.15,-.45);
\draw(2.05,-.75)--(2.15,-.75);
\draw(2.05,-1.05)--(2.15,-1.05);
\draw(2.05,-1.35)--(2.15,-1.35);
\draw(2.05,-1.65)--(2.15,-1.65);
\draw(2.05,-1.95)--(2.15,-1.95);
\draw(2.05,-2.25)--(2.15,-2.25);
\draw(2.05,-2.55)--(2.15,-2.55);
\draw(2.05,-2.85)--(2.15,-2.85);
\draw(2.05,-3.15)--(2.15,-3.15);

\draw(0.5,-.15)node{$\cdots$};
\draw(1,-.15)node{$\cdots$};
\draw(0.5,-1.6)node{$\cdots$};
\draw(1,-1.6)node{$\cdots$};
\draw(0.5,-2.6)node{$\cdots$};
\draw(1,-2.6)node{$\cdots$};
\draw(.5,-3.2)node{$\cdots$};
\draw(1,-3.2)node{$\cdots$};

\filldraw [black] (0,-.15) circle (2.5pt);
\filldraw [black] (0,-.45) circle (2.5pt);

\filldraw [black] (1.6,-.15) circle (2.5pt);
\filldraw [black] (1.6,-.45) circle (2.5pt);
\filldraw [black] (1.6,-.75) circle (2.5pt);
\filldraw [black] (1.6,-1.05) circle (2.5pt);
\filldraw [black] (1.6,-1.35)circle (2.5pt);
\filldraw [black] (1.6,-1.65) circle (2.5pt);
\filldraw [black] (1.6,-1.95)circle (2.5pt);
\filldraw [black] (1.6,-2.25) circle (2.5pt);

\filldraw [black] (2.1,-.15) circle (2.5pt);
\filldraw [black] (2.1,-.45) circle (2.5pt);
\filldraw [black] (2.1,-.75) circle (2.5pt);
\filldraw [black] (2.1,-1.05) circle (2.5pt);
\filldraw [black] (2.1,-1.35)circle (2.5pt);
\filldraw [black] (2.1,-1.65) circle (2.5pt);
\filldraw [black] (2.1,-1.95) circle (2.5pt);
\filldraw [black] (2.1,-2.55) circle (2.5pt);
\filldraw [black] (2.1,-2.85) circle (2.5pt);

\draw(1.5,-1.05)--(1.5,-2.25);
\draw(1.35,-1.8)node{$s$};
\draw(2.2,-1.05)--(2.2,-2.85);
\draw(2.4,-1.95)node{$r$};
\end{tikzpicture}}
\hfill 
{\begin{tikzpicture}
\draw(-.5,-1.625)node{$\lambda=$};

 \draw (0,0)--(0,-3.25);
\draw (1.6,0)--(1.6,-3.25);
\draw(2.1,0)--(2.1,-3.25);

\draw(-.05,-.15)--(.05,-.15);
\draw(-.05,-.45)--(.05,-.45);
\draw(-.05,-.75)--(.05,-.75);
\draw(-.05,-1.05)--(.05,-1.05);
\draw(-.05,-1.35)--(.05,-1.35);
\draw(-.05,-1.65)--(.05,-1.65);
\draw(-.05,-1.95)--(.05,-1.95);
\draw(-.05,-2.25)--(.05,-2.25);
\draw(-.05,-2.55)--(.05,-2.55);
\draw(-.05,-2.85)--(.05,-2.85);
\draw(-.05,-3.15)--(.05,-3.15);

\draw(1.55,-.15)--(1.65,-.15);
\draw(1.55,-.45)--(1.65,-.45);
\draw(1.55,-.75)--(1.65,-.75);
\draw(1.55,-1.05)--(1.65,-1.05);
\draw(1.55,-1.35)--(1.65,-1.35);
\draw(1.55,-1.65)--(1.65,-1.65);
\draw(1.55,-1.95)--(1.65,-1.95);
\draw(1.55,-2.25)--(1.65,-2.25);
\draw(1.55,-2.55)--(1.65,-2.55);
\draw(1.55,-2.85)--(1.65,-2.85);
\draw(1.55,-3.15)--(1.65,-3.15);

\draw(2.05,-.15)--(2.15,-.15);
\draw(2.05,-.45)--(2.15,-.45);
\draw(2.05,-.75)--(2.15,-.75);
\draw(2.05,-1.05)--(2.15,-1.05);
\draw(2.05,-1.35)--(2.15,-1.35);
\draw(2.05,-1.65)--(2.15,-1.65);
\draw(2.05,-1.95)--(2.15,-1.95);
\draw(2.05,-2.25)--(2.15,-2.25);
\draw(2.05,-2.55)--(2.15,-2.55);
\draw(2.05,-2.85)--(2.15,-2.85);
\draw(2.05,-3.15)--(2.15,-3.15);

\draw(0.5,-.15)node{$\cdots$};
\draw(1,-.15)node{$\cdots$};
\draw(0.5,-1.6)node{$\cdots$};
\draw(1,-1.6)node{$\cdots$};
\draw(0.5,-2.6)node{$\cdots$};
\draw(1,-2.6)node{$\cdots$};
\draw(.5,-3.2)node{$\cdots$};
\draw(1,-3.2)node{$\cdots$};

\filldraw [black] (0,-.15) circle (2.5pt);
\filldraw [black] (0,-.45) circle (2.5pt);
\filldraw [black] (0,-1.35)circle (2.5pt);
\filldraw [black] (0,-1.65) circle (2.5pt);
\filldraw [black] (0,-1.95) circle (2.5pt);
\filldraw [black] (0,-2.25) circle (2.5pt);
\filldraw [black] (0,-2.85) circle (2.5pt);
\filldraw [black] (0,-3.15) circle (2.5pt);

\filldraw [black] (1.6,-.15) circle (2.5pt);
\filldraw [black] (1.6,-.45) circle (2.5pt);
\filldraw [black] (1.6,-.75) circle (2.5pt);
\filldraw [black] (1.6,-1.05) circle (2.5pt);
\filldraw [black] (1.6,-1.35)circle (2.5pt);
\filldraw [black] (1.6,-1.65) circle (2.5pt);
\filldraw [black] (1.6,-1.95)circle (2.5pt);
\filldraw [black] (1.6,-2.25) circle (2.5pt);

\filldraw [black] (2.1,-.15) circle (2.5pt);
\filldraw [black] (2.1,-.45) circle (2.5pt);
\filldraw [black] (2.1,-.75) circle (2.5pt);

\draw(1.5,-1.05)--(1.5,-2.25);
\draw(1.35,-1.8)node{$s$};
\draw(-.1,-1.35)--(-.1,-3.15);
\draw(-.25,-2.25)node{$r$};
\end{tikzpicture}}
\hspace*{\fill}

As in Case 1, $\phi_{i-1}(\la)=s$, $\beps_{i-1}(\la)=0$ and so 
we have (\ref{E060520}) for $\xi:=\tf^s_{i-1}\la$. 
Note that $\phi_i(\xi)=\phi-r+1$. Setting $\rho:=\tf_i^{\phi-r+1}\xi$ we obtain, as in Case 1, an embedding (\ref{E060520_2}). 

Now $\eps_{i-1}(\rho)=s-1$. Let $\tau:= \te_{i-1}^{(s-1)}\rho$. We claim that $S^\tau=D^\tau$. To see this note that $\Ga(\tau)$ is obtained from $\Ga(\tf_i^\phi\mu)$  by removing the bead from position $p-2+pc$ and placing it to position $p(c+1)$. We have that $\tau^{(0)}$ is the partition obtained by removing the first column of $(\tf_i^{\phi}\mu)^{(0)}$. It follows that $\Ga(\tau)$ satisfies the properties of Lemma~\ref{lem:irr.Sp.aba} and so $S^\tau=D^\tau$.
 
 Now if $\tau$ is not $p$-restricted, then since the $(i-1)$-addable nodes of $\tau$ do not correspond to addable positions on the non-restricted runner $0$,  by Lemma~\ref{lem:reg2} we obtain that  $\Ext^1_{\Si_{n+s+\phi-r}}(D^\rho, D^\rho)=0$ and so $\Ext^1_{\Si_n}(D^\la,D^\la)=0$. If $\tau$ is $p$-restricted, and so a core, the result follows by Proposition~\ref{prop:Sp.Core}. 
\end{proof}

\subsection{The case $\nu$ is $p$-restricted but not $p$-regular}
\label{subsec:Sp.res}
In this case we have $\mu=\nu^\reg\neq \nu$. 

\begin{Proposition}
\label{prop:Sp.res}
If $\nu$ be $p$-restricted and $p$-singular then $\Ext^1_{\Si_n}(D^\la,D^\la)=0.$
\end{Proposition}

\begin{proof}
By Lemma~\ref{lem:irr.Sp.mul} we have $e_{-i}^{(\eps_i(\la))}D^\laM\cong D^\muM\cong S^{\nu'}$. Since $\nu$ is $p$-restricted, $\nu'$ is $p$-regular, so $\nu'=\muM$. The result now follows directly by Proposition~\ref{prop:Sp.reg}, applied to $\laM$ instead of $\la$,  and (\ref{ESignExt}). 
\end{proof}

\subsection{The case $\nu$ is not $p$-restricted and not $p$-regular}
\label{subsec:Sp.gener}
In this case an abacus display $\Ga(\nu)$ has a non-restricted runner $j$ and a non-regular runner $k$ as in  Lemma~\ref{lem:irr.Sp.aba}. 
Note that $\mu=\nu^\reg$, and let $\Ga(\mu)=\Ga(\nu)^\reg$ be as in Lemma~\ref{lem:reg1}. Suppose that the $i$-addable nodes of $\mu$ correspond to addable positions on the runner $m$ of $\Ga(\mu)$. 
With this notation we have:

\begin{Lemma}
\label{lem:gener1}
We have $\Ext^1_{\Si_n}(D^\la,D^\la)=0$ unless $k=m=j+1$.
\end{Lemma}

\begin{proof}
By Lemma~\ref{lem:reg1} $\mu^{(k)}=\varnothing$ and if the last bead on runner $k$ occurs at position $k+pa$, then every position $b<k+pa$ is occupied. Therefore the addable nodes of $\mu$ occur at unoccupied positions $b>k+pa$. Now if $m\neq k$, then $\mu^B$ is $p$-regular for every $i$-addable node $B$ of $\mu$. In this case Corollary~\ref{CRepeated} gives $\Ext^1_{\Si_n}(D^\la,D^\la)=0$. So we may assume that $m=k$. 

Suppose $k\neq j+1$. 
 By Lemma~\ref{lem:irr.Sp.mul}, 
$e_{-i}^{(r)}D^{\laM}\cong D^{\muM}\cong S^{\nu'}$, with $\nu'$ is not $p$-regular and not $p$-restricted. Note that in $\Ga(\nu')=\Ga(\nu)'$, the runner $k':=p-1-j$ is the non-regular runner and $j':=p-1-k$ is the non-restricted runner. Moreover the $(-i)$-addable nodes of $\nu'$ are on runner $m'=p-1-(m-1)=p-1-(k-1)\neq p-1-j=k'$. By the previous paragraph, $\Ext^1_{\Si_n}(D^\laM,D^\laM)=0$, and in view of (\ref{ESignExt}) we are done. \end{proof}

\begin{Remark}
\label{rem:gener2}
From now on we choose $\Ga(\nu)$ so that $j=p-1$. In view of Lemma~\ref{lem:gener1}, we may assume that $k=m=0$. 
Let $p-1+pd$ be the first unoccupied position on the non-restricted runner $p-1$ in $\Ga(\nu)$. We may assume the following additional property:

\vspace{2mm}
\noindent
{\sf (P)} Position $p-2+pd$ in $\Ga(\nu)$ is unoccupied.

\vspace{ 2mm}
Indeed, suppose that position $p-2+pd$ in $\Ga(\nu)$ is occupied. We show that there exists $\ka\in\Parreg_t$ for some $t$ such that 
\begin{equation}\label{E070520}
\Ext^1_{\Si_n}(D^\la,D^\la)\cong \Ext^1_{\Si_t}(D^\ka,D^\ka),
\end{equation} 
all the assumptions on $\la$ hold for $\ka$, and the corresponding property {\sf (P)} holds for $\ka$.

Assume first that there exists a runner $l$ in $\Ga(\nu)$ with $1\leq l<p-2$ such that position $l+pd$ is unoccupied. We may assume that $l$ is maximal with this property.  Set $\ell:=i-(p-1-l) \pmod{p}$. By Lemma~\ref{lem:irr.Sp.aba}, we have $s:=\beps_{\ell}(\nu)>0$ and $\bphi_{\ell}(\nu)=0$. By Lemma~\ref{lem:rem.ad1}, we have then that $\eps_{\ell}(\mu)=s$ and $\phi_{\ell}(\mu)=0$. Since $\ell\neq i-1,i,i+1,$ we have that $\eps_{\ell}(\la)=s$ and $\phi_{\ell}(\la)=0$. Then $\Ext^1_{\Si_n}(D^\la,D^\la)\cong\Ext^1_{\Si_{n-s}}(D^{\te_\ell^s\la},D^{\te_\ell^s\la})$ by Lemma~\ref{L200509}. 
Moreover, $e_i^{(r)}D^{\te_\ell^s\la}\cong S^{\he_\ell^s\nu}$. In fact, $\Ga(\he_\ell^s\nu)$ is obtained by swapping the runners $l$ and $l+1$ of $\Ga(\nu)$. 
Repeat this process now for $\te_\ell^s\la$ and the runner $l+1$. Eventually we get (\ref{E070520}) for $t=n-s(p-2-l)$ for $\ka$ with 
$\eps_i(\ka)=r$, $e_i^{(r)}D^{\kappa}\cong S^\tau$ and  $\Ga(\tau)$ obtained from $\Ga(\nu)$ by swapping the runners $l$ and $p-2$.

Now assume that positions $l+pd$ are occupied in $\Ga(\nu)$ for all $1\leq l\leq  p-2$. If the last bead on runner $0$ appears at position $pb$ then position $1+pb$ is also occupied. By  Lemma~\ref{lem:irr.Sp.mul}, we have  $e_{-i}^{(r)}D^{\laM}\cong D^{\muM}\cong S^{\nu'}.$ 
Runner $0$ in $\Ga(\nu')=\Ga(\nu)'$ is its non-regular runner and runner $p-1$ is  its non-restricted runner. Let $p-1+pc$ be the first unoccupied position on runner $p-1$ of $\Ga(\nu')$. Since position $1+pb$ is occupied in $\Ga(\nu)$, position $p-2+pc$ is unoccupied in $\Ga(\nu')$, and so by (\ref{ESignExt}), we have (\ref{E070520}) with $\ka=\laM$. 
\end{Remark}

 \begin{Proposition}
\label{prop:Sp.gen}
If $\nu$ is neither $p$-regular nor $p$-restricted then $\Ext^1_{\Si_n}(D^\la,D^\la)=0.$
\end{Proposition}

\begin{proof}
As $D^\mu\cong S^\nu$, we have $\mu=\nu^\reg$. Let $\phi:=\phi_i(\mu)$. In view of Lemma~\ref{lem:shap}, we may assume that $\phi>r\geq 1$. By Lemma~\ref{lem:rem.ad1}, we have $\bphi_i(\nu)=\phi$.  

By Remark~\ref{rem:gener2}, we may assume that $k=m=0$, $j=p-1$, and if $p-1+pd$ is the first unoccupied position on runner $p-1$ then position $p-2+pd$ is also unoccupied. 
As $m=0$, an abacus display $\Ga(\hf_i^r\nu)$ is obtained from $\Ga(\nu)$ by sliding $r$ beads, call them $a_1,\dots,a_r$, from runner $p-1$ to $0$. Let $(p-1)+pc_t$ be the position of $a_t$. We assume that $c:=c_1<\dots<c_r$. Let also $p-2+pb$ be the position of the last bead on runner $p-2$.

\vspace{1mm}
{\sf{Case 1:}} $b\geq c$. Let $s:=\max\{t\mid c_t\leq b\}$.  
Note that $1\leq s\leq b-c+1$. An example of such a configuration is:   

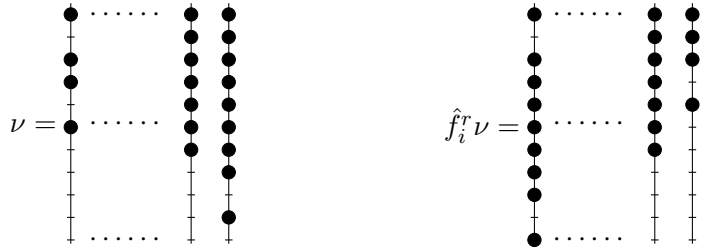
\vspace{1mm}

\hspace*{\fill}
{\begin{tikzpicture}
\draw(-.5,-1.625)node{$\nu=$};

 \draw (0,0)--(0,-3.2);
\draw (1.6,0)--(1.6,-3.2);
\draw(2.1,0)--(2.1,-3.2);

\draw(-.05,-.15)--(.05,-.15);
\draw(-.05,-.45)--(.05,-.45);
\draw(-.05,-.75)--(.05,-.75);
\draw(-.05,-1.05)--(.05,-1.05);
\draw(-.05,-1.35)--(.05,-1.35);
\draw(-.05,-1.65)--(.05,-1.65);
\draw(-.05,-1.95)--(.05,-1.95);
\draw(-.05,-2.25)--(.05,-2.25);
\draw(-.05,-2.55)--(.05,-2.55);
\draw(-.05,-2.85)--(.05,-2.85);
\draw(-.05,-3.15)--(.05,-3.15);

\draw(1.55,-.15)--(1.65,-.15);
\draw(1.55,-.45)--(1.65,-.45);
\draw(1.55,-.75)--(1.65,-.75);
\draw(1.55,-1.05)--(1.65,-1.05);
\draw(1.55,-1.35)--(1.65,-1.35);
\draw(1.55,-1.65)--(1.65,-1.65);
\draw(1.55,-1.95)--(1.65,-1.95);
\draw(1.55,-2.25)--(1.65,-2.25);
\draw(1.55,-2.55)--(1.65,-2.55);
\draw(1.55,-2.85)--(1.65,-2.85);
\draw(1.55,-3.15)--(1.65,-3.15);

\draw(2.05,-.15)--(2.15,-.15);
\draw(2.05,-.45)--(2.15,-.45);
\draw(2.05,-.75)--(2.15,-.75);
\draw(2.05,-1.05)--(2.15,-1.05);
\draw(2.05,-1.35)--(2.15,-1.35);
\draw(2.05,-1.65)--(2.15,-1.65);
\draw(2.05,-1.95)--(2.15,-1.95);
\draw(2.05,-2.25)--(2.15,-2.25);
\draw(2.05,-2.55)--(2.15,-2.55);
\draw(2.05,-2.85)--(2.15,-2.85);
\draw(2.05,-3.15)--(2.15,-3.15);

\draw(0.5,-.15)node{$\cdots$};
\draw(1,-.15)node{$\cdots$};
\draw(0.5,-1.6)node{$\cdots$};
\draw(1,-1.6)node{$\cdots$};
\draw(.5,-3.15)node{$\cdots$};
\draw(1,-3.15)node{$\cdots$};

\filldraw [black] (0,-.15) circle (2.5pt);
\filldraw [black] (0,-.75) circle (2.5pt);
\filldraw [black] (0,-1.05) circle (2.5pt);
\filldraw [black] (0,-1.65) circle (2.5pt);

\filldraw [black] (1.6,-.15) circle (2.5pt);
\filldraw [black] (1.6,-.45) circle (2.5pt);
\filldraw [black] (1.6,-.75) circle (2.5pt);
\filldraw [black] (1.6,-1.05) circle (2.5pt);
\filldraw [black] (1.6,-1.35)circle (2.5pt);
\filldraw [black] (1.6,-1.65) circle (2.5pt);
\filldraw [black] (1.6,-1.95) circle (2.5pt);

\filldraw [black] (2.1,-.15) circle (2.5pt);
\filldraw [black] (2.1,-.45) circle (2.5pt);
\filldraw [black] (2.1,-.75) circle (2.5pt);
\filldraw [black] (2.1,-1.05) circle (2.5pt);
\filldraw [black] (2.1,-1.35)circle (2.5pt);
\filldraw [black] (2.1,-1.65) circle (2.5pt);
\filldraw [black] (2.1,-1.95) circle (2.5pt);
\filldraw [black] (2.1,-2.25) circle (2.5pt);
\filldraw [black] (2.1,-2.85) circle (2.5pt);

\end{tikzpicture}}
\hfill 
{\begin{tikzpicture}
\draw(-.7,-1.625)node{$\hf_i^r\nu=$};

 \draw (0,0)--(0,-3.2);
\draw (1.6,0)--(1.6,-3.2);
\draw(2.1,0)--(2.1,-3.2);

\draw(-.05,-.15)--(.05,-.15);
\draw(-.05,-.45)--(.05,-.45);
\draw(-.05,-.75)--(.05,-.75);
\draw(-.05,-1.05)--(.05,-1.05);
\draw(-.05,-1.35)--(.05,-1.35);
\draw(-.05,-1.65)--(.05,-1.65);
\draw(-.05,-1.95)--(.05,-1.95);
\draw(-.05,-2.25)--(.05,-2.25);
\draw(-.05,-2.55)--(.05,-2.55);
\draw(-.05,-2.85)--(.05,-2.85);
\draw(-.05,-3.15)--(.05,-3.15);

\draw(1.55,-.15)--(1.65,-.15);
\draw(1.55,-.45)--(1.65,-.45);
\draw(1.55,-.75)--(1.65,-.75);
\draw(1.55,-1.05)--(1.65,-1.05);
\draw(1.55,-1.35)--(1.65,-1.35);
\draw(1.55,-1.65)--(1.65,-1.65);
\draw(1.55,-1.95)--(1.65,-1.95);
\draw(1.55,-2.25)--(1.65,-2.25);
\draw(1.55,-2.55)--(1.65,-2.55);
\draw(1.55,-2.85)--(1.65,-2.85);
\draw(1.55,-3.15)--(1.65,-3.15);

\draw(2.05,-.15)--(2.15,-.15);
\draw(2.05,-.45)--(2.15,-.45);
\draw(2.05,-.75)--(2.15,-.75);
\draw(2.05,-1.05)--(2.15,-1.05);
\draw(2.05,-1.35)--(2.15,-1.35);
\draw(2.05,-1.65)--(2.15,-1.65);
\draw(2.05,-1.95)--(2.15,-1.95);
\draw(2.05,-2.25)--(2.15,-2.25);
\draw(2.05,-2.55)--(2.15,-2.55);
\draw(2.05,-2.85)--(2.15,-2.85);
\draw(2.05,-3.15)--(2.15,-3.15);

\draw(0.5,-.15)node{$\cdots$};
\draw(1,-.15)node{$\cdots$};
\draw(0.5,-1.6)node{$\cdots$};
\draw(1,-1.6)node{$\cdots$};
\draw(.5,-3.15)node{$\cdots$};
\draw(1,-3.15)node{$\cdots$};

\filldraw [black] (0,-.15) circle (2.5pt);
\filldraw [black] (0,-.75) circle (2.5pt);
\filldraw [black] (0,-1.05) circle (2.5pt);
\filldraw [black] (0,-1.35) circle (2.5pt);
\filldraw [black] (0,-1.65) circle (2.5pt);
\filldraw [black] (0,-1.95) circle (2.5pt);
\filldraw [black] (0,-2.25) circle (2.5pt);
\filldraw [black] (0,-2.55) circle (2.5pt);
\filldraw [black] (0,-3.15) circle (2.5pt);

\filldraw [black] (1.6,-.15) circle (2.5pt);
\filldraw [black] (1.6,-.45) circle (2.5pt);
\filldraw [black] (1.6,-.75) circle (2.5pt);
\filldraw [black] (1.6,-1.05) circle (2.5pt);
\filldraw [black] (1.6,-1.35)circle (2.5pt);
\filldraw [black] (1.6,-1.65) circle (2.5pt);
\filldraw [black] (1.6,-1.95) circle (2.5pt);

\filldraw [black] (2.1,-.15) circle (2.5pt);
\filldraw [black] (2.1,-.45) circle (2.5pt);
\filldraw [black] (2.1,-.75) circle (2.5pt);
\filldraw [black] (2.1,-1.35)circle (2.5pt);

\end{tikzpicture}}
\hspace*{\fill}


Set $\nu^1:=\hf_i^r\nu$. By Lemma~\ref{lem:irr.Sp.dfilt} we have that $\head S^{\nu^1}\cong D^\la$ and $[S^{\nu^1}:D^\la]=1$.
Now $\beps_{i-1}(\nu^1)=0$ and so $\eps_{i-1}(\la)=0$ by Lemma~\ref{lem:head.nodes}(i). Moreover $\phi_{i-1}(\nu^1)=\bphi_{i-1}(\nu^1)=s$ and so $\phi_{i-1}(\la)=s$, by Lemma~\ref{lem:head.nodes}(ii). Let $\xi:=\tf_{i-1}^{(s)}\la$. Then  $D^\xi\cong f^{(s)}_{i-1}D^\la$ and $D^\la\cong e^{(s)}_{i-1}D^\xi$. 

Set $\nu^2:=\hf_{i-1}^s\nu^1$. It follows using Lemma~\ref{lem:shap} that $\head S^{\nu^2}\cong D^\xi$, $[S^{\nu^2}:D^\xi]=1$ and
\begin{equation}\label{E200508}
\Ext^1_{\Si_n}(D^\la,D^\la)\cong\Ext^1_{\Si_{n+s}}(D^\xi,D^\xi).
\end{equation}
Note that $\phi_{i}(\nu^2)=\bphi_i(\nu^2)=\phi-r$, so $\phi_{i}(\xi)=\phi-r$ by Lemma~\ref{lem:head.nodes}(ii).
Let $\rho:=\tf_i^{\phi-r}\xi$ and $\nu^3:=\hf^{\phi-r}_i\nu^2$. Then $D^\rho$ is in the head of $S^{\nu^3}$. 

By Lemma~\ref{lem:res.Sp} we have that $e_i^{(\phi-r)}S^{\nu^3}$ has a Specht filtration with $S^{\nu^2}$ being the top Specht factor. Moreover, we have $\bphi_{i-1}(\zeta)>0$ for any other Specht factor $S^\zeta$ of the filtration. 
If $B$ is the lowest $(i-1)$-addable node of such $\zeta$, observe that  $B$ is conormal for $\zeta$. Since $\phi_{i-1}(\xi)=0$, we now deduce from Lemma~\ref{lem:ad.node.trick} that $\Hom_{\Si_{n+s}}(S^\zeta,D^\xi)=0$. We have $e_i^{(\phi-r)}S^{\nu^3}\sim X|D^\xi$, with $X\sim Y|\rad S^{\nu^2}$ and $Y$ having a Specht filtration with factors $S^\zeta$ as above. 
We have $\Hom_{\Si_{n+s}}(\rad S^{\nu^2},D^\xi)=0$, since $D^\xi$ is the simple head of $S^{\nu^2}$ and  $[S^{\nu^2}:D^\xi]=1$. Hence, $\Hom_{\Si_{n+s}}(X,D^\xi)=0$. 
Since $D^\rho$ is in the head of $S^{\nu^3}$, we have that $\rad (e_i^{(\phi-r)}D^\rho)$ is a quotient of $X$, and hence 
$\Hom_{\Si_{n+s}}(\rad (e_i^{(\phi-r)}D^\rho),D^\xi)=0$. 
Applying $\Hom_{\Si_{n+s}}(-,D^\xi)$ to the short exact sequence 
$0\to \rad (e_i^{(\phi-r)}D^\rho)\to e_i^{(\phi-r)}D^\rho\to D^\xi\to 0$ and using Lemma~\ref{lem:shap}, 
we now obtain an embedding 
\begin{equation}\label{E200508_2}
\Ext^1_{\Si_{n+s}}(D^\xi,D^\xi)\into \Ext^1_{\Si_{n+s+\phi-r}}(D^\rho,D^\rho).
\end{equation}

Note that $\eps_{i-1}(\nu^3)=\beps_{i-1}(\nu^3)=s$ and so $\eps_{i-1}(\rho)=s$, by Lemma~\ref{lem:head.nodes}(ii). We set $\si:=\te_{i-1}^s\rho$ and $\nu^4:=\he^s_{i-1}\nu^3$. Then $D^\si$ is in the head of $S^{\nu^4}$. Note that $\nu^4=\hf_i^\phi\nu$,  so
$S^{\nu^4}\cong f_i^{(\phi)}S^\nu\cong f_i^{(\phi)}D^\mu\cong D^{\tf_i^\phi\mu}\cong D^\si$ is an irreducible Specht module. The non-regular runner of $\Ga(\nu^4)$ is $p-1$ and the non-restricted runner of $\Ga(\nu^4)$ is $0$. By Lemma~\ref{lem:gener1} we deduce that $\Ext^1_{\Si_{n+s+\phi-r}}(D^\rho,D^\rho)=0$ and so $\Ext^1_{\Si_n}(D^\la,D^\la)=0$ in view of (\ref{E200508}) and (\ref{E200508_2}). 

\vspace{1mm}

{\sf Case 2:} $b<c$. Let $s$ be the number of beads that occur on runner $p-1$ below the position $p-1+pb$. Hence we have that $1\leq r\leq s$. An example of such a configuration is:   

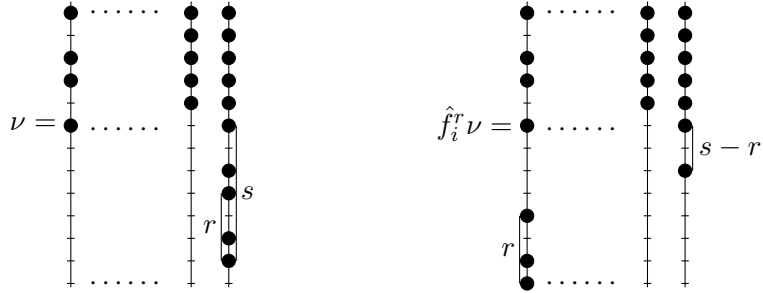
\vspace{2mm}

\hspace*{\fill}
{\begin{tikzpicture}
\draw(-.5,-1.625)node{$\nu=$};

 \draw (0,0)--(0,-3.8);
\draw (1.6,0)--(1.6,-3.8);
\draw(2.1,0)--(2.1,-3.8);

\draw(-.05,-.15)--(.05,-.15);
\draw(-.05,-.45)--(.05,-.45);
\draw(-.05,-.75)--(.05,-.75);
\draw(-.05,-1.05)--(.05,-1.05);
\draw(-.05,-1.35)--(.05,-1.35);
\draw(-.05,-1.65)--(.05,-1.65);
\draw(-.05,-1.95)--(.05,-1.95);
\draw(-.05,-2.25)--(.05,-2.25);
\draw(-.05,-2.55)--(.05,-2.55);
\draw(-.05,-2.85)--(.05,-2.85);
\draw(-.05,-3.15)--(.05,-3.15);
\draw(-.05,-3.45)--(.05,-3.45);
\draw(-.05,-3.75)--(.05,-3.75);

\draw(1.55,-.15)--(1.65,-.15);
\draw(1.55,-.45)--(1.65,-.45);
\draw(1.55,-.75)--(1.65,-.75);
\draw(1.55,-1.05)--(1.65,-1.05);
\draw(1.55,-1.35)--(1.65,-1.35);
\draw(1.55,-1.65)--(1.65,-1.65);
\draw(1.55,-1.95)--(1.65,-1.95);
\draw(1.55,-2.25)--(1.65,-2.25);
\draw(1.55,-2.55)--(1.65,-2.55);
\draw(1.55,-2.85)--(1.65,-2.85);
\draw(1.55,-3.15)--(1.65,-3.15);
\draw(1.55,-3.45)--(1.65,-3.45);
\draw(1.55,-3.75)--(1.65,-3.75);

\draw(2.05,-.15)--(2.15,-.15);
\draw(2.05,-.45)--(2.15,-.45);
\draw(2.05,-.75)--(2.15,-.75);
\draw(2.05,-1.05)--(2.15,-1.05);
\draw(2.05,-1.35)--(2.15,-1.35);
\draw(2.05,-1.65)--(2.15,-1.65);
\draw(2.05,-1.95)--(2.15,-1.95);
\draw(2.05,-2.25)--(2.15,-2.25);
\draw(2.05,-2.55)--(2.15,-2.55);
\draw(2.05,-2.85)--(2.15,-2.85);
\draw(2.05,-3.15)--(2.15,-3.15);
\draw(2.05,-3.45)--(2.15,-3.45);
\draw(2.05,-3.75)--(2.15,-3.75);

\draw(0.5,-.15)node{$\cdots$};
\draw(1,-.15)node{$\cdots$};
\draw(0.5,-1.7)node{$\cdots$};
\draw(1,-1.7)node{$\cdots$};
\draw(.5,-3.75)node{$\cdots$};
\draw(1,-3.75)node{$\cdots$};

\filldraw [black] (0,-.15) circle (2.5pt);
\filldraw [black] (0,-.75) circle (2.5pt);
\filldraw [black] (0,-1.05) circle (2.5pt);
\filldraw [black] (0,-1.65) circle (2.5pt);

\filldraw [black] (1.6,-.15) circle (2.5pt);
\filldraw [black] (1.6,-.45) circle (2.5pt);
\filldraw [black] (1.6,-.75) circle (2.5pt);
\filldraw [black] (1.6,-1.05) circle (2.5pt);
\filldraw [black] (1.6,-1.35)circle (2.5pt);

\filldraw [black] (2.1,-.15) circle (2.5pt);
\filldraw [black] (2.1,-.45) circle (2.5pt);
\filldraw [black] (2.1,-.75) circle (2.5pt);
\filldraw [black] (2.1,-1.05) circle (2.5pt);
\filldraw [black] (2.1,-1.35)circle (2.5pt);
\filldraw [black] (2.1,-1.65) circle (2.5pt);
\filldraw [black] (2.1,-2.25) circle (2.5pt);
\filldraw [black] (2.1,-2.55) circle (2.5pt);
\filldraw [black] (2.1,-3.15) circle (2.5pt);
\filldraw [black] (2.1,-3.45) circle (2.5pt);

\draw(2.2,-1.65)--(2.2,-3.45);
\draw(2.35,-2.55)node{$s$};
\draw(2,-2.55)--(2,-3.45);
\draw(1.85,-3)node{$r$};
\end{tikzpicture}}
\hfill 
{\begin{tikzpicture}
\draw(-.7,-1.625)node{$\hf_i^r\nu=$};

 \draw (0,0)--(0,-3.8);
\draw (1.6,0)--(1.6,-3.8);
\draw(2.1,0)--(2.1,-3.8);

\draw(-.05,-.15)--(.05,-.15);
\draw(-.05,-.45)--(.05,-.45);
\draw(-.05,-.75)--(.05,-.75);
\draw(-.05,-1.05)--(.05,-1.05);
\draw(-.05,-1.35)--(.05,-1.35);
\draw(-.05,-1.65)--(.05,-1.65);
\draw(-.05,-1.95)--(.05,-1.95);
\draw(-.05,-2.25)--(.05,-2.25);
\draw(-.05,-2.55)--(.05,-2.55);
\draw(-.05,-2.85)--(.05,-2.85);
\draw(-.05,-3.15)--(.05,-3.15);
\draw(-.05,-3.45)--(.05,-3.45);
\draw(-.05,-3.75)--(.05,-3.75);

\draw(1.55,-.15)--(1.65,-.15);
\draw(1.55,-.45)--(1.65,-.45);
\draw(1.55,-.75)--(1.65,-.75);
\draw(1.55,-1.05)--(1.65,-1.05);
\draw(1.55,-1.35)--(1.65,-1.35);
\draw(1.55,-1.65)--(1.65,-1.65);
\draw(1.55,-1.95)--(1.65,-1.95);
\draw(1.55,-2.25)--(1.65,-2.25);
\draw(1.55,-2.55)--(1.65,-2.55);
\draw(1.55,-2.85)--(1.65,-2.85);
\draw(1.55,-3.15)--(1.65,-3.15);
\draw(1.55,-3.45)--(1.65,-3.45);
\draw(1.55,-3.75)--(1.65,-3.75);

\draw(2.05,-.15)--(2.15,-.15);
\draw(2.05,-.45)--(2.15,-.45);
\draw(2.05,-.75)--(2.15,-.75);
\draw(2.05,-1.05)--(2.15,-1.05);
\draw(2.05,-1.35)--(2.15,-1.35);
\draw(2.05,-1.65)--(2.15,-1.65);
\draw(2.05,-1.95)--(2.15,-1.95);
\draw(2.05,-2.25)--(2.15,-2.25);
\draw(2.05,-2.55)--(2.15,-2.55);
\draw(2.05,-2.85)--(2.15,-2.85);
\draw(2.05,-3.15)--(2.15,-3.15);
\draw(2.05,-3.45)--(2.15,-3.45);
\draw(2.05,-3.75)--(2.15,-3.75);

\draw(0.5,-.15)node{$\cdots$};
\draw(1,-.15)node{$\cdots$};
\draw(0.5,-1.7)node{$\cdots$};
\draw(1,-1.7)node{$\cdots$};
\draw(.5,-3.75)node{$\cdots$};
\draw(1,-3.75)node{$\cdots$};

\filldraw [black] (0,-.15) circle (2.5pt);
\filldraw [black] (0,-.75) circle (2.5pt);
\filldraw [black] (0,-1.05) circle (2.5pt);
\filldraw [black] (0,-1.65) circle (2.5pt);
\filldraw [black] (0,-2.85) circle (2.5pt);
\filldraw [black] (0,-3.45) circle (2.5pt);
\filldraw [black] (0,-3.75) circle (2.5pt);

\filldraw [black] (1.6,-.15) circle (2.5pt);
\filldraw [black] (1.6,-.45) circle (2.5pt);
\filldraw [black] (1.6,-.75) circle (2.5pt);
\filldraw [black] (1.6,-1.05) circle (2.5pt);
\filldraw [black] (1.6,-1.35)circle (2.5pt);

\filldraw [black] (2.1,-.15) circle (2.5pt);
\filldraw [black] (2.1,-.45) circle (2.5pt);
\filldraw [black] (2.1,-.75) circle (2.5pt);
\filldraw [black] (2.1,-1.05) circle (2.5pt);
\filldraw [black] (2.1,-1.35)circle (2.5pt);
\filldraw [black] (2.1,-1.65) circle (2.5pt);
\filldraw [black] (2.1,-2.25) circle (2.5pt);

\draw(2.2,-1.65)--(2.2,-2.25);
\draw(2.7,-1.95)node{$s-r$};
\draw(-.1,-2.85)--(-.1,-3.75);
\draw(-.25,-3.3)node{$r$};
\end{tikzpicture}}
\hspace*{\fill}


Set $\nu^1:=\hf_i^r\nu$. By Lemma~\ref{lem:irr.Sp.dfilt}, we have $\head S^{\nu^1}\cong D^\la$ and $[S^{\nu^1}:D^\la]=1$. 
Now $\bphi_{i-1}(\nu^1)=0$, so $\phi_{i-1}(\la)=0$ by Lemma~\ref{lem:head.nodes}(i). Moreover $\eps_{i-1}(\nu^1)=\beps_{i-1}(\nu^1)=s-r$, so $\eps_{i-1}(\la)=s-r$, by Lemma~\ref{lem:head.nodes}(ii). 
Let $\xi:=\te_{i-1}^{(s-r)}\la$ and $\nu^2:=\he_{i-1}^{s-r}\nu^1$. Then $D^\xi\cong e_{i-1}^{(s-r)}D^\la$ and $D^\la\cong f_{i-1}^{(s-r)}D^\xi$, hence $\head S^{\nu^2}\cong D^\xi$, $[S^{\nu^2}:D^\xi]=1$, and
\begin{equation}\label{E200508_3}
\Ext^1_{\Si_n}(D^\la,D^\la)\cong\Ext^1_{\Si_{n-s+r}}(D^\xi,D^\xi).
\end{equation}

Note that $\eps_i(\nu^2)=\beps_i(\nu^2)=r$, so $\eps_i(\xi)=r$ by Lemma~\ref{lem:head.nodes}(ii). Let $\rho:=\te_i^{(r)}\xi$ and $\nu^3:=\he^{r}_i\nu^2$. Then  $D^\rho$ is in the head of $S^{\nu^3}$. Arguing as in { \sf Case 1}, but using Lemma~\ref{lem:rem.node.trick} instead of Lemma~\ref{lem:ad.node.trick},
we obtain an embedding 
\begin{equation}\label{E200508_5}
\Ext^1_{\Si_{n-s+r}}(D^\xi,D^\xi)\into \Ext^1_{\Si_{n-s}}(D^\rho,D^\rho).
\end{equation}

Note that $\beps_{i-1}(\nu^3)=\eps_{i-1}(\nu^3)=r$, so $\eps_{i-1}(\rho)=r$ by Lemma~\ref{lem:head.nodes}(ii). We set $\si:=\te_{i-1}^r\rho$ and $\nu^4:=\he^r_{i-1}\nu^3$.  Then $D^\si$ is in the head of $S^{\nu^4}$. Note that $\nu^4=\he_{i-1}^s\nu$, and $s=\beps_{i-1}(\nu)=\eps_{i-1}(\mu)$. So 
$S^{\nu^4}\cong e_{i-1}^{(s)}S^\nu\cong e_{i-1}^{(s)}D^\mu\cong D^{\te_{i-1}^s\mu}\cong D^\si$ is an irreducible Specht module. 
The non-regular runner of $\Ga(\nu^4)$ is $0$ and the non-restricted  runner of $\Ga(\nu^4)$ is $p-2$. By Lemma~\ref{lem:gener1} we deduce $\Ext^1_{\Si_{n-s}}(D^\rho,D^\rho)=0$, hence $\Ext^1_{\Si_n}(D^\la,D^\la)=0$ in view of (\ref{E200508_3}) and (\ref{E200508_5}). 
\end{proof}

\subsection{Concluding remarks on Theorem~\ref{TB}}
\label{subsec:Sp.result}

We note that the assumption that $e_i^{(\eps_i(\la))}D^\la$ is an (irreducible) Specht module in Theorem~\ref{TB} is equivalent to the assumption that $f_i^{(\phi_i(\la))}D^\la$ is an (irreducible) Specht module. Indeed, assume that $e_i^{(\eps_i(\la))}D^\la\cong D^\mu\cong S^\nu$. By Lemma~\ref{lem:rem.ad1}, we have $\phi_i(\mu)=\bphi_i(\nu)$, so $f_i^{(\phi_i(\la))}D^\la\cong f_i^{(\phi_i(\mu))}D^\mu\cong f_i^{(\bphi_i(\nu))}S^\nu\cong S^{\hat f_i^{\bphi_i(\nu)}\nu}$. 
Similarly in the other direction. 


\begin{Example}
\label{exm:res.irr.Sp}
Let 
$$\la=(p^2+1,p+2,(p+1)^{p-2},p,1^{p-1})\quad \text{and}\quad \nu=(p^2,p+1,2^{p-1}, 1^{p(p-1)-1}).$$ 
We have $\eps_0(\la)=2$ and $\te_0^2\la=(p^2,(p+1)^{p-1},p,1^{p-1})$. Then $\nu^\reg=\te_0^2\la$ and $S^\nu$ is irreducible by  Lemma~\ref{lem:irr.Sp.aba}, so $D^{\te_0^2\la}\cong S^\nu$. Hence $\Ext^1_{\Si_{2p(p+1)}}(D^\la,D^\la)=0$ by Theorem~\ref{TB}. 
\end{Example}

\section{Proof of Theorem~\ref{TSmallHeight}}
\label{SSmallHeight}

Throughout this section we assume that $p>2$.

\subsection{Strategy}
\label{subsec:strat.}
Let $h(\la)\leq p+2$. We prove Theorem~\ref{TSmallHeight} by induction on $n$ starting with $n<p$. To perform the inductive step, pick an $i\in I$ with $r:=\eps_i(\la)>0$. We have $h(\te^r_i\la)\leq p+2$. 
If $\phi_i(\la)=0$ then $\Ext^1_{\Si_n}(D^\la,D^\la)\cong\Ext^1_{\Si_{n-r}}(D^{\te^r_i\la}, D^{\te^r_i\la})$ by Lemma~\ref{L200509}. So we may assume that 
$\phi_i(\la)>0$. Let $A$ be the $i$-good node for $\la$ and $B$ be the $i$-cogood node for $\la$. Then by  
Corollary~\ref{C200509}, we have an embedding
$\Ext^1_{\Si_n}(D^\la,D^\la)\into\Ext^1_{\Si_{n-r}}(D^{\te^r_i\la}, D^{\te^r_i\la})$, unless $\la_A^B$ is not $p$-regular. So we may assume that $\la$ is of the form (\ref{E200511}), with 
$A:=(m+1,a+1)$ and $B:=(m+p,a)$. 
By Lemma~\ref{lem:fixed.top}, we may assume further that 
with $m\geq 1$ and $\la_1>a+1$.

We may further assume that $\eps_j(\la)=0$ for all $j\neq i$; in particular the top removable node of $\la$ has residue $i$. 
Indeed, if there is $j\neq i$ with $\eps_j(\la)\neq0$, then reasoning as above for $j$ instead of $i$ we deduce that $h(\la)\geq 2p-2$. If $p>3$, we have $2p-2>p+2$. Let $p=3$. Since $h(\la)\leq p+2=5$, we may assume that $\la$ is of the form $(b,a+1,a,a-1,a-2)$ with $b>a+1\geq 3$ and $(1,b)$ not good, or $(b,c,3,2,1)$ with $b>3$ and $(1,b)$ and $(2,c)$ not good. It is easy to see that in either of these exceptional cases there exists a unique $i$ with $\eps_i(\la)>0$.

Furthermore, we may assume that $h(\la)>p$. Indeed for $h(\la)<p$ we can do this by Proposition~\ref{PKS}, while for $h(\la)=p$, by the previous two paragraphs,  we may assume that $\la$ is of the form $(b,2,1^{p-2})$ with $b\equiv 1 \pmod{p}$ and apply Example~\ref{lem:hp,1}.  

Finally, in view of (\ref{ESignExt}), we may now also assume that $h(\laM)>p$. 

\subsection{Height $p+1$}
\label{subsec:heightp1}
In view of \S\ref{subsec:strat.}, we may assume that one of the following holds:
\begin{enumerate}
\item[{\rm (1)}] $\la=(b,c,2,1^{p-2})$, $i=p-1$, and one of the following  conditions holds: 
\begin{enumerate}
\item[{\rm (1a)}] $b>c>2$ and $\res(1,b)=\res (2,c)=p-1$; 
\item[{\rm (1b)}] $b>c\geq2$, $\res(1,b)=p-1$ and $\res (2,c)=0$;
\end{enumerate}      
\item[{\rm (2)}] $\la=(a+1+t,a+1,a^{p-2},a-1)$ with $t\equiv p-1\pmod{p}$ and $a\geq 2$.  
\end{enumerate}
(The case $b=c>2$ and $\res (2,c)=p-1$ in (1) is not present since then $h(\laM)=p-1$).

\begin{Proposition}
\label{lem:hp11}
If $\lambda$ is as in the case (1), then there exists $u>0$ and $\mu\in \Parreg_{n-u}$, with $h(\mu)\leq p+1$ and 
$\Ext^1_{\Si_n}(D^\lambda,D^\lambda)\into \Ext^1_{\Si_{n-u}}(D^\mu,D^\mu).$
\end{Proposition}

\begin{proof}
{\sf Case (1a)}. In this case we have $\phi_{i-1}(\la)=1$ and $\eps_{i-1}(\la)=0$. So, setting $\nu:=\tf_{i-1}\la={(b,c,2,2,1^{p-3})}$, by Lemma~\ref{L200509} we have  
$\Ext^1_{\Si_n}(D^\la,D^\la)\cong\Ext^1_{\Si_{n+1}}(D^{\nu},D^{\nu})$. Now, $\eps_i(\nu)=2$. So, by Corollary~\ref{C200509} applied to $\la=\nu$,  we obtain an embedding 
$\Ext^1_{\Si_{n+1}}(D^{\nu},D^{\nu})\into \Ext^1_{\Si_{n-1}}(D^{\te_i^2\nu},D^{\te_i^2\nu})$. 

\vspace{1mm}
{\sf Case (1b)}. If $p>3$, then $\phi_{i-1}(\la)=1$ and $\eps_{i-1}(\la)=0$. So, setting $\nu:=\tf_{i-1}\la=(b,c,2,2,1^{p-3})$, by Lemma~\ref{L200509} we have 
$\Ext^1_{\Si_n}(D^\la,D^\la)\cong\Ext^1_{\Si_{n+1}}(D^{\nu},D^{\nu})$.  
Note that $\eps_{i-1}(\te_i\nu)=2$, while $\phi_{i-1}(\te_i\nu)=0$, so using Corollary~\ref{C200509} and Lemma~\ref{L200509}, we get 
$$\Ext^1_{\Si_{n+1}}(D^\nu,D^\nu)\into \Ext^1_{\Si_n}(D^{\te_i\nu},D^{\te_i\nu})\cong \Ext^1_{\Si_{n-2}}(D^{\te_{i-1}^2\te_i\nu},D^{\te_{i-1}^2\te_i\nu}).$$ 
Let $p=3$. We may assume that $c<b-1$ since otherwise $h((\laM))=2=p-1$. 
Using Corollary~\ref{C200509} and Lemma~\ref{L200509}, we get 
\begin{align*}
\Ext^1_{\Si_n}(D^\la,D^\la)\cong \Ext^1_{\Si_{n+2}}(D^{\tf_1^2\la},D^{\tf_1^2\la})\into &\Ext^1_{\Si_{n+1}}(D^{\te_2\tf_1^2\la},D^{\te_2\tf_1^2\la})
\\
\cong &\Ext^1_{\Si_{n-2}}(D^{\te_1^3\te_2\tf_1^2\la},D^{\te_1^3\te_2\tf_1^2\la}),
\end{align*}
where $\te_1^3\te_2\tf_1^2\la=(b-2,c,2,1)$, completing the proof. 
\end{proof}

We now consider the case (2). We require the following two lemmas.

\begin{Lemma}
\label{lem:surj1}
Let $a\geq 2$,  $t\equiv p-1 \pmod{p}$ and $c\equiv p-2  \pmod{p}$. Then there is a surjective homomorphism of Specht modules  
$S^{(a+t,a,2,1^c)}\onto S^{(a+t,a,1^{c+2})}.$
\end{Lemma}

\begin{proof}
Using \cite[Theorem 8.15]{JamesBook}, it is easy to see that
$$\Hom_{\Si_{c+2}}(S^{(2,1^c)}, S^{(1^{c+2})})\cong \Hom_{\Si_{c+2}}(S^{(2,1^c)}, \sgn)\cong \k.$$
So by \cite[Theorem 2.2]{FL}, we have 
\begin{equation}\label{E010621}
\Hom_{\Si_{2a+t+c+2}}(S^{(a+t,a,2,1^c)}, S^{(a+t,a,1^{c+2})})\cong \Hom_{\Si_{c+2}}(S^{(2,1^c)}, S^{(1^{c+2})})\cong \k.
\end{equation}
So there is a unique up to scalar non-zero homomorphism $S^{(a+t,a,2,1^c)}\to S^{(a+t,a,1^{c+2})}$. We show by induction on $a$ that this homomorphism is surjective. For $a=2$, the Specht module $S^{(t+2,2,1^{c+2})}$ is irreducible for example by Lemma~\ref{lem:irr.Sp.aba}, and the result is clear. We now assume that the statement holds for $a$ and show it for $a+1$.  Let $j:=\res(1,a+t+1)=\res (2,a+1)$. The surjection $S^{(a+t,a,2,1^c)}\onto S^{(a+t,a,1^{c+2})}$ yields  a surjection
$f_j^{(2)}S^{(a+t,a,2,1^c)}\onto f_j^{(2)}S^{(a+t,a,1^{c+2})}.$
Moreover, by Lemma~\ref{lem:res.Sp}, $f_j^{(2)}S^{(a+t,a,2,1^c)}$ (resp. $f_j^{(2)}S^{(a+t,a,1^{c+2})})$ has a Specht filtration with top Specht factor $S^{(a+t+1,a+1,2,1^c)}$ (resp. $S^{(a+t+1,a+1,1^{c+2})}$). So there is a surjection 
$$f_j^{(2)}S^{(a+t,a,2,1^c)}\overset{\psi}{\onto} S^{(a+t+1,a+1,1^{c+2})}.$$
If $S^\mu$ is one of the Specht factors in the Specht filtration of $f_j^{(2)}S^{(a+t,a,2,1^c)}$ coming from Lemma~\ref{lem:res.Sp} which is not the top Specht factor, then $\mu \mathbb{4} (a+t+1,a+1,1^{c+2})$ and so $\Hom_{\Si_{2a+t+c+4}}(S^\mu, S^{(a+t+1,a+1,1^{c+2})})=0$ by \cite[Corollary 13.17]{JamesBook}. Hence $\psi$ induces a surjection $S^{(a+t+1,a+1,2,1^c)}\onto S^{(a+t+1,a+1,1^{c+2})}$.
\end{proof} 

\begin{Lemma}
\label{lem:head1}
Let $\nu:=(a+t,a,1^{a(p-1)-1})$ and $\mu:=(a+t,a^{p-1},a-1)$, with  $a\geq 1$ and $t\equiv p-1\pmod{p}$. 
Then 
$[S^\nu:D^\mu]=1$ and $\Hom_{\Si_{|\nu|}}(S^\nu,D^\mu)\cong \k$. 
\end{Lemma}

\begin{proof}
For the decomposition number, it suffices to note that $\mu=\nu^\reg$. 
We now prove that $D^{\mu}$ appears in the head of $S^{\nu}$ by induction on $a$, the case $a=1$ being clear. We now assume that the statement holds for $a$ and show it for $a+1$. Let $j:=\res(1,a+t+1)=\res (2,a+1)$. One should consider three different cases for $j=p-1$, $j=p-2$ and $j\in \{0,\dots, p-3\}$. All three cases can be treated in a similar way, so we only provide full details for the most demanding case $j=p-1$.

Applying $f_j^{{(3})}$ to a surjective map $S^\nu\onto D^\mu$ which exists by the inductive assumption, we get a surjection 
$f_j^{(3)}S^\nu\to f_j^{(3)} D^\mu.$ Now $f_j^{(3)}D^\mu$ has simple head $D^{\tf_j^3\mu}$, where $\tf_j^3\mu=(a+t+1,a+1,a^{p-1})$. Therefore, $D^{\tf_j^3\mu}$ appears in the head of $f_j^{(3)}S^\nu$. By Lemma~\ref{lem:res.Sp}, 
$f_j^{(3)}S^\nu\sim S^{\nu^1}\mid S^{\nu^2}\mid S^{\nu^3}\mid S^{\nu^4}$, where $\nu^4=\hf_j^3\nu=(a+t+1,a+1,2,1^{a(p-1)-2})$, $\nu^3=(a+t+1,a+1,1^{a(p-1)})$, $\nu^2=(a+t+1,a,2,1^{a(p-1)-1})$ and $\nu^1=(a+t,a+1,2,1^{a(p-1)-1})$.

Note that $\eps_{p-2}(\tf_j^3 \mu)=0$. Moreover the node $(1,a+t)$ of $\nu^1$ (resp. the node $(2,a)$ of $\nu^2$) is the highest $(p-2)$-removable node and it is normal. Hence $\Hom_{\Si_{|\nu|+3}}(S^{\nu^1}, D^{\tf_j^3\mu})\cong\Hom_{\Si_{|\nu|+3}}(S^{\nu^2}, D^{\tf_j^3\mu})=0$ by Lemma~\ref{lem:rem.node.trick}.
Thus, if $N\sim S^{\nu^1}\mid S^{\nu^2}\subseteq f_j^{(3)}S^\nu$ is the submodule of $f_j^{(3)}S^\nu$ corresponding to the bottom two Specht factors, then $\Hom_{\Si_{|\nu|+3}}(N, D^{\tf_j^3\mu})=0$, hence $D^{\tf_j^3\mu}$ appears in the head of $(f_j^{(3)}S^\nu)/N\sim S^{\nu^3}\mid S^{\nu^4}$. 

By Lemma~\ref{lem:surj1} we have a surjection, $S^{\la^4}\onto S^{\la^3}$. Using Lemma~\ref{LRowRem} and \cite[Theorem 24.1]{JamesBook}, we get $[S^{\nu^4}:D^{\tf_j^3\mu}]=[S^{\nu^3}:D^{\tf_j^3\mu}]=1$.   So $D^{\tf_j^3\mu}$ is in the head of $S^{\nu^3}$.   Now 
$$S^{(a+t+1,a+1,1^{(a+1)(p-1)-1})}\cong f_1\ldots f_{p-2} S^{\nu^3}\onto f_1\ldots f_{p-2} D^{\tf_j^3\mu}\cong D^{(a+t+1,(a+1)^{p-1},a)}$$ 
by the exactness of $f$'s.
\end{proof}

\begin{Proposition}
\label{lem:hp12}
Let $a,t,\la$ be as in the case (2), and $\mu:=(a+t,a^{p-1},a-1)$.  
We have an embedding
$\Ext^1_{\Si_n}(D^\la,D^\la)\into \Ext^1_{\Si_{n-2}}(D^\mu,D^\mu).$
\end{Proposition}

\begin{proof}
We have $i=\res(1,a+t+1)=\res(2,a+1)$ and $e_i^{(2)}D^\la\cong D^\mu$. From the short exact sequence $0\to \rad f_i^{(2)}D^\mu \to f_i^{(2)}D^\mu\to D^\la\to 0$, using Lemma~\ref{lem:shap}, we get an exact sequence 
$$
\Hom_{\Si_n}(\rad (f_i^{(2)}D^\mu), D^\la)\to\Ext^1_{\Si_n}(D^\la,D^\la)\to
\Ext^1_{\Si_n}(f_i^{(2)}D^\mu,D^\la)\cong \Ext^1_{\Si_{n-2}}(D^\mu,D^\mu).
$$
We complete the proof by showing that  $\Hom_{\Si_n}(\rad (f_i^{(2)}D^\mu), D^\la)=0$. 
Let $\nu=(a+t,a,1^{a(p-1)-1})$. By Lemma~\ref{lem:head1}, 
there is a submodule $X\subseteq S^\nu$ such that 
$S^\nu/X\cong D^\mu$. Then $f^{(2)}_iS^\nu/f^{(2)}_iX\cong f^{(2)}_iD^\mu$, hence there exists a submodule $f^{(2)}_iX\subseteq Y\subseteq f^{(2)}_iS^\nu$ with $Y/f^{(2)}_iX\cong \rad (f_i^{(2)}D^\mu)$ and $f^{(2)}_iS^\nu/Y\cong D^\la$. It suffices to show that $\Hom_{\Si_n}(Y, D^\la)=0$. 

By Lemma~\ref{lem:res.Sp}, $f^{(2)}_iS^\nu\sim S^{\nu^1}\mid\ldots\mid S^{\nu^s}$ with $\nu^s=(a+t+1,a+1,1^{a(p-1)-1})$. Also, using Lemma~\ref{LRowRem}, we have 
$
[S^{\nu^s}:D^\la]=[S^{(1^{a(p-1)-1})}:D^{(a^{p-2},a-1)}]=1. 
$ So it suffices to prove that $\Hom_{\Si_n}(S^{\nu^u}, D^\la)=0$ for $u<s$. Note that $\nu^u_1=a+t$ or $\nu^u_2=a$. 
If $\nu^u_1=a+t$ the $(i-1)$-removable node $(1,a+t)$ of $\nu^u$ is normal. Since $\eps_{i-1}(\la)=0$ we obtain by Lemma~\ref{lem:rem.node.trick} that $\Hom_{\Si_n}(S^{\nu^u}, D^\la)=0$. If $\nu^u_1=a+t+1$ then $\nu^u_2=a$. In this case the $(i-1)$-removable node $(2,a)$ of $\nu^u$ is normal and so we get again that $\Hom_{\Si_n}(S^{\nu^u}, D^\la)=0$.
\end{proof}

\subsection{Height $p+2$}
\label{subsec:heightp2}
In view of \S\ref{subsec:strat.}, we may assume that one of the following holds:
\begin{enumerate}
\item[{\rm (1)}] $\la=(b,c,d,2,1^{p-2})$, $i=p-2$, and one of the following  conditions holds: 
\begin{enumerate}
\item[{\rm (1a)}] $b>c=d=2$ and $\res (1,b)=p-2$;
\item[{\rm (1b)}] $b>c>d\geq 2$ and $\res(1,b)=\res (2,c)=p-2, \res(3,d)=p-1$; 
\item[{\rm (1c)}] $b=c=d>2$, and $\res(3,d)=p-2$;
\item[{\rm (1d)}]  $b=c>d>2$ and $\res(2,c)=\res(3,d)=p-2$.
\item[{\rm (1e)}]  $b=c>d>2$ and $\res(2,c)=p-2, \res(3,d)=0$.
\item[{\rm (1f)}]  $b>c\geq d>2$ and $\res(1,b)=\res(3,d)=p-2, \res(2,c)=p-1$.
\item[{\rm (1g)}]  $b>c=d>2$ and $\res(1,b)=p-2,\res(3,d)=p-1$.
\item[{\rm (1h)}]  $b>c>d>2$ and $\res(1,b)=\res(2,c)=\res(3,d)=p-2$.
\item[{\rm (1i)}]  $b>c>d>2$ and $\res(1,b)=p-2,\res(2,c)=p-1$, $\res(3,d)=0$.
\end{enumerate}
\item[{\rm (2)}] $\la=(a+s+t+1,a+s+1,a+1,a^{p-2},a-1)$ with $a\geq 2$, $t\equiv p-1\pmod{p}$ and $s\equiv p-1\pmod{p}$.
\item[{\rm (3)}] $\la=(a+t+1,a+1,a^{p-2},a-1,b)$ with $b\geq 1$, $t\equiv p-1\pmod{p}$ and $a-b\equiv p-2\pmod{p}$.
\item[{\rm (4)}]  $\la=(a+s+t+1,a+s+1,a+1,a^{p-2},a-1)$ with $a\geq 2$, $t\equiv p-2\pmod{p}$ and $s\equiv 0\pmod{p}$.
\item[{\rm (5)}]  $\la=((a+t+1)^2,a+1,a^{p-2},a-1)$ with $a\geq 2$ and $t\equiv p-1\pmod{p}$.
\end{enumerate}

Moreover, by (\ref{ESignExt}) and what we have already proved in \S\ref{subsec:heightp1}, we may ignore the cases above for which $h(\laM)\leq p+1$.

\begin{Proposition}
\label{lem:hp21}
If $\lambda$ is as in the case (1), then there exists $u>0$ and $\mu\in \Parreg_{n-u}$, with $h(\mu)\leq p+2$ and 
$\Ext^1_{\Si_n}(D^\lambda,D^\lambda)\into \Ext^1_{\Si_{n-u}}(D^\mu,D^\mu).$
\end{Proposition}

\begin{proof}
{\sf Case (1a)}. In this case we always have $h(\laM)\leq p+1$. 

\vspace{1mm}
{\sf Case (1b)}. If $p>3$, then $\phi_{i-1}(\la)=1$ and $\eps_{i-1}(\la)=0$. So setting $\nu:=\tf_{i-1}\la=(b,c,d,2,2,1^{p-3})$, by Lemma~\ref{L200509} we have 
$\Ext^1_{\Si_n}(D^\la,D^\la)\cong\Ext^1_{\Si_{n+1}}(D^{\nu},D^{\nu})$.  
Now, $\eps_i(\nu)=2$ and by Corollary~\ref{C200509} we get $\Ext^1_{\Si_{n+1}}(D^{\nu},D^{\nu})\into \Ext^1_{\Si_{n-1}}(D^{\te_i^2\nu},D^{\te_i^2\nu})$.  

Let $p=3$. We may assume that $c>d+1$ since otherwise $h(\laM)=4=p+1$. 
Setting $\nu:=\tf^2_{0}\la=(b,c,d+1,2,2)$, we get
\begin{align*}
\Ext^1_{\Si_n}(D^\la,D^\la)
\cong
\Ext^1_{\Si_{n+2}}(D^{\nu},D^{\nu})
\into 
\Ext^1_{\Si_{n}}(D^{\te_1^2\nu},D^{\te_1^2\nu})
\cong\Ext^1_{\Si_{n-4}}(D^{\te_0^4\te_1^2\nu},D^{\te_0^4\te_1^2\nu}),
\end{align*}
where the isomorphisms hold by Lemma~\ref{L200509} and  the embedding holds by Corollary~\ref{C200509}. 

\vspace{1mm}
{\sf Case (1c)}. In this case we have $h(\laM)=p-2$.

\vspace{1mm}
{\sf Case (1d)}. If $p>3$, then $\phi_{i-1}(\la)=1$ and $\eps_{i-1}(\la)=0$ and the result follows as in case (1b). If $p=3$ then $h(\laM)=4$. 

\vspace{1mm}
{\sf Case (1e)}.  If $p=3$ then $h(\laM)\leq 4$. 
For $p>3$, setting $\nu:=\tf_{i-1}\la$, we have
\begin{align*}
\Ext^1_{\Si_n}(D^\la,D^\la)
\cong
\Ext^1_{\Si_{n+1}}(D^{\nu},D^{\nu})
\into  
\Ext^1_{\Si_{n}}(D^{\te_i\nu},D^{\te_i\nu})
\cong
\Ext^1_{\Si_{n-2}}(D^{\te_{i-1}^2\te_i\nu},D^{\te_{i-1}^2\te_i\nu}),
\end{align*}
where the isomorphisms hold by Lemma~\ref{L200509} and  the embedding holds by Corollary~\ref{C200509}. 

\vspace{1mm}
{\sf Cases (1f)-(1i)} can be verified in a similar way so we leave the details to the reader.  
\end{proof}

For the case $(2)$ we need the following two lemmas:

\begin{Lemma}
\label{lem:surj2}
Let $t\equiv p-1 \pmod{p}$, $s\equiv p-1 \pmod{p}$ and $c\equiv p-2  \pmod{p}$. For $a\geq 2$,  set $\al(a):=(a+s+t,a+s,a,2,1^c)$ and $\be(a):=(a+s+t,a+s,a,1^{c+2})$. 
Then there is a surjective homomorphism of Specht modules   $S^{\al(a)}\onto S^{\be(a)}$.
\end{Lemma}

\begin{proof}
To ease notation set $l:=a+s+t$. By \cite[Theorem 2.2]{FL} and (\ref{E010621}), we have 
$$\Hom_{\Si_{|\al(a)|}}(S^{\al(a)}, S^{\be(a)})\cong \Hom_{\Si_{{|\al(a)|}-l}}(S^{(a+s,a,2,1^c)}, S^{(a+s,a,1^{c+2})})\cong \k.$$
So there is a unique up to scalar non-zero homomorphism 
$\phi_a: S^{\al(a)}\to S^{\be(a)}.$ 
We show by induction on $a\geq 2$ that $\phi_a$ is surjective. For the induction base, by Lemma~\ref{lem:surj1}, we have a surjection  $S^{(2+s+t,2+s,2,1^c)}\onto S^{(2+s+t,2+s,1^{c+2})}$. Applying $f_{p-1}f_{p-2}$ to this we obtain 
$$S^{\al(2)}\cong f_{p-1}f_{p-2}S^{(2+s+t,2+s,2,1^c)}\onto f_{p-1}f_{p-2}S^{(2+s+t,2+s,1^{c+2})}\cong S^{\be(2)}.$$ 
For the inductive step, let $j:=\res(1,a+s+t+1)=\res(2,a+s+1)=\res(3,a+1)$.
Then the surjection $\phi_a:S^{\al(a)}\onto S^{\be(a)}$ yields  a surjection $f_j^{(3)}S^{\al(a)}\onto f_j^{(3)}S^{\be(a)}.$ 
 The rest of the argument is now similar to the proof of Lemma~\ref{lem:surj1}.
 \end{proof}
 
 \begin{Lemma}
\label{lem:head2}
Let $\nu:=(a+s+t,a+s,a,1^{a(p-1)-1})$ and $\mu:=(a+s+t,a+s, a^{p-1},a-1)$ with $a\geq 1$, $t\equiv p-1\pmod{p}$ and $s\equiv p-1\pmod{p}$. Then
$[S^\nu:D^\mu]=1$ and $\Hom_{\Si_{|\nu|}}(S^\nu,D^\mu)\cong \k$. 
\end{Lemma}

\begin{proof}
For the decomposition number, it suffices to note that $\mu=\nu^\reg$. 
We now prove that $D^{\mu}$ appears in the head of $S^{\nu}$ by induction on $a$, the case $a=1$ being clear. We now assume that the statement holds for $a$ and show it for $a+1$. Let $j:=\res(1,a+s+t+1)=\res (2,a+s+1)=\res(3,a+1)$. One should consider three different cases for $j=p-2$, $j=p-3$ and $j\neq p-2,p-3$. All three cases can be treated in a similar way, so we only provide full details for the most demanding case $j=p-2$.

Applying $f_j^{{(4})}$ to a surjective map $S^\nu\onto D^\mu$ which exists by the inductive assumption, we get a surjection 
$f_j^{(4)}S^\nu\to f_j^{(4)} D^\mu.$ Now $f_j^{(4)}D^\mu$ has simple head $D^{\tf_i^4\mu}$, where $\tf_i^4\mu=(a+s+t+1,a+s+1,a+1,a^{p-1})$. Therefore, $D^{\tf_i^4\mu}$ appears in the head of $f_j^{(4)}S^\mu$. By Lemma~\ref{lem:res.Sp}, 
$f_j^{(4)}S^\nu\sim S^{\nu^1}\mid S^{\nu^2}\mid S^{\nu^3}\mid S^{\nu^4}\mid S^{\nu^5}$, where 
$
\nu^5=(a+s+t+1,a+s+1,a+1,2,1^{a(p-1)-2})$, 
$\nu^4=(a+s+t+1,a+s+1,a+1,1^{a(p-1)})$,  
$\nu^3=(a+s+t+1,a+s+1,a,2,1^{a(p-1)-1})$,  
$\nu^2=(a+s+t+1,a+s,a+1,2,1^{a(p-1)-1})$,  
$\nu^1=(a+s+t,a+s+1,a+1,2,1^{a(p-1)-1})$.

Note that $\eps_{j-1}(\tf_i^4 \mu)=0$. Moreover the node $(1,a+s+t)$ of $\nu^1$ (resp. the node $(2,a+s)$ of $\nu^2$ and the node $(3,a)$ of $\nu^3$) is the highest $(j-1)$-removable node and it is normal. Hence $\Hom_{\Si_{|\nu|+4}}(S^{\nu^1}, D^{\tf_j^4\mu})\cong \Hom_{\Si_{|\nu|+4}}(S^{\nu^2}, D^{\tf_j^4\mu})\cong \Hom_{\Si_{|\nu|+4}}(S^{\nu^3}, D^{\tf_j^4\mu})=0$ by Lemma~\ref{lem:rem.node.trick}.
Therefore, if $N\sim S^{\nu^1}\mid S^{\nu^2}\mid S^{\nu^3}\subseteq f_j^{(4)}S^\nu$ is the submodule of $f_j^{(4)}S^\nu$ corresponding to the bottom three Specht factors, then $\Hom_{\Si_{|\nu|+4}}(N, D^{\tf_j^4\mu})=0$. So $D^{\tf_j^4\mu}$ appears in the head of $(f_j^{(4)}S^\nu)/N\sim S^{\nu^4}\mid S^{\nu^5}$. 

By Lemma~\ref{lem:surj2} we have a surjection, $S^{\nu^5}\onto S^{\nu^4}$. Using Lemma~\ref{LRowRem} and \cite[Theorem 24.1]{JamesBook}, we get $[S^{\nu^5}:D^{\tf_j^4\mu}]=[S^{\nu^4}:D^{\tf_j^4\mu}]=1$.   So $D^{\tf_j^4\mu}$ is in the head of $S^{\nu^4}$.   Now 
\begin{align*}
f_0\ldots f_{p-3} S^{\nu^4}&\cong S^{(a+s+t+1,a+s+1,a+1,1^{(a+1)(p-1)-1})},\cr
f_0\ldots f_{p-3} D^{\tf_j^4\mu}&\cong D^{(a+s+t+1,a+s+1,(a+1)^{p-1},a)}.
\end{align*}
So, we have a surjection $S^{(a+s+t+1,a+s+1,a+1,1^{(a+1)(p-1)-1})}\onto D^{(a+s+t+1,a+s+1,(a+1)^{p-1},a)}$ by the exactness of $f$'s.
\end{proof}

\begin{Proposition}
\label{lem:hp22}
Let $a,s,t,\la$ be as in the case (2). Set 
$\mu:=(a+s+t,a+s, a^{p-1},a-1)$. Then we have an embedding
$\Ext^1_{\Si_n}(D^\la,D^\la)\into \Ext^1_{\Si_{n-3}}(D^\mu,D^\mu).$
\end{Proposition}

\begin{proof}
We have $i=\res(1,a+s+t+1)=\res(2,a+s+1)=\res(3,a+1)$. Then $e_i^{(3)}D^\la\cong D^\mu$ and from the short exact sequence $0\to \rad f_i^{(3)}D^\mu \to f_i^{(3)}D^\mu\to D^\la\to 0$, using Lemma~\ref{lem:shap}, we get an exact sequence 
$$
\Hom_{\Si_n}(\rad (f_i^{(3)}D^\mu), D^\la)\to\Ext^1_{\Si_n}(D^\la,D^\la)\to
\Ext^1_{\Si_n}(f_i^{(3)}D^\mu,D^\la)\cong \Ext^1_{\Si_{n-3}}(D^\mu,D^\mu).
$$
We complete the proof by showing that  $\Hom_{\Si_n}(\rad (f_i^{(3)}D^\mu), D^\la)=0$. 

Let $\nu:=(a+s+t,a+s,a,1^{a(p-1)-1})$. By Lemma~\ref{lem:head2}, 
there is a submodule $X\subseteq S^\nu$ such that 
$S^\nu/X\cong D^\mu$. Then $f^{(3)}_iS^\nu/f^{(3)}_iX\cong f^{(3)}_iD^\mu$, hence there exists a submodule $f^{(3)}_iX\subseteq Y\subseteq f^{(3)}_iS^\nu$ with $Y/f^{(3)}_iX\cong \rad (f_i^{(3)}D^\mu)$ and $f^{(3)}_iS^\nu/Y\cong D^\la$. It suffices to show that $\Hom_{\Si_n}(Y, D^\la)=0$. 

By Lemma~\ref{lem:res.Sp}, $f^{(3)}_iS^\nu\sim S^{\nu^1}\mid\ldots\mid S^{\nu^k}$ with $\nu^k=(a+s+t+1,a+s+1,a+1,1^{a(p-1)-1})$. Also, using Lemma~\ref{LRowRem}, we have 
$
[S^{\nu^k}:D^\la]=[S^{(1^{a(p-1)-1})}:D^{(a^{p-2},a-1)}]=1. 
$ So it suffices to prove that $\Hom_{\Si_n}(S^{\nu^u}, D^\la)=0$ for $u<k$. Note that $\nu^u_1=a+s+t$ or $\nu^u_2=a+s$ or $\nu^u_3=a$. 
If $\nu^u_1=a+s+t$ the $(i-1)$-removable node $(1,a+s+t)$ of $\nu^u$ is normal. Since $\eps_{i-1}(\la)=0$ we obtain by Lemma~\ref{lem:rem.node.trick} that $\Hom_{\Si_n}(S^{\nu^u}, D^\la)=0$. If $\nu^u_1=a+s+t+1$ and $\nu^u_2=a+s$, then the $(i-1)$-removable node $(2,a+s)$ of $\nu^u$ is normal and so we get again that $\Hom_{\Si_n}(S^{\nu^u}, D^\la)=0$. Finally if $\nu^u_1=a+s+t+1$ and $\nu^u_2=a+s+1$ then $\nu^u_3=a$ and the $(i-1)$-removable node $(3,a)$ of $\nu^u$ is normal and so $\Hom_{\Si_n}(S^{\nu^u}, D^\la)=0$.
\end{proof}

For the case $(3)$ we require the following three lemmas:

\begin{Lemma}
\label{lem:surj3}
Let  $t\equiv p-1 \pmod{p}$, $s\equiv p-4 \pmod{p}$ and $c\equiv p-2  \pmod{p}$. For $a\geq2$, we set $\al(a):=(a+s+t,a+s,a,2,1^c)$ and $\be(a):=(a+s+t,a+s,a,1^{c+2})$. 
Then there is a surjective homomorphism of Specht modules  
$S^{\al(a)}\onto S^{\be(a)}$.
\end{Lemma}

\begin{proof}
If $p=3$, the lemma is equivalent to Lemma~\ref{lem:surj2}, so we assume that $p>3$. 
By \cite[Theorem 2.2]{FL} and (\ref{E010621}) we have 
$$\Hom_{\Si_{|\al(a)|}}(S^{\al(a)}, S^{\be(a)})\cong \Hom_{\Si_{c+2}}(S^{(2,1^c)}, S^{(1^{c+2})})\cong \k.$$
So there is a unique up to scalar non-zero homomorphism 
$\phi_a:S^{\al(a)}\to S^{\be(a)}.$ 
We show by induction on $a\geq 2$ that $\phi_a$ is surjective. For the induction base, by Lemma~\ref{lem:surj1}, we have a surjective map $S^{(2+s+t,2+s,2,1^c)}\onto S^{(2+s+t,2+s,1^{c+2})}$. Applying $f_{p-1}f_{p-2}$ we obtain 
$$S^{\al(2)}\cong f_{p-1}f_{p-2}S^{(2+s+t,2+s,2,1^c)}\onto f_{p-1}f_{p-2}S^{(2+s+t,2+s,1^{c+2})}\cong S^{\be(2)}.$$ 
For the inductive step, let $j:=\res(1,a+s+t+1)=\res(2,a+s+1)$. Then $\res(3,a+1)=j+3$.
The surjection $\phi_a:S^{\al(a)}\onto S^{\be(a)}$ yields  a surjection $f_j^{(2)}S^{\al(a)}\onto f_j^{(2)}S^{\be(a)}.$
 Arguing as in the the proof of Lemma~\ref{lem:surj1} we get a surjection $S^{(a+s+t+1,a+s+1,a,2,1^c)}\overset{\psi}{\onto}S^{(a+s+t+1,a+s+1,a,1^{c+2})}.$
 Applying $f_{j+3}$ to $\psi$ and arguing as in Lemma~\ref{lem:surj1} we get a surjection $S^{\al(a+1)}\onto S^{\be(a+1)}.$
\end{proof}

\begin{Lemma}
\label{lem:head3}
Let $\nu:=(a+t,a,b+1,1^{a(p-1)-2})$ and $\mu:=(a+t,a^{p-1},a-1,b)$ with $a>b\geq 0$, $t\equiv p-1\pmod{p}$ and  $a-b\equiv p-3\pmod{p}$. Then
$[S^\nu:D^\mu]=1$ and $\Hom_{\Si_{|\nu|}}(S^\nu,D^\mu)\cong \k$. 
\end{Lemma}

\begin{proof}
For the equality $[S^\nu:D^\mu]=1$ , it suffices to note that $\mu=\nu^\reg$. We now write $\nu=\nu(b)=(b+s+t,b+s,b+1,1^{(b+s)(p-1)-2})$ and $\mu=\mu(b)=(b+s+t,(b+s)^{p-1},b+s-1,b)$ with  $s>0$ and $s \equiv p-3\pmod{p}$, and proceed by induction on $b$. For $b=0$ the result follows by Lemma~\ref{lem:head1} with  $a=s$. Suppose there is a non-zero homomorphism $\phi_b:S^{\nu(b)}\to D^{\mu(b)}$ (which is automatically surjective). Let $j:=\res(1,b+s+t+1)=\res (2,b+s+1)$. In constructing $\phi_{b+1}$, 
one should consider three different cases: $j=p-2$, $j=p-3$, or  $j\neq p-2,p-3$. All cases are treated in a similar way, so we only provide full details for the most demanding case $j=p-2$.

Assume first that $p>3$. 
Applying $f_j^{{(3})}$ to $\phi_b$, we get a surjection 
$f_i^{(3)}S^{\nu(b)}\onto f_i^{(3)} D^{\mu(b)}.$ Now $f_i^{(3)}D^{\mu(b)}$ has simple head $D^{\tf_j^3\mu(b)}$, where $\tf_j^3\mu(b)=(b+s+t+1,b+s+1,(b+s)^{p-1},b)$. Therefore, $D^{\tf_j^3\mu(b)}$ appears in the head of $f_j^{(3)}S^{\nu(b)}$. By Lemma~\ref{lem:res.Sp}, 
$f_j^{(3)}S^{\nu(b)}\sim S^{\nu^1}\mid S^{\nu^2}\mid S^{\nu^3}\mid S^{\nu^4}$, with $\nu^4=(b+s+t+1,b+s+1,b+1,2,1^{(b+s)(p-1)-3})$, $\nu^3=(b+s+t+1,b+s+1,b+1,1^{(b+s)(p-1)-1})$, $\nu^2=(b+s+t+1,b+s,b+1,2,1^{a(p-1)-2})$, $\nu^1=(b+s+t,b+s+1,b+1,2,1^{a(p-1)-2})$.
Note that $\eps_{j-1}(\tf_i^3 \mu(b))=0$. Moreover the node $(1,b+s+t)$ of $\nu^1$ (resp. the node $(2,b+s)$ of $\nu^2$) is the highest $(j-1)$-removable node and it is normal. Hence $\Hom_{\Si_{|\nu|+3}}(S^{\nu^1}, D^{\tf_j^3\mu(b)})=\Hom_{\Si_{|\nu|+3}}(S^{\nu^2}, D^{\tf_j^3\mu(b)})=0$ by Lemma~\ref{lem:rem.node.trick}.
Therefore, if $N\sim S^{\nu^1}\mid S^{\nu^2}\subseteq f_j^{(3)}S^{\nu(b)}$ is the submodule of $f_j^{(3)}S^{\nu(b)}$ corresponding to the bottom two Specht factors, then $\Hom_{\Si_{|\nu|+3}}(N, D^{\tf_j^3\mu(b)})=0$. Hence $D^{\tf_j^3\mu(b)}$ appears in the head of $M:=(f_j^{(3)}S^{\la(b)})/N\sim S^{\nu^3}\mid S^{\nu^4}$. Thus we have a surjective map $M\overset{\psi}{\onto} D^{\tf_j^3\mu(b)}$. We apply now $f_{j+3}$ and $\psi$ yields a surjective map
$$f_{j+3}M\onto f_{j+3}D^{\tf_j^3\mu(b)}\cong D^{\tf_{j+3}\tf_j^3\mu(b)},$$
where ${\tf_{j+3}\tf_j^3\mu(b)}=(b+s+t+1,b+s+1,(b+s)^{p-1},b+1)$. Note that  
$$f_{j+3}S^{\nu^3}\cong S^{(b+s+t+1,b+s+1,b+2,1^{(b+s)(p-1)-1})},\ f_{j+3}S^{\nu^4}\cong S^{(b+s+t+1,b+s+1,b+2,2,1^{(b+s)(p-1)-3})}.$$ 
By Lemma~\ref{lem:surj3}, we have a surjection 
$f_{j+3}S^{\nu^4}\onto f_{j+3}S^{\nu^3}$. Using Lemma~\ref{LRowRem} one easily sees that $[f_{j+3}S^{\nu^4}:D^{\tf_{j+3}\tf_j^3\mu(b)}]=[f_{j+3}S^{\nu^3}:D^{\tf_{j+3}\tf_j^3\mu(b)}]=1$. Hence $D^{\tf_{j+3}\tf_j^3\mu}$ is in the head of $f_{j+3}S^{\nu^3}$. Now $e_{j+3}f_{j+3}S^{\nu^3}\cong S^{\nu^3}$ and $e_{j+3}D^{\tf_{j+3}\tf_j^3\mu(b)}\cong D^{\tf_j^3\mu(b)}$, and so $D^{\tf_j^3\mu(b)}$ is in the head of $S^{\nu^3}$. By the exactness of $f$'s, we get a surjection 
 \begin{align*}
S^{\nu(b+1)}\cong  f_0f_1^{(2)}f_2\ldots f_{p-3} S^{\nu^3}\onto 
 f_0f_1^{(2)}f_2\ldots f_{p-3} D^{\tf_j^3\mu(b)}\cong D^{\mu(b+1)}.
 \end{align*}

Assume now that $p=3$. Applying $f_j^{{(4})}$ to $\phi_b$, we get a surjection 
$f_j^{(4)}S^{\nu(b)}\to f_j^{(4)} D^{\mu(b)}.$ Now $f_j^{(4)}D^{\mu(b)}$ has simple head $D^{\tf_j^4\mu(b)}$, where $\tf_j^4\mu(b)=(b+s+t+1,b+s+1,(b+s),b+1)$. Now the result is deduces similarly to the previous case.
\end{proof}

\begin{Lemma}
\label{lem:head4}
Let $\nu:=(a+t,a,b+1,1^{a(p-1)-2})$ and $\mu:=(a+t,a^{p-1},a-1,b)$ with $a>b+1$, $t\equiv p-1\pmod{p}$ and  $a-b\equiv p-2\pmod{p}$. Then
$[S^\nu:D^\mu]=1$ and $\Hom_{\Si_{|\nu|}}(S^\nu,D^\mu)\cong \k$. 
\end{Lemma}

\begin{proof}
For the equality $[S^\nu:D^\mu]=1$, it suffices to note that $\mu=\nu^\reg$. 

Assume first that $p>3$. Consider the partitions $\xi=(a+t,a,b+2,1^{a(p-1)-2})$ and $\tau:=(a+t,a^{p-1},a-1,b+1)$ and let $j:=\res(3,b+2)$. Note that $e_jS^\xi\cong S^\nu$ and $e_jD^\tau\cong D^\mu$. By Lemma~\ref{lem:head3} we have $\Hom_{\Si_{|\xi|}}(S^\xi,D^\tau)\cong \k$. Applying $e_j$ to the surjection $S^\xi\onto D^\tau$ we get $\Hom_{\Si_{|\nu|}}(S^\nu,D^\mu)\cong \k$.

Assume now that $p=3$. Then $\nu:=(a+t,a,a-x+1,1^{2(a-1)})$ and $\mu:=(a+t,a^2,a-1,a-x)$ with $x\leq a$ and $x=a-b\equiv 1\pmod{3}$. In this case the lemma follows similarly to Lemma \ref{lem:head2} (the case $a=x$ holding by Lemma \ref{lem:head1}).
\end{proof}

\begin{Proposition}
\label{lem:hp23}
If $\la$ is as in the case (3), then there exists $u>0$ and $\mu\in \Parreg_{n-u}$, with $h(\mu)\leq p+2$ such that
$\Ext^1_{\Si_n}(D^\lambda,D^\lambda)\into \Ext^1_{\Si_{n-u}}(D^\mu,D^\mu).$
\end{Proposition}

\begin{proof}
If $b=a-1$ then $p=3$, and $h(\laM)\leq 4$. 
So we may assume that $a>b+1$. We have  $i=\res(1,a+t+1)=\res(2,a+1)$. Then $e_i^{(2)}D^\la\cong D^\mu$ with $\mu:=(a+t,a^{p-1},a-1,b)$. From the short exact sequence $0\to \rad f_i^{(2)}D^\mu \to f_i^{(2)}D^\mu\to D^\la\to 0$, using Lemma~\ref{lem:shap}, we get an exact sequence 
$$
\Hom_{\Si_n}(\rad (f_i^{(2)}D^\mu), D^\la)\to\Ext^1_{\Si_n}(D^\la,D^\la)\to
\Ext^1_{\Si_n}(f_i^{(2)}D^\mu,D^\la)\cong \Ext^1_{\Si_{n-2}}(D^\mu,D^\mu).
$$
Hence it is enough to prove that $\Hom_{\Si_n}(\rad (f_i^{(2)}D^\mu), D^\la)=0$. By Lemma~\ref{lem:head4} we have $\Hom_{\Si_{n-2}}(S^\nu,D^\mu)\cong \k$ for $\nu:=(a+t,a,b+1,1^{a(p-1)-2})$. Note that since $a>b+1$ the node $(2,a)$ is removable. Now, argue as in the proofs of Propositions~\ref{lem:hp12} and~\ref{lem:hp22}.  
\end{proof}

\begin{Proposition}
 \label{lem:hp24}
 If $\la$ is as in the case (4), then there exists $u>0$ and $\mu\in \Parreg_{n-u}$, with $h(\mu)\leq p+2$ and 
$\Ext^1_{\Si_n}(D^\lambda,D^\lambda)\into \Ext^1_{\Si_{n-u}}(D^\mu,D^\mu).$
 \end{Proposition}

\begin{proof}
Assume first that $t>p-2$.  
We have $i=\res(1,a+s+t+1)=\res(3,a+1)$, so $\res(2,a+s+1)=i+1$. Note that the $(i+2)$-addable node $(2,a+s+2)$ of $\la$ is cogood,  so $\phi_{i+2}(\la)>0$. 
We consider the sequence of partitions $\la:=\la^0,\la^1,\ldots,\la^{p-2}$, with $\la^{l+1}:=\tf_{i+2+l}^{(\phi_{i+2+l}(\la^l))}\la^l$ for  $l=0,\dots,p-3$. Then $\la^{p-2}\in\Parreg_{n+m}$ for some $m$. 
Note for all $0\leq l<p-2$ we have that $\la^{l+1}$ is an $(i+2+l)$-reflection of $\la^l$, so $\Ext^1_{\Si_n}(D^\la, D^\la)\cong \Ext^1_{\Si_{n+m}}(D^{\la^{p-2}}, D^{\la^{p-2}})$ by repeated application of Lemma~\ref{L200509}.
Note also that $\la^{p-2}_1=a+s+t+1>a+s+p-1=\la^{p-2}_2$.  
By Corollary~\ref{C200509} we have an embedding 
$$\Ext^1_{\Si_{n+m}}(D^{\la^{p-2}}, D^{\la^{p-2}})\into\Ext^1_{\Si_{n+m-1}}(D^{\te_i\la^{p-2}}, D^{\te_i\la^{p-2}}).$$
Consider the partitions  $\mu^0,\ldots,\mu^{p-2}:=\te_i\la^{p-2}$ with $\mu^{l-1}:=\te_{i+1+l}^{(\eps_{i+1+l}(\mu^{l}))}\mu^{l}$ for $l=1,\dots,p-2$.
Note that $\mu^0\in\Parreg_{n-p+1}$, and 
for all $1\leq l\leq p-2$ we have $\mu^{l-1}$ is an $(i+1+l)$-reflection of $\mu^{l}$. Hence 
$$\Ext^1_{\Si_{n+m-1}}(D^{\te_i\la^{p-2}}, D^{\te_i\la^{p-2}})\cong\Ext^1_{n-p+1}(D^{\mu^0},D^{\mu^0}).$$
by Lemma~\ref{L200509}.

If $t=p-2$ the proof is similar, using the following steps. If $a\geq p-1$ first add normal nodes of residues $i+2,i+3,\ldots,i-1$ (in this order). Then those of residue $i+1$ and $i$. Then remove conormal nodes of residues $i+1$, $i$, $i+1$ and then $i-1,i-2,\ldots,i+2$. If instead $a\leq p-2$ then first add normal nodes of residues $i+2,i+3,\ldots,i-1$ and then $i$. Then remove conormal nodes of residues $i-1,i-2,\ldots,i+2$ and then $i$.

\end{proof}

\begin{Proposition}
 If $\la$ is as in the case (5), then there exists $u>0$ and $\mu\in \Parreg_{n-u}$, with $h(\mu)\leq p+2$ and 
$\Ext^1_{\Si_n}(D^\lambda,D^\lambda)\into \Ext^1_{\Si_{n-u}}(D^\mu,D^\mu).$
 \end{Proposition}

\begin{proof}
The proof is similar to that of the previous proposition, using the same sequence of residues of nodes to be added/removed as for case (4) with $t=p-2$.
\end{proof}

\section{Proof of Theorem~\ref{TSmallWeight}}
\label{SSmallWeight}

Throughout the section we assume that $p>2$. Recall the notation $\wt(\la)$, $\quot(\la)$, $\wt_j(\Ga)$, $r_j(\Ga)$, $i_j(\Ga)$ and $\Ga_j$ from \S\ref{subsec:comb} and \S\ref{SSRA}.

\subsection{Main Tricks}
In this subsection we reduce the proof of Theorem~\ref{TSmallWeight} to some purely combinatorial facts  about partitions. Recall the notions of reflection from \S\ref{subsec:tran.fun} and of $i$-difficult from Definition~\ref{DDiff}.

\begin{Lemma}\label{T110320}
Let $w\in\Z_{\geq 0}$, and assume that for each $p$-regular partitions $\la$ with $\wt(\la)\leq w$ there exists a sequence $\la^0=\la,\la^1,\ldots,\la^k$ of $p$-regular partitions such that  $\la^\ell$ is a reflection of $\la^{\ell-1}$ for  $\ell=1,\dots k$, and  $\la^k$ satisfies {\em any} of the following conditions:
\begin{enumerate}
\item[{\rm (i)}] $|\la^k|<|\la|$;
\item[{\rm (ii)}] for $\mu\in\{\la^k,(\la^k)^\Mull\}$ there exists  $i\in I$ such that $\eps_i(\mu),\phi_i(\mu)>0$ and $\mu$ is not $i$-difficult;
\item[{\rm (iii)}] $\Ext^1_{\Si_{|\la^k|}}(D^{\la^k},D^{\la^k})=0$.
\end{enumerate}
Then $\Ext^1_{\Si_{|\la|}}(D^\la,D^\la)=0$ for all $p$-regular partitions $\la$ with $\wt(\la)\leq w$.
\end{Lemma}

\begin{proof}
By Lemma~\ref{L200509}, we have $\Ext^1_{\Si_n}(D^\la,D^\la)\cong\Ext^1_{\Si_{|\la^k|}}(D^{\la^k},D^{\la^k})$. If we are in case (ii) then by (\ref{ESignExt}) and Corollary \ref{C200509}, we have 
$$\dim\Ext^1_{\Si_n}(D^\la,D^\la)=\dim\Ext^1_{\Si_m}(D^\mu,D^\mu)\leq\dim\Ext^1_{\Si_r}(D^\nu,D^\nu)$$
where $\nu:=\tilde e_i^{\eps_i(\mu)}\mu$. 

Applying Lemma \ref{Lwt} multiple times with $\{\la,\mu\}=\{\la^\ell,\la^{\ell+1}\}$ we deduce $\wt(\la)=\wt(\la^k)$. If we are in case (ii) we also have that $\wt(\nu)<\wt(\mu)$ by the same lemma applied to $\la=\nu$ and $\mu=\mu$ and as the Mullineux bijection preserves the weight, so $\wt(\nu)<\wt(\la)$. The lemma now follows by induction on $n$ and $w$.
\end{proof}

\begin{Lemma}\label{T110320_2}
If $\Ext^1_{\Si_{|\la|}}(D^\la,D^\la)\not=0$ for some $p$-regular partition $\la$, then there exists a $p$-regular partition $\mu$ such that $\Ext^1_{\Si_{|\mu|}}(D^\mu,D^\mu)\not=0$ and the following conditions hold: 
\begin{enumerate}
\item[{\rm (1)}] $\wt(\mu)\leq\wt(\la)$ and $|\mu|\leq |\la|$;
\item[{\rm (2)}] for any $i\in I$, either $\eps_i(\mu)=0$ or 
$\mu$ is $i$-difficult. 
\end{enumerate}
\end{Lemma}

\begin{proof}
If the condition (2) holds for $\la$ we can take $\mu=\la$. So 
we may assume that for some $i\in I$, we have $\eps_i(\la)>0$  and $\la$ is not $i$-difficult. 
Let $\mu:=\te_i^{\eps_i(\la)}\la$. 
Applying Lemma \ref{Lwt} with $\la=\mu$ and $\mu=\la$, we get 
$\wt(\mu)\leq\wt(\la)$. Moreover,   
$\Ext^1_{\Si_{|\la|}}(D^\la,D^\la)\into\Ext^1_{\Si_{|\mu|}}(D^\mu,D^\mu)$ 
by Lemma \ref{L200509} and Corollary \ref{C200509}. 
Repeating the above argument if needed, we obtain $\mu$ as wanted.
\end{proof}


\begin{Lemma}\label{L110320_3}
Let $\la\in\Parreg_n$, $\Gamma=\Ga(\la)$ with the corresponding   $\quot(\la)=(\la^{(0)},\dots,\la^{(p-1)})$, 
$1\leq j\leq p-1$, $i:=i_j(\Ga)$. Suppose 
$r_{j-1}:=r_{j-1}(\Gamma)<r_j(\Gamma)=:r_j$ and $\wt_{j-1}(\Gamma)+\wt_j(\Gamma)\leq 7$. If $\la$ is $i$-difficult then 
 $(\la^{(j-1)},\la^{(j)})$ is as in Table I (the cases are labeled as $B_{\wt_{j-1}(\Gamma)+\wt_j(\Gamma)}^k$ for $k=1,2,\dots$):
\end{Lemma}

{\small
\[\begin{array}{|c|c|c|c||c|c|c|c||c|c|c|c|}
\hline
&\hspace{-3pt}\la^{(j-1)}\hspace{-4pt}&\la^{(j)}&\hspace{-3pt}r_j-r_{j-1}\hspace{-4pt}&&\hspace{-3pt}\la^{(j-1)}\hspace{-4pt}&\la^{(j)}&\hspace{-3pt}r_j-r_{j-1}\hspace{-4pt}&&\hspace{-3pt}\la^{(j-1)}\hspace{-4pt}&\la^{(j)}&\hspace{-3pt}r_j-r_{j-1}\hspace{-4pt}\\
\hline\hline
B_2^1&\varnothing&(1^2)&1&B_6^7&(1)&(2,1^3)&2&B_7^{12}&(2)&(2,1^3)&2\\
\hline
B_3^1&\varnothing&(1^3)&2&B_6^8&\varnothing&(4,1^2)&2&B_7^{13}&(1^2)&(1^5)&2\\
\hline
B_3^2&\varnothing&(2,1)&1&B_6^9&\varnothing&(3,2,1)&2&B_7^{14}&(1)&(3,1^3)&2\\
\hline
B_4^1&\varnothing&(1^4)&3&B_6^{10}&(3)&(1^3)&1&B_7^{15}&(1)&(2^3)&2\\
\hline
B_4^2&\varnothing&(2,1^2)&2&B_6^{11}&(2)&(2,1^2)&1&B_7^{16}&(1)&(2^2,1^2)&2\\
\hline
B_4^3&(1)&(1^3)&1&B_6^{12}&(1^2)&(2^2)&1&B_7^{17}&\varnothing&(5,1^2)&2\\
\hline
B_4^4&\varnothing&(3,1)&1&B_6^{13}&(1^2)&(1^4)&1&B_7^{18}&\varnothing&(4,2,1)&2\\
\hline
B_5^1&\varnothing&(1^5)&4&B_6^{14}&(1)&(3,2)&1&B_7^{19}&\varnothing&(3^2,1)&2\\
\hline
B_5^2&\varnothing&(2,1^3)&3&B_6^{15}&(1)&(3,1^2)&1&B_7^{20}&(4)&(1^3)&1\\
\hline
B_5^3&(1)&(1^4)&2&B_6^{16}&\varnothing&(5,1)&1&B_7^{21}&(3)&(2,1^2)&1\\
\hline
B_5^4&\varnothing&(3,1^2)&2&B_6^{17}&\varnothing&(2^3)&1&B_7^{22}&(2,1)&(1^4)&1\\
\hline
B_5^5&\varnothing&(2^2,1)&2&B_7^1&\varnothing&(1^7)&6&B_7^{23}&(1^3)&(2^2)&1\\
\hline
B_5^6&(2)&(1^3)&1&B_7^2&\varnothing&(2,1^5)&5&B_7^{24}&(2)&(3,1^2)&1\\
\hline
B_5^7&(1)&(2^2)&1&B_7^3&(1)&(1^6)&4&B_7^{25}&(2)&(2^2,1)&1\\
\hline
B_5^8&(1)&(2,1^2)&1&B_7^4&\varnothing&(3,1^4)&4&B_7^{26}&(1^2)&(3,2)&1\\
\hline
B_5^9&\varnothing&(4,1)&1&B_7^5&\varnothing&(2^2,1^3)&4&B_7^{27}&(1^2)&(2^2,1)&1\\
\hline
B_6^1&\varnothing&(1^6)&5&B_7^6&(2)&(1^5)&3&B_7^{28}&(1^2)&(2,1^3)&1\\
\hline
B_6^2&\varnothing&(2,1^4)&4&B_7^7&(1)&(2,1^4)&3&B_7^{29}&(1)&(4,2)&1\\
\hline
B_6^3&(1)&(1^5)&3&B_7^8&\varnothing&(4,1^3)&3&B_7^{30}&(1)&(4,1^2)&1\\
\hline
B_6^4&\varnothing&(3,1^3)&3&B_7^9&\varnothing&(3,2,1^2)&3&B_7^{31}&\varnothing&(6,1)&1\\
\hline
B_6^5&\varnothing&(2^2,1^2)&3&B_7^{10}&\varnothing&(2^3,1)&3&B_7^{32}&\varnothing&(3,2^2)&1\\
\hline
B_6^6&(2)&(1^4)&2&B_7^{11}&(3)&(1^4)&2&B_7^{33}&\varnothing&(2^3,1)&1\\
\hline
\end{array}\]

\nopagebreak
\centerline{\sc Table I}
}

\begin{proof}
By Lemma \ref{L110320}, we may assume that $r_j(\Ga)-r_{j-1}(\Ga)<\wt_{j-1}(\Ga)+\wt_j(\Ga)$. The lemma now follows from Lemma~\ref{L110320_2} by checking all possible cases for $(\Gamma_{j-1},\Gamma_j)$.
\end{proof}


\begin{Lemma}\label{L120320}
Let $\la\in\Parreg_n$ with $\wt(\la)\leq 7$, $\Gamma=\Ga(\la)$, $2\leq j\leq p-1$, $i:=i_j(\la)$. Suppose that $r_{j-2}(\Gamma)<r_{j-1}(\Gamma)<r_j(\Gamma)$. 
If $\la$ is $i$-difficult and $(i-1)$-difficult then $(\la^{(j-2)},\la^{(j-1)},\la^{(j)})$ is as in Table II:
\end{Lemma}

{\small
\[\begin{array}{|c|c|c|c|c|c|}
\hline
&\la^{(j-2)}&\la^{(j-1)}&\la^{(j)}&r_{j-1}-r_{j-2}&r_j-r_{j-2}\\
\hline\hline
C_1&\varnothing&(1^2)&(1^4)&1&2\\
\hline
C_2&\varnothing&(1^2)&(1^5)&1&3\\
\hline
C_3&\varnothing&(1^2)&(2,1^3)&1&2\\
\hline
C_4&\varnothing&(2,1)&(1^4)&1&2\\
\hline
\end{array}\]

\vspace{1mm}
\centerline{\sc Table II}
}

\begin{proof}
By Lemma \ref{L110320_3}, $((\Gamma_j,\Gamma_{j+1}),(\Gamma_{j+1},\Gamma_{j+2}))$ is one of the following:
$(B_2^1,B_6^{12})$, $(B_2^1,B_6^{13})$, $(B_2^1,B_7^{13})$, $(B_2^1,B_7^{26})$, $(B_2^1,B_7^{27})$, $(B_2^1,B_7^{28})$, $(B_3^1,B_7^{23})$, $(B_3^2,B_7^{22})$.
The lemma now follows from Lemma~\ref{L110320_2}.
\end{proof}

Let $\la\in\Parreg_n$. We will use the following terminology:

\begin{enumerate}
\item[{\bf Trick 1.}] We say that {\em Trick 1 applies to $\la$} if $\eps_{i}(\la)>0$ and $\la$ is not $i$-difficult for some $i\in I$. 

\item[{\bf Trick 2.}]  We say that {\em Trick 2 applies to $\la$} if 
there exist $i\in I$ and $m\in \Z_{\geq 2}$ such that, setting $\la^1=\la$, $\la^{\ell+1}=\tilde f_{i+\ell}^{\phi_{i+\ell}(\la^\ell)}\la^\ell$ for $\ell=1,\dots,m-1$, we have $\eps_{i+\ell}(\la^\ell)=0$ for all $\ell=1,\dots,m-1$, $\eps_{i+m}(\la^m),\phi_{i+m}(\la^m)>0$, and $\la^m$ is not $(i+m)$-difficult. 
\end{enumerate}


Suppose now that $w:=\wt(\la)\leq 7$. 
In view of Lemma~\ref{T110320}, to prove Theorem \ref{TSmallWeight}, it suffices to show that either $\min\{h(\la),h(\la^\Mull)\}\leq p+2$ or that Trick 1 or Trick 2 applies to $\la$ or $\la^\Mull$. Indeed, if Trick 1 applies to $\mu\in\{\la,\la^\Mull\}$ then either $\phi_i(\mu)=0$, in which case the condition (i)  of 
Lemma~\ref{T110320} with $k=1$ is satisfied, or $\phi_i(\la)>0$, in which case the condition (ii) of Lemma~\ref{T110320} with $k=0$ is satisfied. 
If Trick 2 applies to $\mu\in\{\la,\la^\Mull\}$ then the condition (ii) of 
 Lemma~\ref{T110320} with $k=m-1$ is satisfied. Finally, if $\min\{h(\la),h(\la^\Mull)\}\leq p+2$, then by (\ref{ESignExt}) and Theorem \ref{TSmallHeight} we have $\Ext^1_{\Si_n}(D^\la,D^\la)=0$, so in this case the condition (iii) of Lemma~\ref{T110320} with $k=0$ is satisfied.

Choose an abacus configuration $\Gamma=\Ga(\la)$ with $r_0(\Gamma)\leq r_j(\Gamma)$ for all $j\neq 0$. Let 
$\quot(\la)=(\la^{(0)},\dots,\la^{(p-1)})$ be the corresponding quotient. Let $r_j:=r_j(\Gamma)$, $\wt_j:=\wt_j(\Gamma)$ and $i_j:=i_j(\Ga)$. 
Given a fixed numer $r$ (usually $r=r_0$), we will write
\begin{equation}\label{E110221}
\Ga=((\la^{(0)},r_0-r),\ldots,(\la^{(p-1)},r_{p-1}-r)).
\end{equation}
This defines $\Ga$ up to adding/removing full rows of beads at the top, which does not affect $\quot(\la)$. Since we will (usually) work with general configurations, where some consecutive runners might be repeated and certain runners are not always present, we will use exponents to indicate the multiplicities of runners. For example
\[((\varnothing,0)^2,((1),1),(\varnothing,2),((1),2)^0,((1),0))=((\varnothing,0),(\varnothing,0),((1),1),(\varnothing,2),((1),0))\]
indicates the following abacus configuration:

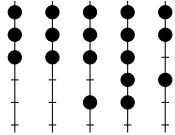
\vspace{1mm}
\hspace*{\fill}
{\begin{tikzpicture}

\draw (0,0)--(0,-1.75);
\draw (0.5,0)--(0.5,-1.75);
\draw(1,0)--(1,-1.75);
\draw (1.5,0)--(1.5,-1.75);
\draw (2,0)--(2,-1.75);

\draw(-.05,-.15)--(.05,-.15);
\draw(-.05,-.45)--(.05,-.45);
\draw(-.05,-.75)--(.05,-.75);
\draw(-.05,-1.05)--(.05,-1.05);
\draw(-.05,-1.35)--(.05,-1.35);
\draw(-.05,-1.65)--(.05,-1.65);

\draw(0.45,-.15)--(.55,-.15);
\draw(0.45,-.45)--(.55,-.45);
\draw(0.45,-.75)--(.55,-.75);
\draw(0.45,-1.05)--(.55,-1.05);
\draw(0.45,-1.35)--(.55,-1.35);
\draw(0.45,-1.65)--(.55,-1.65);

\draw(0.95,-.15)--(1.05,-.15);
\draw(0.95,-.45)--(1.05,-.45);
\draw(0.95,-.75)--(1.05,-.75);
\draw(0.95,-1.05)--(1.05,-1.05);
\draw(0.95,-1.35)--(1.05,-1.35);
\draw(0.95,-1.65)--(1.05,-1.65);

\draw(1.45,-.15)--(1.55,-.15);
\draw(1.45,-.45)--(1.55,-.45);
\draw(1.45,-.75)--(1.55,-.75);
\draw(1.45,-1.05)--(1.55,-1.05);
\draw(1.45,-1.35)--(1.55,-1.35);
\draw(1.45,-1.65)--(1.55,-1.65);

\draw(1.95,-.15)--(2.05,-.15);
\draw(1.95,-.45)--(2.05,-.45);
\draw(1.95,-.75)--(2.05,-.75);
\draw(1.95,-1.05)--(2.05,-1.05);
\draw(1.95,-1.35)--(2.05,-1.35);
\draw(1.95,-1.65)--(2.05,-1.65);

\filldraw [black] (0,-.15) circle (2.5pt);
\filldraw [black] (0,-.45) circle (2.5pt);
\filldraw [black] (0,-.75)circle (2.5pt);

\filldraw [black] (0.5,-.15) circle (2.5pt);
\filldraw [black] (0.5,-.45) circle (2.5pt);
\filldraw [black] (0.5,-.75)circle (2.5pt);

\filldraw [black] (1,-.15) circle (2.5pt);
\filldraw [black] (1,-.45) circle (2.5pt);
\filldraw [black] (1,-.75)circle (2.5pt);
\filldraw [black] (1,-1.35)circle (2.5pt);

\filldraw [black] (1.5,-.15) circle (2.5pt);
\filldraw [black] (1.5,-.45) circle (2.5pt);
\filldraw [black] (1.5,-.75)circle (2.5pt);
\filldraw [black] (1.5,-1.05)circle (2.5pt);
\filldraw [black] (1.5,-1.35)circle (2.5pt);

\filldraw [black] (2,-.15) circle (2.5pt);
\filldraw [black] (2,-.45) circle (2.5pt);
\filldraw [black] (2,-1.05)circle (2.5pt);

\end{tikzpicture}}
\hspace*{\fill}
\vspace{1mm}

If $1\leq j\leq p-1$ and $r_{j-1}<r_j$, we say that $j$ is an {\em increase (for $\Ga$)}. In view of Lemmas \ref{T110320_2} and \ref{L110320_3} if $j$ is an increase we may assume that $(\Gamma_{j-1},\Gamma_j)$ is one of the pairs $B_{\overline{w}}^k$ from Table I for some $\overline{w}\leq w$. This will be used without further reference.

Assume that $\eps_{i_j}(\la)>0$ for some runner $j$ and let $ap+j$ be the removable position corresponding to the $i_j$-good node of $\la$. Adding full rows of beads on top of $\Ga$ if necessary we may assume that $ap+j>p$. By Lemma~\ref{L110320_2}, Trick 1 does not apply to runner $j$ if and only if $\phi_{i_j}(\la)>0$, $(a-1)p+j-1$ is the addable position corresponding to the $i_j$-cogood node of $\la$, and the positions $x$ satisfying $(a-1)p+j+1\leq x\leq ap+j-2$ are occupied in $\Ga$. For $j\geq 1$ this last condition corresponds to position $ap+k$ (resp. $(a-1)p+k$) on runner $k$ being occupied if $k\leq j-2$ (resp. $k\geq j+1$). If on the other hand $j=0$ then the last condition corresponds to position $(a-1)p+k$ on runner $k$ being occupied for $1\leq k\leq p-2$. This fact will be used in the following whenever using Trick 1 to restrict abacus configurations that have to be considered.

Recall that by assumption $w\leq 7$. This will be used throuout the proof for example to exclude some cases by $p$-regularity or when applying Tricks 1 and 2 (if we know how some runners look like and that the total weight on these runners is $w'$ then there can be at most weight $7-w'$ on the remaining runners).


\subsection{No increases}
This is the case where all $r_j$ are equal, so we will shorten notation \eqref{E110221} 
to $(\la^{(0)},\ldots,\la^{(p-1)})$.

{\sf Case 1:} $\max\{h(\la^{(k)})\}\leq 1$ or $\max\{\la^{(k)}_1\}\leq 1$. In the first case $h(\la)\leq p$. In the second case if $h=\max\{h(\la^{(k)})\}$ and $j$ is minimal with $h(\la^{(j)})=h$ then the rim of $\la$ (and so also the $p$-rim) has at most $p(h+1)-j-1$ nodes and $h(\la)=ph-j$. So $h(\la^\Mull)\leq p$ provided $\la$ is $p$-regular.

{\sf Case 2:} $\max\{h(\la^{(k)})\}\geq 2$ and $\max\{\la^{(k)}_1\}\geq 2$. Assume that there is a runner $j$ with $\la^{(j)}=\varnothing$ and let $j$ be maximal such. By Trick 1 applied to runner $j+1$ if $j<p-1$ or by Trick 2 applied to the last runner if $j=p-1$ we may assume that $\la^{(k)}\not=\varnothing$ if $k<j$. In particular we may assume that there is at most one runner with 0 weight.

Let $W$ be the multiset $\{\wt_j\}$. Then we may assume that $p=3$ and $W$ is one of
\begin{align*}
&\{6,1,0\},\,\,\{5,2,0\},\,\,\{5,1^2\},\,\,\{5,1,0\},\,\,\{4,3,0\},\,\,\{4,2,1\},\,\,\{4,2,0\},\,\,\{4,1^2\},\,\,\{4,1,0\},\\
&\{3^2,1\},\,\,\{3^2,0\},\,\,\{3,2^2\},\,\,\{3,2,1\},\,\,\{3,2,0\},\,\,\{3,1^2\},\,\,\{3,1,0\},\,\,\{2^3\},\,\,\{2^2,1\},\,\,\{2^2,0\}
\end{align*}
or $p=5$ and $W$ is one of
\begin{align*}
&\{4,1^3,0\},\,\,\{3,2,1^2,0\},\,\,\{3,1^4\},\,\,\{3,1^3,0\},\,\,\{2^3,1,0\},\,\,\{2^2,1^3\},\,\,\{2^2,1^2,0\}.
\end{align*}
Using $p$-regularity and excluding cases where Trick 1 can be applied to one of the runners with non-zero weight or Trick 2 to the last runner we may assume that $\Ga$ is one of the following:
\begin{align*}
&((1),\varnothing,(4,2)),\,\,((3^2),(1),\varnothing),\,\,(\varnothing,(1),(2^3)),\,\,((1^2),\varnothing,(3,2)),\,\,((1),\varnothing,(3,2)),\\
&((2,1),(2^2),\varnothing),\,\,(\varnothing,(2^2),(2,1)),\,\,(\varnothing,(2^2),(1^3)),\,\,(\varnothing,(2^2),(1^2)),\,\,((1^2),\varnothing,(2^2)),\\
&((2^2),\varnothing,(1)),\,\,((1),\varnothing,(2^2)),\,\,((3),(2,1),(1)),\,\,((2),(1^2),(2)),\,\,((1^2),(2),(1^2)),\\
&((2),(1^2),\varnothing),\,\,((1)^3,\varnothing,(2^2)),\,\,((2,1),(1),(1^2),(1),\varnothing),\,\,((2),(1^2),(1)^2,\varnothing).
\end{align*}
If $h(\la),h(\la^\Mull)\geq p+3$ then $\Ga=((3^2),(1),\varnothing)$, in which case $\la^\Mull$ has an abacus configuration $((4,3),\varnothing^2)$ to which Trick 1 applies to runner 0.

\subsection{Increase(s) but not consecutive}\label{incr.no.cons} Recall the notation \eqref{E110221}. 
Assume that 
$$\Ga=((\varnothing,0)^{a+1},(\mu,k),(\varnothing,k)^{a_k},(\varnothing,k-1)^{a_{k-1}},\ldots,(\varnothing,0)^{a_0})$$ with $\mu\in\Par_w$, $k\geq 1$ and $a,a_1,\dots,a_k\in\Z_{\geq 0}$. Then we may assume that $a=0$ by Trick 1 (with $j=a+1$) and that $a_{k-1},\ldots,a_0=0$ by Trick 2 (with $j=p-1$ and $m=a+2$). Further if $\mu\not=(1^w)$ and the last beads on runners of the forms $(\mu,k)$ and $(\varnothing,k)$ are in rows $x$ and $y$ respectively, then $x>y+1$, so we also have $a_k=0$ and then $p=2$, giving a contradiction. Thus $\mu=(1^w)$ and then $k=w-1$ since else Trick 1 applies with $j=1$. Then $\Ga=((\varnothing,0),((1^w),w-1),(\varnothing,w-1)^{p-2})$ and so $h(\la^\Mull)=w$. In particular we may assume that $w\geq p+3$. Thus $p=3$ and $w=6$ or $7$. In either case $\la^\Mull$ has an abacus configuration $((\varnothing,0)^2,((1^w),w-1))$, to which Trick 1 applies with $j=2$.

Now assume that 
$$\Ga=((\varnothing,0)^a,((1^b),0),(\mu,k),(\varnothing,k)^{a_k},(\varnothing,k-1)^{a_{k-1}},\ldots,(\varnothing,0)^{a_0})$$ with $\mu\in\Par_{b-w}$ and $b,k\geq 1$. We may assume that Tricks 1 and 2 do not apply. Further, $\mu\not=\varnothing$. 
Then $a_{k-1},\ldots,a_0=0$. If $\mu\not=(1^{w-b})$ then  $a_k=0$. Assume first that $b\geq 2$ and $\mu=(1^{w-b})$. Then $a=0$, so $\Ga=(((1^b),0),((1^{w-b}),k),(\varnothing,k)^{p-2})$, with $k=w-2b-1$ or $w-b$. In either case $\la$ is not $p$-regular, since positions $p(r_0-1)+j$ are occupied for each $j$ but position $p(r_0-b)$ is not occupied. Assume now that $b=1$ and $\mu=(1^{w-1})$. Then $a=1$ (since we may assume that $a\leq 1$ by Trick 1 applied to the increase and if $a=0$ then $\la$ is not $p$-regular). If $\mu=(1^{w-1})$ then $\Ga=((\varnothing,0),((1),0),((1^{w-1}),k),(\varnothing,k)^{p-2})$ with $k=w-3$ or $w-1$. In the second case Trick 1 applies with $j=2$, so we may assume that $k=w-3$, in which case $\la^\Mull$ has an abacus configuration $((\varnothing,0)^{p-2},((1^{w-2}),w-3),((2),w-3))$ to which Trick 1 applies for $j=p-1$. If $\mu\not=(1^{w-b})$ then $p=3$ and $\Ga=((\varnothing,0),((1),0),(\mu,k))$, in which case $h(\la)\leq w$ (using Lemma \ref{L110320_3} to limit the possibilities for $\mu$).

So if there is an increase $j$ with $\wt_{j-1}+\wt_j=w$ we may assume that $\la^{(j-1)}\not=(1^b)$ with $b\geq 0$.

Suppose that $j<k$ are increases and $\la^{(j-1)}=\la^{(k-1)}=\varnothing$. Then we can apply Trick 1 with $i=i_k$, unless $r_{k-1}<r_{j-1}$. Since $r_{j-1},r_{k-1}\geq r_0$ by assumption, there exists an increase $\ell<j$. Repeating the above argument if needed, we may assume that $\la^{(\ell-1)}\not=0$. This will be used without further reference to reduce the number of cases that have to be studied. In particular, if there is more than one increase, we may assume by Lemma \ref{L110320_3} that $w=6$ or $7$ and that there are exactly 2 increases.

We will provided details for the first cases. The cases with only one increase $j$ for which $\wt_{j-1}+\wt_j\geq 5$ are easier, since the total weight on the other runners is at most 2. Reductions for such cases can be obtained similarly to the previous cases (using Tricks 1 and 2 as well as $p$-regularity of $\la$) and details are left to the reader. Many runner configurations can be ruled out by applying Trick 2 to the last runner. Runners with no weight can in many cases be ruled out by a combination of Trick 1 applied to the following runner and Trick 2.

We remind that we will implicitly use the fact that $w\leq 7$ in the proof. When considering abacus configurations with certain given runners this will in particular bound the total weight on other runners, thus reducing the forms that other runners might have and this will thus in some cases allow to exclude some configurations by either $p$-regularity or by applying Tricks 1 or 2.

{\sf Case $B_4^3$ and $B_2^1$}. If Trick 1 does not applies to runners corresponding to increases, then
\[\Ga=(((1),0),((1^3),1),(\varnothing,1)^a,((1),0)^b,(\varnothing,0),((1^2),1),(\varnothing,1)^c,((1),1)^d,(\varnothing,1)^e,(\varnothing,0)^f))\]
with $b+d\leq 1$. We may assume that $a,c,f=0$ by Tricks 1 or 2. 
So $p\geq 5$ and $h(\la^\Mull)=2b+5\leq 7$.
%

{\sf Case $B_5^3$ and $B_2^1$}. 
Similarly to the previous case we may first assume that
\[\Ga=(((1),0),((1^4),2),(\varnothing,2)^a,(\varnothing,1)^b,(\varnothing,0),((1^2),1),(\varnothing,1)^c,(\varnothing,0)^d)\]
and then that $b,c,d=0$. So $p\geq 5$ and $h(\la^\Mull)=6$.

{\sf Case $B_5^6$ and $B_2^1$}. 
We may first assume that
\[\Ga=(((2),0),((1^3),1),(\varnothing,1),((1^2),2),(\varnothing,2)^a),(\varnothing,1)^b)\]
and then that $b=0$. So $p\geq 5$ and $h(\la^\Mull)=7$.

{\sf Case $B_5^7$ and $B_2^1$}. No case needs to be considered (by Trick 1 applied to increases).

{\sf Case $B_5^8$ and $B_2^1$}. 
We may first assume that
\[\Ga=(((1),0),((2,1^2),1),(\varnothing,1)^a,(\varnothing,0),((1^2),1),(\varnothing,1)^b,(\varnothing,0)^c)\]
and then that $a,b,c=0$. So $p=4$ giving a contradiction.

{\sf Case $B_4^3$ and $B_3^1$}. 
We may first assume that
\[\Ga=(((1),0),((1^3),1),(\varnothing,1)^a,(\varnothing,0),((1^3),2),(\varnothing,2)^b,(\varnothing,1)^c,(\varnothing,0)^d)\]
and then that $a,c,d=0$. So $p\geq 5$ and $h(\la^\Mull)=8$. For $p=5$, it can be computed that $\Ga(\la^\Mull)=((\varnothing,0)^2,((3,1^2),2),(\varnothing,1),((1^2),2))$ to which Trick 1 applies with $j=2$.

{\sf Case $B_4^3$ and $B_3^2$}. 
We may first assume that
\[\Ga=(((1),0),((1^3),1),(\varnothing,1)^a,(\varnothing,0),((2,1),1),(\varnothing,1)^b,(\varnothing,0)^c)\]
and then that $a,b,c=0$. So $p=4$ giving a contradiction.

{\sf Case $B_2^1$}. Then $w\geq 3$. 
By Trick 1 applied to the increase we may assume that $\Ga=(A,(\varnothing,0),((1^2),1),B^1,B^0)$ where all runners of $A$ are of the form $(\mu,0)$ with $\mu$ not one of the following
\[\varnothing,\,\,(2),\,\,(3),\,\,(2,1),\,\,(4),\,\,(3,1),\,\,(2,1^2),\,\,(5),\,\,(4,1),\,\,(3,1^2),\,\,(2,1^3),\]
those of $B^1$ of the form $(\mu,1)$ with $\mu$ not one of the following
\[(1^2),\,\,(2,1),\,\,(3,1),\,\,(2^2),\,\,(4,1),\,\,(3,2)\]
and those of $B^0$ of the form $(\mu,0)$ with $\mu$ not one of the following
\[(1),\,\,(2),\,\,(3),\,\,(4),\,\,(2^2),\,\,(5),\,\,(3,2),\,\,(2^2,1).\]

By $p$-regularity of $\la$ we have that $A$ has no runner with $\mu$ one of the following
\[(1^3),\,\,(1^4),\,\,(2^2,1),\,\,(1^5),\]
$B^1$ has no runner with $\mu$ one of the following
\[(1^4),\,\,(2,1^3),\,\,(1^5),\]
$B^0$ has no runner with $\mu$ one of the following
\[(2,1^2),\,\,(1^4),\,\,(3,1^2),\,\,(2,1^3),\,\,(1^5),\]
that if $A$ has a runner $((1^2),0)$ then it also has a runner $((1),0)$ before the last such runner or $B^1$ has a runner $((1^3),1)$ or $B^0$ has a runner $((x,1),0)$ with $1\leq x\leq 2$ and that if $B^0$ has a runner $((1^3),0)$ then $B^0=(\ldots,((1^2),0),\ldots,((1^3),0),\ldots)$.

By Trick 1 (applied to the corresponding or following runners), we may assume that $A$ does not have one of the following forms:
\[(\ldots,((1),0)^2,((1^2),0),\ldots),\,\,(((1),0),((1^2),0)^2).\]
Further we may assume that $B^1=(\overline{B}^1,(\varnothing,1)^m)$ and that $\overline{B}^1$ only has runners $((1),1)$ and $((1^3),1)$, with $((1^3),1)$ the first runner of of $\overline{B}^1$ if it appears. Using also Trick 2 to rule out runners $(\varnothing,0)$, it can also be checked that $B^0$ only has possible runners $((1^2),0)$, $((2,1),0)$ or $((1^3),0)$, with $((2,1),0)$ appearing at the beginning of $B^0$ if it is present. Further $B^1$ has no runner $(\varnothing,1)$ if $B^0=(((2,1),0),\ldots)$.

So since $3\leq w\leq 7$ we may assume that
\begin{align*}
\Ga=&(((1^2),0)^a,((1),0)^b,((1^2),0)^c,((x-2,2),0)^d,((1),0)^e,(\varnothing,0),((1^2),1),((1^3),1)^f,\\
&((1),1)^g,(\varnothing,1)^h,((2,1),0)^i,((1^2),0)^k,((1^3),0)^l)
\end{align*}
with $4\leq x\leq 5$, $l\leq k$, $1\leq 2a+b+2c+dx+e+3f+g+3i+2k+3l\leq 5$ and $h\cdot i=0$. Further $a\leq c\leq b\leq 1$ or $f+i+k\geq 1$. If $h(\la),h(\la^\Mull)\geq p+3$ then $\Ga$ is of the form 
\[\Ga=(((1^2),0)^C,((1),0),((1^2),0),((1),0)^D,(\varnothing,0),((1^2),1),(\varnothing,1)^E)\]
with $(C,D,E)\in\{(1,0,0),(0,1,0),(0,2,1)\}$. In the first case we have  $$\Ga(\la^\Mull)=((\varnothing,0),((4,1),1),(\varnothing,1),((2),1),(\varnothing,1)),$$ to which Trick 1 applies with $j=3$. In the other two cases $$\Ga(\la^\Mull)=((\varnothing,0)^D,((w-2,1),1),(\varnothing,1)^{D+1},((1),1)),$$ to which Trick 1 applies with $j=p-1$.

{\sf Case $B_3^1$}. Then $w\geq 4$. 
We may assume by Trick 1 applied to the increase that $\Ga=(A,(\varnothing,0),((1^3),2),B^2,B^1,B^0)$, where all runners of $A$ are of the form $(\mu,0)$, all runners of $B^i$ of the form $(\mu,i)$ for $0\leq i\leq 2$, $A$ has no runner $(\mu,0)$ with $\mu$ one of the following:
\[\varnothing,\,\,(2),\,\,(3),\,\,(2,1),\,\,(4),\,\,(3,1),\,\,(2,1^2)\]
and $(B^2,B^1,B^0)$ has no runners of the following forms:
\begin{align*}
&((1^3),2),\,\,((2,1^2),2),\,\,((1^2),1),\,\,((2,1),1),\,\,((3,1),1),\,\,((2^2),1),\,\,((1),0),\,\,((2),0),\\
&((3),0),\,\,((4),0),\,\,((2^2),0).
\end{align*}
Since $\la\in\Parreg_n$, $A$ has no runner $((1^x),0)$ with $3\leq x\leq 4$ and $B$ has no runner of the forms:
\begin{align*}
((1^4),1),\,\,((1^3),0),\,\,((2,1^2),0),\,\,((1^4),0).
\end{align*}
From Tricks 1 and 2 we can also see that $B^2=(\overline{B}^2,(\varnothing,2)^m)$ and that $\overline{B}^2$, $B^1$ and $B^0$ have no runners of the form $(\varnothing,i)$. Further by Trick 1, $\overline{B}^2$ has no runner of the forms
\[((2),2),\,\,((3),2),\,\,((2,1),2),\,\,((4),2),\,\,((3,1),2),\,\,((2^2),2)\]
and by Trick 2, $B^1$ has no runner of the form $((x),1)$ with $1\leq x\leq 4$.

If $A=(((1^2),0),\ldots)$ and $B$ has no runner of the form $((1^2),0)$ then $\la\not\in\Parreg_n$, while if $A=(((1),0)^2,((1^2),0))$ or $\overline{B}^2=(\ldots,((1),2),((1^2),2),\ldots)$ then Trick 1 applies.

So we may assume
\begin{align*}
\Ga=&(((1),0)^a,((1^2),0)^b,((1),0)^c,((2^2),0)^d,(\varnothing,0),((1^3),2),((1^2),2)^e,((1),2)^f,\\
&((1^4),2)^g,(\varnothing,2)^h,((x-2,1^2),1)^i,((y-1,1),0)^k)
\end{align*}
with $3\leq x\leq 4$, $2\leq y\leq 4$ and $1\leq a+2b+c+4d+2e+f+4g+ix+ky\leq 4$. Further $b\leq a\leq 1$ or $(y,k)=(2,1)$. If $(x,i)$ or $(y,k)=(4,1)$, then we may further assume that $h=0$.

If $h(\la),h(\la^\Mull)\geq p+3$ then
\[\Ga=(((2^2),0),(\varnothing,0),((1^3),2))\,\,\text{or}\,\,((1^2),0),(\varnothing,0),((1^3),2),(\varnothing,2),((1^2),0))\]
or $\Ga$ is of the form
\begin{align*}
&(((1),0),((1^2),0)^b,((1),0)^c,(\varnothing,0),((1^3),2),((1^2),2)^e,((1),2)^f,(\varnothing,2)^h,((1^2),0)^k,((2,1),0)^l)
\end{align*}
with $b+c+e+f+k+3\leq p\leq 3b+3c+2k+3l+3$ and $(c,p)\not\in\{(0,3),(2,5)\}$. It can be checked that in either case $\la^\Mull$ has an abacus configuration of the form $((\varnothing,0)^m,(\mu,2),C)$, where $1\leq m\leq p-2$, $\mu_1\geq 3$, $h(\mu)=3$ and all runners of $C$ are of the form $(\varnothing,2)$, $((1),2)$, $((1^2),2)$ or $(\varnothing,1)$. Further $m\geq 2$ if $C$ has a runner $((1^2),2)$. So Trick 2 applies to the last runner or Trick 1 applies to runner $m$.

{\sf Case $B_3^2$}. Then $w\geq 4$. We may assume by Trick 1 applied to the increase that $\Ga=(A,(\varnothing,0),((2,1),1),B^1,B^0)$, where all runners of $A$ are of the form $(\mu,0)$, all runners of $B^i$ of the form $(\mu,i)$ for $0\leq i\leq 1$, no runner of $A$ is of the form $(\mu,0)$ with $\mu$ one of the following
\[\varnothing,\,\,(2),\,\,(3),\,\,(2,1),\,\,(4),\,\,(3,1),\,\,(2,1^2)\]
and $(B^1,B^0)$ has no runner of the forms
\[((1^2),1),\,\,((2,1),1),\,\,((3,1),1),\,\,((2^2),1),\,\,((1),0),\,\,((2),0),\,\,((3),0),\,\,((4),0),\,\,((2^2),0).\]

By $p$-regularity we further have that $A$ has no runner of the forms $((1^3),0)$ or $((1^4),0)$, that $B^1$ has no runner of the form $((1^4),1)$ and $B^0$ has no runner of the forms $((1^3),0)$, $((2,1^2),0)$ or $((1^4),0)$. Further if $A$ has a runner of the form $((1^2),0)$ then there is a runner of the form $((1),0)$ before it or $B^0$ has a runner of the form $((1^2),0)$.

By Trick 1 we may also exclude that $B^1$ has runners of the forms $((3),1)$ or $((4),1)$ and, also using Trick 2, that $B^1$ has a runner $(\varnothing, 1)$ or $B^0$ a runner $(\varnothing,0)$. Further $A\not=(((1),0)^2,((1^2),0))$, $B^1\not=(\ldots,((1),1),((2),1),\ldots)$ and $B^1\not=(\ldots,((1),1),((1^3),1),\ldots)$.

So we may assume
\begin{align*}
\Ga=&(((1),0)^a,((1^2),0)^b,((1),0)^c,((2^2),0)^d,(\varnothing,0),((2,1),1),((2),1)^e,((x-2,1^2),1)^f,\\
&((1),1)^g,((y-1,1),0)^h)
\end{align*}
with $3\leq x\leq 4$, $2\leq y\leq 4$ and $1\leq a+2b+c+4d+2e+fx+g+hy\leq 4$. Further $b\leq a\leq 1$ or $(b,f,y)=(1,1,2)$. It follows that $p-2=a+b+c+d+e+f+g+h=1$ or 3. If $h(\la),h(\la^\Mull)\geq p+3$ then $p=5$ and $\Ga=(((1),0),((1^2),0),((1),0),(\varnothing,0),((2,1),1))$, in which case $\la^\Mull$ has an abacus configuration $((\varnothing,0)^2,((1^2),0),(\varnothing,-2),((4,1),-1))$, to which Trick 1 applies with $j=2$.

{\sf Case $B_4^1$}. Then $w\geq 5$. Again we may assume that $\Ga=(A,(\varnothing,0),((1^4),3),B^3,B^2,B^1,B^0)$, where all runners of $A$ are of the form $(\mu,0)$ and those of $B^i$ of the form $(\mu,i)$ for $0\leq i\leq 3$. By Trick 1 applied to the increase we may assume that $A$ has no runner $(\mu,0)$ with $\mu$ one of the following
\[\varnothing,\,\,(2),\,\,(3),\,\,(2,1)\]
and that $(B^3,B^2,B^1,B^0)$ has no runner of the forms
\[((1^3),2),\,\,((1^2),1),\,\,((2,1),1),\,\,((1),0),\,\,((2),0),\,\,((3),0).\]

By $p$-regularity, neither $A$ nor $B^0$ has runners of the form $((1^3),0)$ and further that if $A$ has a runner $((1^2),0)$ then $A=(((1),0),((1^2),0))$.

By Trick 1 we may assume that $B^3=(\overline{B}^3,(\varnothing,0)^m)$ where $\overline{B}^3$ can only have runners $((1),3)$, $((1^2),3)$ or $((1^3),3)$ and that $B^2$ has no runner $((3),2)$. Further $\overline{B}^3\not=(((1),3),((1^2),3))$. Applying also Trick 2 we may further assume that $B^0$ has no runner $(\varnothing,0)$, that $B^1=(((1^3),1)^m)$ and that $B^2$ has no runner.

So we may assume
\begin{align*}
\Ga=(((1),0)^a,((1^2),0)^b,(\varnothing,0),((1^4),3),((1^x),3)^c,((1),3)^d,(\varnothing,3)^e,((1^3),1)^f,((x-1,1),0)^g)
\end{align*}
with $2\leq x\leq 3$, $1\leq a+2b+cx+d+3f+gx\leq 3$ and $b\leq a$. So $h(\la^\Mull)=4a+4b+2f+3g+4$. For $p\leq 4a+4b+2f+3g+1$, $\la$ has an abacus configuration of the form $((\varnothing,0)^m,(\mu,3),C)$, where $1\leq m\leq p-2$, $\mu\in\{(4,1^3),(3,2,1^2),(3,1^3),(2^3,1),(2^2,1^2),(2,1^3)\}$ and all runners of $C$ are of the forms $((1),3)$, $(\varnothing,3)$, $(\varnothing,2)$ or $(\varnothing,1)$. So Trick 2 applies to the last runner.

{\sf Case $B_4^2$}. Then $w\geq 5$. We may again assume that $\Ga=(A,(\varnothing,0),((2,1^2),2),B^2,B^1,B^0)$, where all runners of $A$ are of the form $(\mu,0)$ and those of $B^i$ of the form $(\mu,i)$ for $0\leq i\leq 2$. As in the previous case, by Trick 1 applied to the increase we may assume that $A$ has no runner $(\mu,0)$ with $\mu$ one of the following
\[\varnothing,\,\,(2),\,\,(3),\,\,(2,1)\]
and that $(B^2,B^1,B^0)$ has no runner of the forms
\[((1^3),2),\,\,((1^2),1),\,\,((2,1),1),\,\,((1),0),\,\,((2),0),\,\,((3),0).\]

By $p$-regularity, again neither $A$ nor $B^0$ has runners of the form $((1^3),0)$ and further that if $A$ has a runner $((1^2),0)$ then $A=(((1),0),((1^2),0))$.

By Tricks 1 and 2 we may further assume that $B^0$ has no runner $(\varnothing,0)$, that $B^1$ has no runner, that $B^2$ has no runner $(\varnothing,2)$ or $((3),2)$ and that $B^2\not=(((1),2),((2),2))$ or $(((1),2),((1^2),2))$.

So we may assume
\begin{align*}
\Ga=&(((1),0)^a,((1^2),0)^b,(\varnothing,0),((2,1^2),2),((2),2)^c,((x-1,1),2)^d,((1),2)^e,\\
&((x-1,1),0)^f)
\end{align*}
with $2\leq x\leq 3$, $1\leq a+2b+2c+dx+e+fy\leq 3$ and $b\leq a$. So $a+b+c+d+e+f=1$ or $3$. Then $b=0$ and either $a+c+d+e+f+g=1$ or $a+e=3$. It can be checked that $h(\la)$ or $h(\la^\Mull)\leq p+2$ unless $\Ga=(((1),0)^2,(\varnothing,0),((2,1^2),2),((1),2))$, in which case case $\la^\Mull$ has an abacus configuration $((\varnothing,0)^2,((3,2,1),2),((1),2),(\varnothing,2))$, to which Trick 1 applies with $j=2$.

{\sf Case $B_4^3$}. Then $w\geq 5$.  We may assume that $\Ga=(A,((1),0),((1^3),1),B^1,B^0)$, where all runners of $A$ are of the form $(\mu,0)$ and those of $B^i$ of the form $(\mu,i)$ for $0\leq i\leq 1$. By Trick 1 applied to the increase we may assume that $A$ has no runner $((1),0))$, $((2),0)$ or $((3),0)$, that $B^1$ has no runner $((1^3),1)$ and that $B^0$ has no runner $((1^2),0)$ or $((2,1),0)$.

By $p$-regularity $A$ has no runner $((1^3),0)$.

By Trick 1 $A=(\overline{A},(\varnothing,0)^m)$ where $\overline{A}$ can only have runners $((1^2),0)$ and $((2,1),0)$, $B^1=(\overline{B}^1,(\varnothing,1)^m)$ where $\overline{B}^1$ can only have runners $((1),1)$ and $((1^2),1)$, $B^0$ has no runner $((3),0)$ and $B^0\not=(\ldots,(\varnothing,0),((1^3),0),\ldots)$. By Trick 2 we then have that $B^0=(((1^3),0)^m,(\varnothing,0)^x)$.

If $B^0=(((1^3),0),(\varnothing,0)^x)$ then $A=((\varnothing,0)^r)$ with $r\geq 1$ by $p$-regularity and so Trick 2 applies to the last runner. Thus $B^0=((\varnothing,0)^x)$.

So we may assume
\[\Ga=(((x,1),0)^a,(\varnothing,0)^b,((1),0),((1^3),1),((1),1)^c,((1^2),1)^d,((1),1)^e,(\varnothing,1)^f,(\varnothing,0)^g)\]
with $1\leq x\leq 2$, $1\leq ax+c+2d+e\leq 3$ and $g\leq 1$ with $g=0$ unless $(a,x)=(1,3)$. Further $b\leq 1$ by Trick 1 applied to runner $a+b$ and then $b=1$ unless $(a,x)=(1,3)$ or $(c,d)=(1,1)$ by $p$-regularity. In each case $h(\la^\Mull)\leq p+1$.

{\sf Case $B_4^4$}. Then $w\geq 5$. We may again assume that $\Ga=(A,(\varnothing,0),((3,1),1),B^1,B^0)$, where all runners of $A$ are of the form $(\mu,0)$ and those of $B^i$ of the form $(\mu,i)$ for $0\leq i\leq 1$. By Trick 1 applied to the increase we may assume that $A$ has no runner $(\varnothing,0)$, $((2),0)$, $((3),0)$ or $((2,1),0)$, $B^1$ has no runner $((1^2),1)$ or $((2,1),1)$ and $B^0$ has no runner $((1),0)$, $((2),0)$ or $((3),0)$.

By $p$-regularity neither $A$ nor $B^0$ has a runner $((1^3),0)$ and if $A$ has a runner $((1^2),0)$ then $A=(((1),0),((1^2),0))$.

By Trick 2 $B^0$ has no runner. Further by Trick 1 (applied to the following runner) or 2 (applied to the last runner), $B^1$ has no runner $(\varnothing,1)$, $((1),1)$ or $((1^3),1)$.

So we may assume
$\Ga=(((1),0)^a,((1^2),0)^b,(\varnothing,0),((3,1),1),((x),1)^c)$ 
with $2\leq x\leq 3$, $1\leq a+2b+cx\leq 3$ and $b\leq a$. Then $p-2=a+b+c=1$ or $3$, so $b=0$, $(a,c)\in\{(1,0),(3,0),(0,1)\}$ and then $h(\la)=p+1$.

{\sf Case $B_5^1$}. Then $w\geq 6$. We may assume
\[\Ga=(((1),0)^a,(\varnothing,0),((1^5),4),((1^x),4)^b,(\varnothing,4)^c,((1^2),0)^d)\]
with $1\leq x\leq 2$ and $1\leq a+bx+2d\leq 2$. So $h(\la^\Mull)=5a+4d+5$. For $p\leq 5a+4d-2$, $\Ga(\la^\Mull)=((\varnothing,0)^{p-a-d-1},((a+d+1,d+xb+1,1^3),4),(\varnothing,4-d)^{a+d})$, to which Trick 2 applies to the last runner.

{\sf Case $B_5^2$}. Then $w\geq 6$. We may assume
\[\Ga=(((1),0)^a,(\varnothing,0),((2,1^3),3),((2),3)^b,((1^x),3)^c,((1^2),0)^d)\]
with $1\leq x\leq 2$ and $1\leq a+2b+cx+2d\leq 2$. So $a+b+c+d=1$ and $p=3$. If $a=1$ then $\la^\Mull$ has an abacus configuration $((\varnothing,0),((2^2,1^2),3),(\varnothing,3))$, to which Trick 2 applies with $j=2$. If $b+c=1$ then $h(\la^\Mull)=4$. If $d=1$ then $\la^\Mull$ has an abacus configuration $((\varnothing,0),((2^3,1),3),(\varnothing,2))$, to which Trick 2 applies to the last runner.

{\sf Case $B_5^3$}. Then $w\geq 6$. We may assume
$$\Ga=(((1^2),0)^a,(\varnothing,0)^b,((1),0),((1^4),2),((1^x),2)^c,(\varnothing,2)^d)$$
with $1\leq x\leq 2$ and $1\leq 2a+cx\leq 2$. Further $b=1$ by $p$-regularity and Trick 1 applied to runner $a+b$. So $p\geq 5$ and $h(\la^\Mull)=6$ if $a=0$ or $h(\la^\Mull)=9$ if $a=1$. If $p=5$ and $\Ga=(((1^2),0),(\varnothing,0),((1),0),((1^4),2),(\varnothing,2))$ then $\la^\Mull$ has an abacus configuration $(((2,1^2),0),((3),0),(\varnothing,0),(\varnothing,-3)^2)$, to which Trick 1 applies to runner 1.

{\sf Case $B_5^4$}. Then $w\geq 6$. We may assume
\[\Ga=(((1),0)^a,(\varnothing,0),((3,1^2),2),((2),2)^b,((1^2),0)^c)\]
with $1\leq a+2b+2c\leq 2$. So $a+b+c=1$, $p=3$ and $h(\la)=5$ or $h(\la^\Mull)\leq 5$.

{\sf Case $B_5^5$}. Then $w\geq 6$. We may assume $\Ga=(((1),0)^a,(\varnothing,0),((2^2,1),2),((1^2),2)^b)$ with $1\leq a+2b\leq 2$. So $a+b=1$, $p=3$ and $h(\la)=5$ or $h(\la^\Mull)=3$.

{\sf Case $B_5^6$}. We may assume that
\[\Ga=(((1^2),0)^a,((2),0),((1^3),1),((1^x),1)^b,(\varnothing,1)^c,((1),1)^d,(\varnothing,1)^e,((y),0)^f,(\varnothing,0)^g)\]
with $1\leq x,y\leq 2$, $c\leq d\leq 1$, $g\leq 1$ and $2a+bx+d+fy\leq 2$. In either case $h(\la^\Mull)\leq p+1$.

{\sf Case $B_5^7$}. Then $w\geq 6$. We may assume
$\Ga=(((2),0)^a,((1),0),((2^2),1),((1^2),1)^b,((1),0)^c)$
with $1\leq 2a+2b+c\leq 2$. Then $a+b+c=1$, $p=3$ and $h(\la)\leq 5$.

{\sf Case $B_5^8$}. Then $w\geq 6$. We may assume
$$\Ga=(((1^2),0)^a,(\varnothing,0)^b,((1),0),((2,1^2),1),((x),1)^c,((1^2),1)^d)$$
with $1\leq x\leq 2$ and $1\leq 2a+xc+d\leq 2$. Further we may assume that $b=1$ by $p$-regularity and Trick 1 applied to runner $a+b$. So $a=0$, $x=1$, $c=2$ and $d=1$, in which case $p=5$ and $h(\la^\Mull)=4$.

{\sf Case $B_5^9$}. Then $w\geq 6$. We may assume $\Ga=(((1),0)^a,(\varnothing,0),((4,1),1))$ with $1\leq a\leq 2$, so $a=1$, $p=3$ and $h(\la)=4$.

{\sf Case $B_6^1$}. Then $w=7$. We may assume $\Ga=(((1),0)^a,(\varnothing,0),((1^6),5),((1),5)^b,(\varnothing,5)^c)$ with $a+b=1$. So $h(\la^\Mull)=5a+6$. For $p\leq 5a+3$, $\Ga(\la^\Mull)=((\varnothing,0)^{p-1-a},((2,1^5),5),(\varnothing,5)^a)$, to which Trick 1 or 2 applies.

{\sf Case $B_6^2$}. Then $w=7$. We may assume $\Ga=(((1),0)^a,(\varnothing,0),((2,1^4),4),((1),4)^b)$ with $a+b=1$. If $b=1$ then $h(\la^\Mull)=5$. If $a=1$ then $\Ga(\la^\Mull)=((\varnothing,0),(\varnothing,-5),((2^2,1^3),-1))$, to which Trick 2 applies with $j=0$.

{\sf Case $B_6^3$}. Then $w=7$. We may assume that $\Ga=((\varnothing,0)^a,((1),0),((1^5),3),((1),3),(\varnothing,3)^b)$. By $p$-regularity and Trick 1 applied to runner $a$ we may assume that $a=1$. So $p\geq 5$ and $h(\la^\Mull)=8$. For $p=5$, $\Ga(\la^\Mull)=((\varnothing,0)^3,((1^4),3),((2,1),3))$, to which Trick 1 applies with $j=3$.

{\sf Case $B_6^4$}. Then $w=7$. We may assume that $\Ga=(((1),0),(\varnothing,0),((3,1^3),3))$. So $p=3$ and $\Ga(\la^\Mull)=((\varnothing,0),((2^2,1^2),3),((1),3))$, to which Trick 2 applies with $j=2$.

{\sf Case $B_6^5$}. Then $w=7$. We may assume $\Ga=(((1),0),(\varnothing,0),((2^2,1^2),3)$. So $p=3$ and $\Ga(\la^\Mull)=((\varnothing,0),((2^3,1),3),(\varnothing,3))$, to which Trick 2 applies with $j=2$.

{\sf Case $B_6^6$}. We may assume $\Ga=(((2),0),((1^4),2),((1),2)^a,(\varnothing,2)^b,(\varnothing,0)^c)$ with $a\leq 1$. Further we may assume that $c=1$ by $p$-regularity and Trick 1 applied to runner 0. So $h(\la^\Mull)=5$.

{\sf Case $B_6^7$}. Then $w=7$. We may assume that $\Ga=((\varnothing,0)^a,((1),0),((2,1^3),2),((1),2))$ with $a\leq 1$. So $a=0$ and then $\la\not\in\Parreg_n$, leading to a contradiction.

{\sf Case $B_6^8$}. Then $w=7$. We may assume $\Ga=(((1),0),(\varnothing,0),((4,1^2),2))$. So $p=3$ and $h(\la)=5$.

{\sf Case $B_6^9$}. Then $w=7$. We may assume $\Ga=(((1),0),(\varnothing,0),((3,2,1),2)$. So $p=3$ and $h(\la)=5$.

{\sf Case $B_6^{10}$}. We may assume that $\Ga=(((3),0),((1^3),1),((1),1)^a,(\varnothing,1)^b,((1),0)^c)$ with $a+c\leq 1$ and $b\leq 1$. So $p=3$. If $a=1$ then $\la\not\in\Parreg_n$, while if $a=0$ then $h(\la^\Mull)=3+b$.

{\sf Case $B_6^{11}$}. We may assume that $\Ga=(((2),0),((2,1^2),1),((1),1)^a,((1),0)^b,(\varnothing,0)^c)$ with $a+b\leq 1$ and $b\leq c\leq 1$. So $p=3$ and $h(\la^\Mull)\leq 5$.

{\sf Case $B_6^{12}$}. By Trick 1 applied to the increase we may assume that $$\Ga=(((1^2),0),((2^2),1),(\varnothing,1)^a),((1),0)^b),$$
in which case $\la\not\in\Parreg_n$ giving a contradiction.

{\sf Case $B_6^{13}$}. Then $w=7$. By $p$-regularity $$\Ga=((\varnothing,0)^a,((1),0),(\varnothing,0)^b,((1^2),0),((1^4),1),(\varnothing,1)^c,(\varnothing,0)^d)$$
with $a\geq 1$. We may assume that $a=1$ and $b,d=0$, in which case $p\geq 5$ and $h(\la^\Mull)=6$.

{\sf Case $B_6^{14}$}. Then $w=7$ and by Trick 1 applied to the increase we may assume $\Ga=(((1),0),((3,2),1),(\varnothing,1)^{p-3},((1),0))$. So $h(\la^\Mull)=p+2$.

{\sf Case $B_6^{15}$}. Then $w=7$. By Trick 1 applied to the increase we may assume that there is a runner with non-zero weight after the increase. In particular the increase is not at the end and then Trick 2 applies to the last runner.

{\sf Case $B_6^{16}$}. Then $w=7$. We may assume that $\Ga=(((1),0),(\varnothing,0),((5,1),1))$, so $p=3$ and $h(\la^\Mull)=4$.

{\sf Case $B_6^{17}$}. Then $w=7$. We may assume that $\Ga=(((1),0),(\varnothing,0),((2^3),1))$, so $p=3$ and $h(\la^\Mull)=5$.

{\sf Case $B_7^6$}. We may assume $\Ga=(((2),0),((1^5),3),(\varnothing,3)^a,(\varnothing,0)^b)$. Further we may assume that $b=1$ by $p$-regularity and Trick 1 applied to runner 0. So $h(\la^\Mull)=7$. For $p=3$, $\la^\Mull$ has an abacus configuration $((\varnothing,0),((1^4),3),((3),2))$, to which Trick 1 applies with $j=2$.

{\sf Case $B_7^{11}$}. We may assume $\Ga=(((3),0),((1^4),2),(\varnothing,2)^a,(\varnothing,1)^b)$. Further we may assume that $b=1$ by $p$-regularity and Trick 1 applied to runner 0. Then $h(\la^\Mull)=4$.

{\sf Case $B_7^{12}$}. We may assume $\Ga=(((2),0),((2,1^3),2),(\varnothing,0)^a)$ with $a\leq 1$. So $p=3$ and $h(\la^\Mull)=5$.

{\sf Case $B_7^{20}$}. We may assume $\Ga=(((4),0),((1^3),1))$. So $p=2$ giving a contradiction. 

{\sf Case $B_7^{21}$}. We may assume $\Ga=(((3),0),((2,1^2),1))$. So $p=2$ giving a contradiction.

{\sf Case $B_7^{22}$}. Then $\la\not\in\Parreg_n$ giving a contradiction.

{\sf Case $B_7^{24}$}. We may assume $\Ga=(((2),0),((3,1^2),1),(\varnothing,0)^a)$ with $a\leq 1$. So $p=3$ and $h(\la^\Mull)=5$.

{\sf Case $B_7^{25}$}. We may assume $\Ga=(((2),0),((2^2,1),1))$. So $p=2$ giving a contradiction.

\subsection{Consecutive increases}\label{cons.incr}
Here we assume that there exists $2\leq j\leq p-1$ such that $j-1$ and $j$ are increases. By Lemma \ref{L120320}, one of the configurations $C_k$ appears in $\Ga$ and we may assume there is no further increase in $r_j$ (by Lemmas \ref{L110320_3} and \ref{L120320}). 
The reductions can be obtained similarly to those in \S\ref{incr.no.cons} (use Tricks 1 and 2 and $p$-regularity) and details are left to the reader.

{\sf Case $C_1$}. We may assume that
$\Ga=((\varnothing,0),((1^2),1),((1^4),2),((1),2)^a,(\varnothing,2)^b)$ 
with $a\leq 1$. So $h(\la^\Mull)=5$.

{\sf Case $C_2$}. We may assume $\Ga=((\varnothing,0),((1^2),1),((1^5),3),(\varnothing,3)^a)$. So $h(\la^\Mull)=7$. For $p=3$, $\la^\Mull$ has an abacus configuration $(((1^2),0),((2^2),1),((1),-3))$, to which we can apply Trick 1 with $j=1$.

{\sf Case $C_3$}. We may assume $\Ga=((\varnothing,0),((1^2),1),((2,1^3),2))$. So $p=3$ and $h(\la^\Mull)=5$.

{\sf Case $C_4$}. We may assume $\Ga=((\varnothing,0),((2,1),1),((1^4),2),(\varnothing,2)^a)$. So $h(\la^\Mull)=5$.

\end{document}